\documentclass{amsart}
\usepackage{amssymb}
\usepackage{amsfonts}
\usepackage{latexsym}

\newtheorem{theorem}{Theorem}[section]
\newtheorem{lemma}[theorem]{Lemma}
\newtheorem{proposition}[theorem]{Proposition}
\newtheorem{corollary}[theorem]{Corollary} 
\theoremstyle{definition}  
\newtheorem{definition}[theorem]{Definition}

\newtheorem{example}[theorem]{Example}
\newtheorem{question}[theorem]{Question}
\newtheorem{conjecture}[theorem]{Conjecture}  
\newtheorem{remark}[theorem]{Remark}

\newcommand{\Tr}{\text{Tr}}
\newcommand{\id}{\text{id}}

\newcommand{\End}{\text{End}} 

\newcommand{\Hom}{\text{Hom}} 
\newcommand{\Ad}{\text{Ad}}
\newcommand{\Id}{\text{Id}}

\newcommand{\Rep}{\text{Rep}}

\newcommand{\eps}{\varepsilon}

\newcommand{\I}{\mathcal{I}}

\newcommand{\g}{\mathfrak{g}}

\newcommand{\actl}{\rightharpoonup}
\newcommand{\actr}{\leftharpoonup}

\newcommand{\la}{\langle\,} 
\newcommand{\ra}{\,\rangle}

\newcommand{\1}{_{(1)}} 
\newcommand{\2}{_{(2)}} 
\newcommand{\3}{_{(3)}} 
\newcommand{\4}{_{(4)}}

\renewcommand{\I}{^{(1)}} 
\newcommand{\II}{^{(2)}}

\newcommand{\R}{\mathcal{R}}

\newcommand{\ben}{\begin{enumerate}}
\newcommand{\een}{\end{enumerate}}

\newcommand{\mC}{{\mathcal C}}
\newcommand{\mB}{{\mathcal B}}
\newcommand{\mD}{{\mathcal D}}
\newcommand{\trace}{{\rm tr}}
\newcommand{\CC}{{\mathbb{C}}}
\newcommand{\mM}{{\mathcal M}}
\newcommand{\mN}{{\mathcal N}}

\begin{document}

\title{On fusion categories}
\begin{abstract}Using a variety of methods developed in the literature 
(in particular, the theory of weak Hopf algebras),
we prove a number of general results about fusion categories
in characteristic zero. 
We show that the global dimension 
of a fusion category is always positive, and that the S-matrix 
of any (not necessarily hermitian) modular category is unitary. 
We also show that the category of module functors between two module categories 
over a fusion category is semisimple, and that fusion categories
and tensor functors between them are undeformable (generalized
Ocneanu rigidity). In particular the number of such categories
(functors) realizing a given fusion datum is finite. 
Finally, we develop the theory of Frobenius-Perron dimensions in an 
arbitrary fusion category.
At the end of the paper we generalize some of these results to
positive characteristic.  
\end{abstract}

\author{Pavel Etingof}
\address{Department of Mathematics, Massachusetts Institute of Technology,
Cambridge, MA 02139, USA}
\email{etingof@math.mit.edu}

\author{Dmitri Nikshych}
\address{Department of Mathematics and Statistics,
University of New Hampshire,  Durham, NH 03824, USA}
\email{nikshych@math.unh.edu}

\author{Viktor Ostrik}
\address{Department of Mathematics, Massachusetts Institute of Technology,
Cambridge, MA 02139, USA}
\email{ostrik@math.mit.edu}

\maketitle


\section{Introduction}

Throughout this paper (except for Section 9), 
$k$ denotes an algebraically closed field of characteristic zero. 
By a {\em fusion category} $\mC$ over $k$ we mean 
a $k$-linear semisimple rigid tensor (=monoidal) 
category with finitely many simple 
objects and finite dimensional spaces of morphisms, 
such that the endomorphism algebra of the neutral object is $k$
(see \cite{BaKi}). 
Fusion categories arise in several areas of mathematics and physics -- 
conformal field theory, operator algebras, representation theory of 
quantum groups, and others. 

This paper is devoted to the study of general properties
of fusion categories. This has been an area
of intensive research for a number of years, 
and many remarkable results have been 
obtained. However, many of these results were proved not in general but 
under various assumptions on the category.   
The goal of this paper is to remove such assumptions, and 
to give an account of the theory of fusion categories in full generality. 

The structure of the paper is as follows. 

In Section 2, we 
give our main results about squared norms, global dimensions,
weak Hopf algebras, and Ocneanu rigidity. 
This section contains many results which are partially 
or fully due to other authors, and our own results are blended in at 
appropriate places. 

Sections 3-7 are mostly devoted to review of the 
technical tools and to proofs of the results of Section 2.
Namely, in section 3, 
we prove one of the main theorems of this paper, saying that 
any fusion category has a nonzero global dimension (in fact, we show that 
for $k=\CC$, the dimension is positive). 
The proof relies in an essential way 
on the theorem that in any fusion category, the identity functor 
is isomorphic to ${****}: V \mapsto V^{****}$
(where $V \mapsto V^*$ is the duality in $\mC$), as a tensor functor, 
which is proved using the theory of weak Hopf algebras
(more specifically, the formula for the fourth power of the antipode
from \cite{N}). We also prove that the $SL_2({\mathbb Z})$-representation 
attached to any modular category over $\CC$ is unitary.

In Section 4, we give a short review of the theory of weak Hopf 
algebras, which are both a 
tool and an object of study in this paper, and prove the isomorphism 
between the identity functor ${\rm Id}$ and ${****}$, 
which is crucial for the main theorem.  

In Section 5, we prove a formula for the trace of squared antipode 
of a semisimple connected regular weak Hopf algebra. 
Then we proceed to prove a number of results (partially due to M\"uger)
about semisimplicity of the category of module functors between 
module categories over a fusion category, in particular of the dual and 
the Drinfeld center. We also generalize these results to multi-fusion 
categories (the notion obtained from that of a fusion category by relaxing 
the condition that $\End(\mathbf 1)=k$). In particular, we prove that a 
semisimple weak Hopf algebra over $k$ is cosemisimple. 
Finally, we prove a categorical version of the class equation 
of Kac and Zhu.

In Section 6, we consider a pair of semisimple weak
Hopf algebras $B\subset A$ and study the subcomplex of 
$B$-invariants in the co-Hochschild complex of $A$
(this complex is a weak analog of the complex considered in \cite{Sch},
\cite{EG1} for Hopf algebras). We prove that 
this complex is acyclic in positive degrees. 

In Section 7, we show that the complex considered in 
Section 6 coincides with 
the deformation complex of a monodial functor between two 
multi-fusion categories introduced independently by Yetter 
\cite{Y1,Y2} and Davydov \cite{Da}. 
Using this, we establish ``Ocneanu rigidity''
(absence of deformations) for arbitrary nondegenerate 
multi-fusion categories and 
tensor functors between them. The idea of the proof of this result 
is due to Ocneanu-Blanchard-Wassermann, 
and the fusion category case was worked out completely 
by Blanchard-Wassermann  in \cite{Wa,BWa} under the assumption that
the global dimension is non-zero. 

In Section 8 we discuss the notion of Frobenius-Perron dimensions 
of objects in a fusion category, show that they are additive and
multiplicative, and have similar properties to categorical
dimensions in a pivotal category.  
We prove a categorical analogue of the Nichols-Zoeller freeness
theorem (saying that a finite dimensional Hopf algebra is a free 
module over a Hopf subalgebra), in particular prove the freeness
theorem for semisimple quasi-Hopf algebras. We also show that 
the global dimension of a fusion category is divisible by its
Frobenius-Perron dimension (in the ring of algebraic integers), 
and the ratio is $\le 1$. In particular, if the Frobenius-Perron
dimension of a fusion category is an integer, then it coincides
with the global dimension. This result may be regarded as a
categorical version of the well known theorem of Larson and
Radford, saying that the antipode of a semisimple Hopf algebra 
is involutive. 
Further, we show that the Frobenius-Perron dimension 
of a fusion category is divisible by the Frobenius-Perron
dimension of its full subcategory, and use this result to show
that any fusion category of Frobenius-Perron dimension $p$ (a prime)
is the category of representations of a cyclic group with a
3-cocycle. We also classify fusion categories of dimension $p^2$.
Finally, we show that the property of a fusion category 
to have integer Frobenius-Perron 
dimensions (which is equivalent to being the representation
category of a quasi-Hopf algebra) is stable under basic operations with categories, in 
particular is ``weak Morita invariant'' (i.e. invariant under passing to the
the opposite of the 
dual category). At the end of the section we define group-theoretical 
fusion categories, which is a subclass of fusion categories with integer
Frobenius-Perron dimensions; 
they are constructed explicitly from finite groups. 

Many of the results of this paper have analogs in positive 
characteristic. These generalizations are given in Section 9. 
In particular, we show that a fusion category of nonzero global
dimension over a field of positive characteristic can be lifted
to characteristic zero. 

Throughout the paper, we are going to freely use the theory 
of rigid tensor categories. We refer the reader to the textbooks 
\cite{K},\cite{BaKi} for details. We also recommend the reader 
the expository paper \cite{CE}, where much of the content 
of this paper is discussed in detail. 

{\bf Acknowledgments.} The research of P.E. was partially supported by
the NSF grant DMS-9988796, and was done in part for the Clay Mathematics 
Institute. 
The research of D.N. was supported by the NSF grant DMS-0200202.
The research of V.O. was supported by the NSF 
grant DMS-0098830. We are grateful to A.~Davydov, 
A.~Kirillov Jr., M.~M\"uger, 
and A.~Wassermann for useful discussions,
and to E.~Blanchard for giving us the formulation of the result of \cite{BWa} 
before publication. 

\section{Results on squared norms, global dimensions, weak Hopf
algebras, and Ocneanu rigidity}

Let $k$ be an algebraically closed field. 
By a {\em multi-fusion} category over $k$ 
we mean a rigid semisimple $k$-linear 
tensor category $\mC$ with finitely many simple objects and finite dimensional 
spaces of morphisms. If the unit object $\bold 1$ of $\mC$ is simple, 
then the category $\mC$ is said to be a {\em fusion} category.   
Otherwise (if $\bold 1$ is not simple), it is easy to see that 
we have $\mathbf 1=\oplus_{i\in J}\mathbf 1_i$, where 
$\mathbf 1_i$ are pairwise nonisomorphic simple objects. 

Let us list a few examples to keep in mind when thinking about
fusion and multi-fusion categories. 

\vskip .05in

{\bf Examples of fusion categories.} 

1. The category ${\rm Vec}_{G}$ 
of finite dimensional vector spaces graded by a finite group $G$
(or, equivalently, finite dimensional modules over the function
algebra ${\rm Fun}(G,k)$.) 
Simple objects in this category
are evaluation modules $V_g$, $g\in G$, and 
the tensor product is given by $V_g\otimes V_h=V_{gh}$, with the
associativity morphism being the identity. 

More generally, pick a 3-cocycle $\omega\in
Z^3(G,k^\times)$. To this cocycle we can attach a
twisted version ${\rm Vec}_{G,\omega}$ of ${\rm Vec}_{G}$: the simple objects 
and the tensor product functor are the
same, but the associativity isomorphism is given by
$\Phi_{V_g,V_h,V_k}=\omega(g,h,k)\textrm{id}$. The
pentagon axiom then follows from the cocycle
condition
$$\omega(h,k,l)\omega(g,hk,l)\omega(g,h,k)=\omega(gh,k,l)\omega(g,h,kl).$$
Note that cohomologous cocycles define equivalent fusion
categories. 

2. The category of finite dimensional $k$-representations
of a finite group $G$, whose order is relatively prime to the
characteristic of $k$. 

3. The category of integrable modules
(from category $\mathcal O$) over the affine algebra $\widehat{sl}_2$ at
level $l$ (see \cite{BaKi}). 
The tensor product in this category 
is the fusion product, defined at the level of objects by the
Verlinde fusion rule
 $$V_i\otimes
V_j=\sum_{\substack{k=\vert i-j\vert
\\ k\equiv i+j~\textrm{mod}~2}}^{l-\vert i+j-l\vert}V_k
$$
\vskip .05in

{\bf Examples of multi-fusion categories.}

1. The category of finite dimensional bimodules 
over a finite dimensional semisimple $k$-algebra,
with bimodule tensor product. It has simple objects 
$M_{ij}$ with ``matrix'' tensor product $M_{ij}\otimes
M_{jk}=M_{ik}$; thus the identity object is 
$\bold 1=\oplus_i M_{ii}$. 

2. The category of 
finite dimensional 
modules over the function algebra on a finite groupoid. 

\subsection{Squared norms and global dimensions}

Let us introduce the notion of the 
global dimension of a fusion category. First of all, 
we have the following known result (see e.g., \cite{O}).

\begin{proposition}\label{doubledual}
In a fusion category, any simple object $V$ is isomorphic to its double dual 
$V^{**}$. 
\end{proposition}

To prove this, it suffices to note that for any simple object $V$, the 
right dual $V^*$ is the unique simple object $X$ 
for which $V\otimes X$ contains the neutral object, while 
the left dual $^*V$ is the unique simple object $X$ 
for which $V\otimes X$ projects to the neutral object;
so $^*V=V^*$ by semisimplicity of the category, and hence $V=V^{**}$
for any simple $V$. 

Next, recall 
(\cite{BaKi}, p.39) 
that if $V\in \mC$, and $g:V\to V^{**}$ is a morphism, then
one can define its ``quantum trace'' $\Tr_V(g)\in k$ by the formula
$$
\Tr_V(g)={\rm ev}_{V^*} \circ (g\otimes 1_{V^*})\circ {\rm coev}_V,
$$ 
where  $\mathbf 1$ is the unit object of  $\mC$,
${\rm coev}_V: \mathbf 1\to V\otimes V^*$ and 
${\rm ev}_V: V^*\otimes V\to \mathbf 1$ are the coevaluation and 
evaluation maps (see \cite{BaKi}). 
It is easy to show (and well known, see e.g. \cite{BaKi}) 
that for any simple object $V$ and a nonzero morphism 
$a: V\to V^{**}$ one has $\Tr_V(a)\ne 0$. 
Indeed, otherwise there is a sequence of nonzero 
maps $\mathbf 1\to V\otimes V^*\to \mathbf 1$ with zero composition, which 
would imply that the multiplicity of $\mathbf 1$ in 
$V\otimes V^*$ is at least 2
-- a contradiction. 

Now, following \cite{Mu1},
for every simple object $V$ of a fusion category $\mC$, define 
the {\em squared norm} 
$|V|^2\in k^\times$ of $V$ as follows. Fix an isomorphism 
$a: V\to V^{**}$ (which exists by Proposition \ref{doubledual}), and let 
$|V|^2=\Tr_V(a)\Tr_{V^*}((a^{-1})^*)$. 
It is clearly independent on the choice of $a$ 
(since $a$ is uniquely determined up to a scaling), and nonzero by the 
explanation in the previous paragraph. For example,  
for the neutral object $\mathbf 1$ we have $|\mathbf 1|^2=1$. 

\begin{definition}\cite{Mu1} The {\em global dimension} of a fusion category $\mC$ 
is the sum of squared norms of its simple objects. It is denoted by 
$\dim(\mC)$. 
\end{definition}

One of our main results is 

\begin{theorem} \label{posi}
If $k=\CC$ then $|V|^2>0$ for all simple objects $V\in \mC$;
therefore $\dim(\mC)\ge 1$, and is $>1$ for any nontrivial $\mC$. 
In particular, for any fusion category $\mC$ 
one has $\dim(\mC)\ne 0$.  
\end{theorem}

\begin{remark} Note that the second statement 
immediately follows from the first one,
since $\mC$ is always defined over a finitely generated 
subfield $k'$ of $k$, which can be embedded into $\CC$. \end{remark}

The proof of Theorem \ref{posi}, which relies on Theorem 
\ref{fourstars} below, is given in Section~3.
Theorem \ref{fourstars} is proved in Section~4. 

\begin{remark} In the course of proof of Theorem \ref{posi}
we show that for any simple object $V$, the number $|V|^2$ is an 
eigenvalue of an integer matrix. In particular, if $k=\CC$ then 
$|V|^2$ is an algebraic integer. Thus, for $k=\CC$ 
Theorem \ref{posi} actually  
implies that $|V|^2$ is a ``totally positive'' algebraic integer,
i.e. all its conjugates are real and 
positive (since one can twist a fusion category 
by an automorphism of $\CC$ and get a new fusion category). 
Similarly, $\dim(\mC)$ is ``totally'' $\ge 1$, i.e. all its 
conjugates are $\ge 1$.\end{remark} 

One of the main tools in the proof of 
Theorem \ref{posi} is 

\begin{theorem} \label{fourstars} 
In any fusion category, the identity functor is isomorphic 
to the functor **** as a tensor functor. 
\end{theorem}

For the category of representations of a Hopf algebra, this result follows 
from Radford's formula for $S^4$ \cite{R1}. 
In general, it follows 
from the analog of Radford's formula for weak Hopf algebras, 
which was proved by the second author in \cite{N}
(see Section 4). 

\begin{definition} A 
pivotal structure on a fusion category $\mC$ is an isomorphism of 
tensor functors $i:\Id\to **$. A fusion category equipped 
with a pivotal structure 
is said to be a pivotal fusion category. 
\end{definition}

We conjecture a stronger form of Theorem \ref{fourstars}:

\begin{conjecture} 
Any fusion category admits a pivotal structure. 
\end{conjecture}

For example, this is true for the representation category 
of a semisimple Hopf algebra over $k$, since 
in this case by the Larson-Radford theorem \cite{LR2}, 
the squared antipode is $1$ and hence $\Id=**$. 
Furthermore, in Section 8 we will show that the conjecture 
is true for the representation category of a 
semisimple quasi-Hopf algebra.  

\subsection{Results on pivotal, spherical, and modular categories}

Now let $\mC$ be a pivotal fusion category. 
In such a category, one can define the dimension of an object
$V$ by $\dim(V)=\Tr_V(i)$, and we have the 
following result which justifies the notation 
$|V|^2$.  

\begin{proposition}\label{pivo} In a pivotal fusion category one has 
$|V|^2=\dim(V)\dim(V^*)$ for any simple object $V$. 
Moreover, if $k=\CC$ then $\dim(V^*)=\overline{\dim(V)}$, so 
$|V|^2=|\dim(V)|^2$.
\end{proposition}

Proposition \ref{pivo} is proved 
in Section 3. 

This result can be further specialized to spherical 
(in particular, ribbon) categories, 
(see \cite{Mu1}). Namely, 
recall from \cite{Mu1} that a pivotal fusion category is spherical if
and only if $\dim(V)=\dim(V^*)$ for all simple objects $V$. 
Thus we have the following corollary. 

\begin{corollary} \label{spheri}
In a spherical category, $|V|^2=\dim(V)^2$.
In particular, if $k=\CC$, then $\dim(V)$ is (totally) real. 
\end{corollary} 

Corollary \ref{spheri} readily follows from Proposition \ref{pivo}. 

\begin{remark} We note that in \cite{Mu1}, the number $|V|^2$ is 
called the squared dimension of $V$, and denoted $d^2(V)$. We do not use 
this terminology for the following reason. 
In a pivotal category over $\CC$, the dimensions 
of simple objects do not have to be real numbers
(it suffices to consider the category of representations 
of a finite group $G$, where $i$ is given
by a nontrivial central element of $G$). Thus, in general
$\dim(V^*)\ne \dim(V)$, and $|V|^2\ne \dim(V)^2$.
Therefore, the term ``squared dimension'', while being adequate 
in the framework of spherical categories of \cite{Mu1}, is no longer adequate
in our more general setting.\end{remark}  

Now assume that $\mC$ is a modular category
(i.e., a ribbon fusion category with a nondegenerate $S$-matrix, 
see \cite{BaKi}). In this case one can define
matrices $S$ and $T$ of algebraic numbers, which yield a projective 
representation of the modular group 
$SL_2(\mathbb Z)$ (\cite{BaKi}). We note that the matrix $S$ is defined only 
up to a sign, since in the process of defining $S$ it is necessary to extract a 
square root of the global dimension of $\mC$. So to be more precise
we should say that by a modular category we mean the underlying 
ribbon category together with a choice of this sign. 

\begin{proposition} \label{modu}
If $k=\CC$ then the projective representation 
of $SL_2(\mathbb Z)$ associated to $\mC$ is unitary in the standard hermitian 
metric (i.e. the matrices $S$ and $T$ are unitary). 
\end{proposition}

Proposition \ref{modu} is proved in Section 3. 

\begin{remark} As before, the proposition actually means that 
this representation is totally unitary, i.e.  
the algebraic conjugates of this representation are unitary as well. 
\end{remark}

\begin{remark}
It is interesting whether Proposition \ref{modu} generalizes to mapping class 
groups of higher genus Riemann surfaces.\end{remark}

We note that the results of this subsection are known for 
hermitian categories \cite{BaKi}. Our results imply that these results
are valid without this assumption. 

\subsection{Module categories, the dual category, the Drinfeld center}

Let $\mM$ be an indecomposable 
left module category over a rigid tensor category 
$\mC$ \cite{O} (all module categories we consider are assumed semisimple).
This means, $\mM$ is a module category which cannot be split in a direct sum 
of two nonzero module subcategories. 
In this case one can define the dual category $\mC^*_{\mM}$ to be the 
category of module functors from $\mM$ to itself:
$\mC^*_{\mM}={\rm Fun}_\mC(\mM,\mM)$ \cite{O}. 
This is a rigid tensor category (the tensor product is the composition
of functors, the right and left duals are the right and left adjoint functors). 

For example, let us consider $\mC$ itself as a module category over 
$\mC\otimes \mC^{op}$, via $(X,Y)\otimes Z=X\otimes Z\otimes Y$. 
Then the dual category is the Drinfeld center $Z(\mC)$ of $\mC$ 
(see \cite{O1, Mu2}; for the basic theory of the Drinfeld center
see \cite{K}). 

The following result 
was proved by M\"uger \cite{Mu1,Mu2}
under the assumption 
$\dim(\mC)\ne 0$ (which, as we have seen, is superfluous in zero characteristic) 
and minor additional assumptions. 

\begin{theorem}\label{dualcat} For any indecomposable module category 
$\mM$ over a fusion category $\mC$,
the category $\mC^*_{\mM}$ is semisimple
(so it is a fusion category), and 
$\dim(\mC^*_{\mM})=\dim(\mC)$. In particular, 
for any fusion category $\mC$ the category $Z(\mC)$
is a fusion category of global dimension $\dim(\mC)^2$.   
\end{theorem}

In fact, we have the following more general results, which 
(as well as the results in the next subsection) are inspired by 
\cite{Mu1,Mu2}.

\begin{theorem}\label{m1m2} For any module categories $\mM_1,\mM_2$ over 
a fusion category $\mC$, 
the category of module functors ${\rm Fun}_{\mC}(\mM_1,\mM_2)$
is semisimple.  
\end{theorem}

Theorems \ref{dualcat} and \ref{m1m2} are proved in Section 5.  

\subsection{Multi-fusion categories}

We say that a multi-fusion category $\mC$ is indecomposable if 
it is not a direct sum of two nonzero multi-fusion categories. 

Let $\mC$ be a multi-fusion category. Then for any simple object 
$X\in \mC$ there exist unique $i,j\in J$ such that 
$\mathbf 1_i\otimes X\ne 0$ and $X\otimes \mathbf 1_j\ne 0$;
moreover, we have $\mathbf 1_i\otimes X=X\otimes \mathbf 1_j=X$. 
Thus, as an additive category, $\mC=\oplus_{m,n}\mC_{mn}$, where 
$\mC_{mn}$ is the full abelian subcategory of $\mC$ with 
simple objects having $i=m,j=n$. It is easy to check 
that $\mC_{ii}$ are fusion categories for all $i$, 
and $\mC_{ij}$ are $(\mC_{ii},\mC_{jj})$-bimodule categories,
equipped with product functors $\mC_{ij}\times \mC_{jl}\to \mC_{il}$ 
satisfying some compatibility conditions. We will refer to $\mC_{ii}$ 
as {\it the component categories} of $\mC$.

Since $\mC_{ii}$ 
are fusion categories, they have well defined global dimensions.
In fact, the off-diagonal subcategories $\mC_{ij}$ can also be assigned global 
dimensions, even though they are not tensor categories. 
To do this, observe that for any simple object $V\in \mC_{ij}$, 
and any morphism $g:V\to V^{**}$ one can define $\Tr_V(g)\in 
\End(\mathbf 1_i)=k$. 
Therefore, we can define $|V|^2$ by the usual formula 
$|V|^2=\Tr_V(g)\Tr_{V^*}((g^{-1})^*)$, and set $\dim(\mC_{ij}):=
\sum_{V\in {\rm Irr}\mC_{ij}}|V|^2$. 

\begin{proposition}\label{samedim}
If $\mC$ is an indecomposable multi-fusion category, then 
all the categories $\mC_{ij}$ have the same global dimensions. 
\end{proposition}

The setup of the previous section can be generalized to the 
multi-fusion case. Namely, we have the following generalization 
of Theorem \ref{m1m2} to the multi-fusion case:

\begin{theorem}\label{multi-fus} 
If $\mC$ is a multi-fusion category, 
and $\mM_1,\mM_2$ are module categories over $\mC$, 
then the category ${\rm Fun}_{\mC}(\mM_1,\mM_2)$ 
is semisimple. In particular, 
the category $\mC^*_\mM$ is semisimple for any module
category $\mM$. 
\end{theorem}

Proposition \ref{samedim} and 
Theorem \ref{multi-fus} is proved in Section 5. 

\begin{remark} We note that it follows from the arguments of
\cite{O,Mu1} that if $\mC$ is a multifusion category and $\mM$ is
a faithful module category of $\mM$, then $\mM$ is a faithful 
module category over $\mC^*_\mM$, and $(\mC^*_\mM)^*_\mM=\mC$.\end{remark} 

\subsection{Results on weak Hopf algebras}

A convenient way to visualize (multi-) fusion categories is 
using weak Hopf algebras (see Section 4 for the definitions). 
Namely, let $\mC$ be a multi-fusion category. 
Let $R$ be a finite dimensional semisimple $k$-algebra, 
and $F$ a fiber functor (i.e. an exact, faithful tensor functor) from 
$\mC$ to $R$-bimod. Let $A=\End_k(F)$ (i.e. the algebra of
endomorphisms of the composition of $F$ with the forgetful
functor to vector spaces). 

\begin{theorem} \label{Schlah} (\cite{Sz})
$A$ has a natural structure of a semisimple  
weak Hopf algebra with base $R$, and $\mC$ is equivalent, as a tensor category,
to the category $\Rep(A)$ of finite dimensional representations of $A$. 
\end{theorem}

\begin{remark} In order 
to lift the naturally existing map $A\to A\otimes_R A$
to a weak Hopf algebra coproduct $A\to A\otimes_k A$, 
one needs to use a separability idempotent 
in $R\otimes R$. We will always use the unique 
{\it symmetric} separability idempotent. 
In this case, the weak Hopf algebra $A$ satisfies 
an additional {\it regularity condition}, saying that
the squared antipode is the identity on the base of $A$.\end{remark} 

One can show that for any multi-fusion category $\mC$,
a fiber functor $F$ exists for a suitable $R$.
Indeed, let $\mM$ be any semisimple faithful module category
over $\mC$, i.e. such that any nonzero object of $\mC$ acts by nonzero 
(for example, $\mC$ itself).
Let $R$ be a semisimple algebra whose blocks are labeled
by simple objects of $\mM$: $R=\oplus_{M\in {\rm Irr}(\mM)}R_M$, 
and $R_M$ are simple. Now define a functor $F:\mC\to R-\text{bimod}$ 
as follows: for any object $X\in \mC$ 
set $$
F(X)=\oplus_{M,N\in {\rm Irr}(\mM)}\Hom_{\mM}(M,X\otimes N)
\otimes B_{NM},
$$ 
where $B_{NM}$ is the simple $(R_N,R_M)$-bimodule. 
This functor has an obvious tensor structure (coming from composition of 
morphisms), and it is clearly exact and 
faithful. So it is a fiber functor. 
 
Therefore, we have

\begin{corollary} \label{Schlah1} (\cite{H},\cite{Sz}) 
Any multi-fusion category is equivalent to the category of 
finite dimensional representations of a (non-unique) regular semisimple 
weak Hopf algebra. This weak Hopf algebra is connected if and only if 
$\mC$ is a fusion category.  
\end{corollary}

\begin{remark} In particular, we can choose $R$ to be a commutative algebra
(i.e., all blocks have size $1$). Thus, any multi-fusion category 
is equivalent to the category of representations of a 
weak Hopf algebra with a commutative base. This result was proved by 
Hayashi in \cite{H}, where weak Hopf algebras with commutative bases 
appear under the name ``face algebras''.\end{remark}

\vskip .05in

The language of weak Hopf algebras is convenient to visualize 
various categorical constructions. One of them is that of 
a dual category. Indeed, let $\mC$ be the category of 
representations of a weak Hopf 
algebra $A$. Let $R$ be the base of $A$. 
Then we have a natural fiber functor from $\mC$ 
to the category of $R$-bimodules -- the forgetful functor. 
This functor defines a natural 
structure of a module category over $\mC$ on the 
category $\mM=R-\text{bimod}$. In this case, the dual category $\mC^*_{\mM}$ 
is simply
the representation category of the dual weak Hopf algebra $A^{*cop}$
with the opposite coproduct (see \cite{O}). 
Furthermore, as we just showed, 
this example is general, in the sense that 
any faithful 
module category over a multi-fusion category can be obtained in this way. 

The global dimension of the category of representations of a 
connected weak Hopf algebra 
is closely related to the trace of squared antipode in this weak Hopf algebra. 
This connection is expressed by the following theorem. 

Let a fusion category $\mC$ be the representation category 
of a regular semisimple weak Hopf algebra $A$. Let $S$ be the
antipode of $A$. Let $q_i$ be the primitive idempotents 
of the center of the base $A_t$, and 
$n_i$ the block sizes of the matrix algebras $q_iA_t$.  

\begin{theorem}\label{traceofS2} 
One has
$\Tr(S^2|_{q_jAS(q_i)})=\dim(\mC)n_i^2n_j^2$. 
In particular, over $\CC$ one always has 
$\Tr(S^2|_{q_jAS(q_i)})>0$. 
\end{theorem}

\begin{remark}\label{tra} Theorem \ref{traceofS2}
implies that $\Tr(S^2|_{A}) = \dim(\Rep(A)) \dim(A_t)^2$.\end{remark}

We also have the following result, which generalizes the Larson-Radford 
theorem for Hopf algebras \cite{LR1}.  

\begin{theorem}\label{sscoss} 
A semisimple weak Hopf algebra over $k$ is also cosemisimple. 
\end{theorem}

Theorem \ref{traceofS2} is proved in Section 5. 
Theorem \ref{sscoss} follows from
Theorem \ref{multi-fus}, since the representation category of the
dual weak Hopf algebra is the dual category. 

\subsection{Vanishing of the Yetter cohomology}

The cohomology theory that controls deformations 
of tensor categories and tensor functors is the cohomology 
theory developed in \cite{CY,Y1,Y2,Da}, which we will call 
the Yetter cohomology theory. 
This theory associates to any tensor category $\mC$ and 
a tensor functor $F$ for $\mC$ to another tensor category 
$\mC'$, a sequence of cohomology spaces $H^i_F(\mC)$
(See Section 7). If $F$ is the 
identity functor, then these groups are simply denoted by 
$H^i(\mC)$. Namely, first order deformations of $\mC$ are elements
of $H^3(\mC)$, while first order deformations of a tensor structure on a 
given functor $F$ between fixed tensor categories $\mC,\mC'$ are
elements of $H^2_F(\mC)$. 
Obstructions to such deformations lie, respectively, in 
$H^4(\mC)$ and $H^3_F(\mC)$.

For multi-fusion categories, the Yetter cohomology often vanishes. 
Namely, we have the following result. 

\begin{theorem}\label{van} 
Let $\mC$, $\mC'$ be multi-fusion categories,
and $F: \mC\to \mC'$ be a (unital) tensor functor 
from $\mC$ to $\mC'$. Then 
$H^i_F(\mC)=0$ for $i>0$.  
\end{theorem}

Theorem \ref{van} is proved in Section 7, using weak Hopf algebras. 
The idea of this proof is due to Ocneanu, Blanchard, and Wassermann.  

\subsection{Ocneanu rigidity}

The statement that a fusion category 
over a field of characteristic zero cannot be nontivially deformed
is known as Ocneanu rigidity, because its formulation and 
proof for unitary categories was suggested (but not published) by Ocneanu. 
Later Blanchard and Wassermann (\cite{BWa,Wa}) proposed an algebraic 
version of Ocneanu's argument, which proves the following result. 

\begin{theorem}
\label{ocnrig} A fusion category
does not admit nontrivial  deformations. In particular, the number of 
such categories (up to equivalence) with a given Grothendieck ring
(=fusion rules) is finite. 
\end{theorem}

\begin{remark} This theorem is proved in \cite{BWa,Wa} 
under the assumption that $\dim(\mC)\ne 0$, but as we know from
Theorem \ref{posi}, this assumption 
is superfluous in zero characteristic.\end{remark} 

In fact, one has the following somewhat more general theorem. 

\begin{theorem}\label{multocnrig}
Theorem \ref{ocnrig} is valid for multi-fusion categories. 
\end{theorem}

This theorem follows directly from Theorem \ref{van} for $i=3$ and $F=\Id$
(this is, essentially, the proof of Ocneanu-Blanchard-Wassermann). 
A sketch of proof is given in Section 7.

Note that Theorem \ref{multocnrig} implies that any multi-fusion category is
defined over an algebraic number field.

{\bf Question.} Is any multi-fusion category defined over a cyclotomic field?

Note that the answer is yes for fusion categories coming from quantum groups
at roots of unity, for the categories of \cite{TY}, 
and for group-theoretical categories discussed in Section 8.
Also, it follows from the results of subsection \ref{cy} that 
the global dimension of a fusion category belongs to a cyclotomic field.

We also have 

\begin{theorem}\label{ocnrigfun}
A (unital) tensor functor between multi-fusion categories
does not have nontrivial deformations. 
In particular, the number of 
such functors (up to equivalence) for fixed source and target 
categories is finite. 
\end{theorem}

\begin{remark} In particular, Theorem \ref{ocnrigfun} implies that a given
fusion category has finitely many fiber functors. This is also shown in
\cite{IK} in the special case of categories coming from subfactors
(see Theorem 2.4 in \cite{IK}).\end{remark}

\begin{remark} Theorem \ref{ocnrigfun} also implies Ocneanu's
result that a given fusion category has finitely many braidings. 
Indeed, any braiding in particular defines an equivalence of tensor
categories between $\mC$ and $\mC^{op}$.\end{remark} 

The proof of Theorem \ref{ocnrigfun} is analogous to the proof of 
Theorem \ref{multocnrig}. Namely, the first part of the result
follows directly from Theorem \ref{van} for $i=2$. The second part is a
consequence of a combinatorial fact that there are only finitely many
homomorphisms between the Grothendieck rings of the categories under
consideration; this fact is a consequence of Proposition 2.1 of \cite{O}
since the second Grothendieck ring can be considered as a based module over
the first one.

\begin{remark} One says that an indecomposable multi-fusion category $\mC'$ is a quotient
category of an indecomposable multi-fusion category $\mC$ if there exists a tensor functor
$F:\mC \to \mC'$ such that any object of $\mC'$ is a subobject of $F(X),\;
X\in \mC$ (See Section 5). It follows from the above that a given multi-fusion category
$\mC$ has only finitely many quotient categories.\end{remark}

\begin{corollary} 
\label{ocnrigmod} A module category $\mM$ over 
a multi-fusion category $\mC$
does not admit nontrivial  deformations. In particular, the number of 
such module categories (up to equivalence) with a given 
number of simple objects is finite. 
\end{corollary}

To prove Corollary \ref{ocnrigmod}, it suffices to 
choose a semisimple algebra $R$ with number of blocks equal to the number of
simple objects in $\mM$, and apply Theorem \ref{ocnrigfun} to the functor
$F:\mC \to R-\text{bimod}$ associated with $\mM$. 

\begin{remark} We note that this result was conjectured in \cite{O}. 
\end{remark}
\section{Fusion categories over $\CC$}

The goal of this subsection is to prove Theorem \ref{posi},
Proposition \ref{pivo}, Proposition \ref{modu}, and discuss 
their consequences.

Let $\mC$ denote a fusion category over $\CC$. We will denote representatives 
of isomorphism classes of simple 
objects of $\mC$ by $X_i$, $i\in I$, with $0\in I$ labeling the 
neutral object $\mathbf{1}$. Denote by $*: I\to I$ 
the dualization map. (We note that at the level of labels, there is
no distinction between right and left dualization, since
by Proposition \ref{doubledual}, $V$
is isomorphic to $V^{**}$ for any simple object $V$). 

\subsection{Proof of Theorem \ref{posi}}

By Theorem \ref{fourstars}, there exists 
an isomorphism of tensor functors 
$g: \Id\to ****$. Let us fix such an isomorphism. 
So for each simple $V\in \mC$ we have 
a morphism $g_V: V\to V^{****}$ such that 
$g_{V\otimes W}=g_V\otimes g_W$, and $g_{V^*}=(g_V^{-1})^*$. 
For all simple $V$, let us fix an isomorphism $a_V : V\to V^{**}$
such that $a_V^2=g_V$ ($a_V^2$ makes sense as 
${\rm Hom}(V,V^{**})$ is canonically isomorphic 
to ${\rm Hom}(V^{**},V^{****})$, and so $a_V$ can also be regarded 
as an element of ${\rm Hom}(V^{**},V^{****})$). 
Then $a_{V^*}=\epsilon_V (a_V^{-1})^*$, where $\epsilon_V=\pm 1$, 
and $\epsilon_V\epsilon_{V^*}=1$. Define 
$d(V)=\Tr_V(a_V)$, in particular $d_i=d(X_i)$. 
These numbers are not at all canonical -- they depend on the 
choices we made. Nevertheless, it is easy to see 
from the definition that $|V|^2=\epsilon_V d(V)d(V^*)$.  

Now let us take any simple objects $V,W\in \mC$ and 
consider the tensor product $a_V\otimes a_W: V\otimes W\to (V\otimes W)^{**}$. 
We have $V\otimes W=\oplus_U {\rm Hom}(U,V\otimes W)\otimes U$, 
so we have $a_V\otimes a_W=\oplus b_{U}^{VW}\otimes a_U$, where 
$b_U^{VW}: {\rm Hom}(U,V\otimes W)\to {\rm Hom}(U^{**},V^{**}\otimes W^{**})$.
The source and target of this map are actually canonically isomorphic, so we 
can treat $b_U^{VW}$ as an automorphism of ${\rm Hom}(U,V\otimes W)$.
Also, since $a_V^2=g_V$ and $g_{V\otimes W}=g_V\otimes g_W$, we get 
$(b_U^{VW})^2=1$. 

Let $b_j^{ik}:=b_{X_j}^{X_iX_k}$. Let $T_j^{ik}=\trace(b_j^{ik})$. 
Since 
$(b_j^{ik})^2=1$, its eigenvalues are $\pm 1$, hence $T_j^{ik}$ are integers. 

By the definition of $b_j^{ik}$, we have 
$$
d_id_k=\sum_{j\in I}T_j^{ik}d_j
$$
Let $T_i$ be the matrix such that $(T_i)_{kj}=T_j^{ik}$, and let 
$\mathbf d$ be the vector with components $d_i$. This vector is clearly 
nonzero (e.g., $d_0\ne 0$), and we have 
$$
T_i\mathbf d=d_i\mathbf d.
$$
Thus, $\mathbf d$ is a common eigenvector of $T_i$ with eigenvalue $d_i$. 

Let $\epsilon_i=\epsilon_{X_i}$. Then we have: 
$\mathbf d$ is an eigenvector of the matrix $A_i=\epsilon_i T_iT_{i^*}$
with eigenvalue $|X_i|^2$. 

We claim that $T_{i^*}=\epsilon_i T_i^t$. This implies the theorem, as 
in this case $A_i=T_iT_i^t$, which is nonnegative definite
(since the entries of $T_i$ are integer, hence real), so its 
eigenvalues, including $|X_i|^2$, have to be nonnegative. 

To prove that $T_{i^*}=\epsilon_i T_i^t$, we need to show that 
$T_{j}^{i^*k}=\epsilon_i T_k^{ij}$. But this, in turn, follows from
the equality $b_j^{i^*k}=\epsilon_i (b_k^{ij})^*$
(which makes sense since we have canonical 
isomorphisms ${\rm Hom}(U,V\otimes W)\to {\rm Hom}(V^*\otimes U,W)\to 
{\rm Hom}(W,V^*\otimes U)^*$, the latter by semisimplicity). 

It remains to prove the last equality. 
Let $f: X_j\to X_i\otimes X_k$ be a morphism. 
By the definition, $b_j^{ik}f=(a_i\otimes a_k)fa_j^{-1}$. 
On the other hand, $f$ gives rise to an element 
$f': X_{i^*}\otimes X_j\to X_k$, given by 
$f'=({\rm ev}_i\otimes 1)(1\otimes f)$, and hence to a linear function 
on ${\rm Hom}(X_k,X_{i^*}\otimes X_j)$ given by 
$f''(h)=({\rm ev}_i\otimes 1)(1_{i^*}\otimes f)h\in {\rm End}(X_k)=\CC$
(for brevity we omit associativity maps and 
use the same notation for a morphism and its  
double dual).
Now, we have 
\begin{eqnarray*}
(b_j^{ik}f)''(h)
&=& ({\rm ev}_i\otimes 1)(1_{i^*}\otimes (a_i\otimes a_k)fa_j^{-1})h \\
&=& \epsilon_i a_k({\rm ev}_i\otimes 1)(a_{i^*}^{-1}\otimes fa_j^{-1})h \\
&=& \epsilon_i  ({\rm ev}_i\otimes 1)(1\otimes f)(a_{i^*}^{-1}\otimes a_j^{-1})ha_k \\
&=& \epsilon_i ({\rm ev}_i\otimes 1)(1\otimes f)(a_{i^*}\otimes a_j)ha_k^{-1},
\end{eqnarray*}
where the last equality due to the facts that $g_{V\otimes W}=g_V\otimes g_W$,
and $a_V=g_Va_V^{-1}$. 
Thus, 
$$
(b_j^{ik}f)''(h)=
f''(b_k^{i^*j}h),
$$
as desired. Theorem \ref{posi} is proved. 

\begin{remark}\label{pivotalization}
Here is a modification of the above proof, which is
simpler but requires an extension of the category $\mC$.

The first step is to prove the result in the case of pivotal
categories. This is a substantially simplified version of the
above proof, since in the pivotal case 
$\epsilon_i$ and $b_j^{ik}$ are 1, so $T_j^{ik}$ are equal to 
$N_j^{ik}$, the multiplicity of $X_j$ in $X_i\otimes X_k$. 

The second step is reduction to the pivotal case. To do this, 
take the isomorphism $g: Id\to ****$ as tensor functors
constructed in the proof of Theorem \ref{fourstars}. 
Define the fusion category $\tilde \mC$, whose simple  
objects are pairs $(X,f)$, where $X$ is a simple object of $\mC$,
and $f:X\to X^{**}$ is an isomorphism, such that $f^{**}f=g$. 
It is easy to see that this fusion category has
a canonical pivotal structure. Thus, $|X|^2>0$ in 
$\tilde\mC$. However, we have the forgetful functor 
$F:\tilde \mC\to \mC$, $F(X,f)=X$, which obviously 
preserves squared norms. Hence, $|V|^2>0$ for simple objects 
$V$ in $\mC$ (since for any $V$ there exists $f:V\to
V^{**}$ such that $f^{**}f=g$), which completes the proof. 
\end{remark} 

\subsection{Proof of Proposition \ref{pivo}}

The first statement is straightforward. To prove the second statement, 
let $N_i$ be the matrix whose entries are given by
$(N_i)_{kj}={\rm dim}{\rm Hom}(X_j,X_i\otimes X_k)$. 
Then, like in the proof of Theorem \ref{posi}, 
$N_i\mathbf d=d_i\mathbf d$. 
Let us multiply this equation on the left by $\mathbf d^*$.
(the hermitian adjoint). We have 
$\mathbf d^*N_i=(N_i^t\mathbf d)^*=
(N_{i^*}\mathbf d)^*=\overline{d_{i^*}}\mathbf d^*$. 
Thus the multiplication will yield
$$
\overline{d_{i^*}}|\mathbf d|^2=d_i|\mathbf d|^2. 
$$
Since $|\mathbf d|^2>0$, this yields the second statement. 
The proposition is proved. 

\subsection{Proof of Proposition \ref{modu}
and properties of modular categories}

It is well known (see \cite{BaKi}) that the matrix $T$ is diagonal
with eigenvalues being roots of unity. So we only need to show that 
the S-matrix is unitary. The entries $s_{ij}$ of the S-matrix are known to 
satisfy $s_{ij}=s_{ji}$, and $\sum_j s_{i^*j}s_{jl}=\delta_{il}$ 
(\cite{BaKi}). So it suffices to show that $s_{i^*j}=\overline{s_{ij}}$. 

Recall that $s_{j0}=\pm \dim(V_j)/\sqrt{D}$, where $D$ is the global dimension of 
the category. So $s_{j0}$ is a real number. 
Further, by the Verlinde formula (\cite[3.1]{BaKi}), 
if $\mathbf v_m$ is the $m$-th column of $S$ then 
$$
N_i\mathbf v_m=\frac{s_{im}}{s_{0m}}\mathbf v_m, 
$$
for all $i$. Let $\mathbf v_m^*$ denote the row 
vector, Hermitian adjoint to $\mathbf v_m$. 
Then 
$$
\mathbf v_m^*N_i\mathbf v_m=\frac{s_{im}}{s_{0m}}|\mathbf v_m|^2.
$$
One the other hand, since $N_i^T=N_{i^*}$, and $N_i$ is real, we have
$$
\mathbf v_m^*N_i\mathbf v_m=(N_i^T\mathbf v_m)^*\mathbf v_m=(N_{i^*}\mathbf v_m)^*\mathbf v_m=
(\frac{s_{i^*m}}{s_{0m}}\mathbf v_m)^*\mathbf v_m=
\frac{\overline{s_{i^*m}}}{s_{0m}}|\mathbf v_m|^2.
$$
Since $|\mathbf v_m|^2\ne 0$, the statement is proved. 

\begin{remark}
An important property of modular categories is 

\begin{proposition}\label{moddivi} (\cite{EG4}, Lemma 1.2) The dimension 
of every simple object of a modular category divides 
$\sqrt{D}$, where $D$ is the global dimension 
of the category. 
\end{proposition} 

This implies that in a ribbon category, dimensions of objects 
divide the global dimension (this is seen by taking the Drinfeld center 
and using Theorem \ref{dualcat} and 
the well known fact that the Drinfeld center of 
a braided category contains the original category). This fact is useful 
in classification of modular categories (See section 8).  
\end{remark}

In conclusion of this subsection let us give a simple application 
of Theorem \ref{posi} to modular categories. 

Let $c\in {\mathbb C}/8{\mathbb Z}$ be the Virasoro central
charge of a modular category $\mC$, defined by the relation 
$(ST)^3=e^{2\pi ic/8}S^2$ 
(see \cite{BaKi}), and $D$ be the global dimension of 
$\mC$. Let $f: (-\infty,1]\to [1,\infty)$ be the inverse function to 
the monotonely decreasing function $g(x)=(3-x)\sqrt{x}/2$. 

\begin{proposition} One has
$$
D\ge f({\rm sign}(s_{00})\cos(\pi c/4))
$$
\end{proposition}

\begin{remark} The number $s_{00}$ equals $\pm 1/\sqrt{D}$.\end{remark}

\begin{proof} We need to show that
$$
{\rm sign}(s_{00})\cos(\pi c/4)\ge (3-D)\sqrt{D}/2
$$
To prove this, recall \cite{BaKi} that the number $e^{\pi ic/4}$ can be
expressed via a Gaussian sum:
$$
{\rm sign}(s_{00})e^{\pi i c/4}=\sum |X_i|^2\theta_i/\sqrt{D}, 
$$
where $\theta_i$ are the twists. 
This means that 
$$
{\rm sign}(s_{00})e^{\pi i c/4}-D^{-1/2}=
\sum_{i\ne 0} |X_i|^2\theta_i/\sqrt{D}. 
$$
Taking squared absolute values of both sides, and using that $|\theta_i|=1$, 
and $|X_i|^2>0$ (by Theorem \ref{posi}), we get 
$$
|{\rm sign}(s_{00})
e^{\pi ic/4}-D^{-1/2}|^2\le (\sum_{i\ne 0}|X_i|^2)^2/D=\frac{(D-1)^2}{D}.
$$
After simplifications, this yields the desired inequality.
\end{proof}

\begin{remark} Sometimes this inequality sharp. For example,  
for $c=\pm 2/5,\pm 22/5$, it yields an estimate 
$D\ge (5-\sqrt{5})/2$, 
which becomes an equality for 
the well known Yang-Lee category 
($c=-22/5$, see \cite{DMS}) and three its Galois images.
The same is true for $c=\pm 14/5,\pm 34/5$, in which case the estimate is 
$D\ge (5+\sqrt{5})/2$ (it is attained at 
four Galois images of the Yang-Lee category). 
Also, it goes without saying that the inequality 
is sharp for $c=0$, where it becomes an equality for the trivial category. 
\end{remark}

These considerations inspire the following questions.

{\bf Questions.} Does there exist a sequence of 
nontrivial fusion categories over $\CC$ 
for which the global dimensions tend to $1$? Does there 
exist such a sequence with a bounded  number of simple objects?
Is the set of global dimensions of fusion categories over $\CC$ 
a discrete subset of $\mathbb R$? 

We expect that the answer to the first question is ``no'', i.e. the point $1$ 
of the set of possible dimensions of fusion categories is an isolated point
(``categorical property T''). 

\begin{section}{Weak Hopf algebras and their integrals}

Below we collect the definition and basic properties of
weak Hopf algebras.

\subsection{Definition of a weak Hopf algebra}

\begin{definition}[\cite{BNSz}] 
\label{finite weak Hopf algebra}
A {\em weak Hopf algebra} is a vector space $A$
with the structures of an associative algebra $(A,\,m,\,1)$ 
with a multiplication $m:A\otimes_k A\to A$ and unit $1\in A$ and a 
coassociative coalgebra $(A,\,\Delta,\,\epsilon)$ with a comultiplication
$\Delta:A\to A\otimes_k A$ and counit $\epsilon:A\to k$ such that:
\begin{enumerate}
\item[(i)] The comultiplication $\Delta$ 
is a (not necessarily unit-preserving) homomorphism of algebras:
\begin{equation}
\label{Delta m}
\Delta(hg) = \Delta(h) \Delta(g), \qquad h,g\in A,
\end{equation}
\item[(ii)] The  unit and counit satisfy the following identities: 
\begin{eqnarray}
\label{Delta 1}
(\Delta \otimes \id) \Delta(1) 
& =& (\Delta(1)\otimes 1)(1\otimes \Delta(1)) 
= (1\otimes \Delta(1))(\Delta(1)\otimes 1), \\
\label{eps m}
\epsilon(fgh) &=& \epsilon(fg\1)\, \epsilon(g\2h) = \epsilon(fg\2)\, 
\epsilon(g\1h),
\end{eqnarray}
for all $f,g,h\in A$. 
\item[(iii)]
There is a linear map $S: A \to A$, called an {\em antipode}, such that
\begin{eqnarray}
m(\id \otimes S)\Delta(h) &=&(\epsilon\otimes\id)(\Delta(1)(h\otimes 1)),
\label{S epst} \\
m(S\otimes \id)\Delta(h) &=& (\id \otimes \epsilon)((1\otimes h)\Delta(1)),
\label{S epss} \\
S(h) &=& S(h\1)h\2 S(h\3),
\label{S id S}
\end{eqnarray}
for all $h\in A$.
\end{enumerate} 
\end{definition}

\begin{remark} We use Sweedler's notation for a comultiplication:  
$\Delta(c) = c\1 \otimes c\2$.\end{remark} 

Axioms (\ref{Delta 1})  and (\ref{eps m}) above 
are  analogous to the usual bialgebra axioms of 
$\Delta$ being a unit preserving map and $\epsilon$ being an
algebra homomorphism. Axioms (\ref{S epst}) and (\ref{S epss}) 
generalize the properties of the antipode in a Hopf algebra
with respect to the counit. Also, it is possible to show
that given (\ref{Delta m}) - (\ref{S epss}),
axiom (\ref{S id S}) is equivalent to $S$ being both
anti-algebra and anti-coalgebra map. 

The antipode of a finite dimensional weak Hopf algebra is
bijective \cite[2.10]{BNSz}.

\begin{remark}
A weak Hopf algebra is a Hopf algebra if and only if the comultiplication 
is unit-preserving and if and only if $\epsilon$ is a homomorphism of algebras.
\end{remark}

A {\em morphism} between weak Hopf algebras $A_1$ and $A_2$
is a map $\phi : A_1 \to A_2$ which is both algebra and coalgebra
homomorphism preserving $1$ and $\epsilon$
and which intertwines the antipodes of $A_1$ and $A_2$,
i.e., $ \phi\circ S_1 = S_2\circ \phi$. The image of a morphism is
clearly a weak Hopf algebra.

When $\dim_k A < \infty$, there is a natural weak Hopf algebra
structure on the dual vector space $A^*=\Hom_k(A,k)$ given by
\begin{eqnarray}
& & \phi\psi(h) = (\phi\otimes\psi)(\Delta(h)), \\
& & \Delta(\phi)(h\otimes g) =  \phi(hg) \\
& & S(\phi)(h) = \phi(S(h)),
\end{eqnarray}
for all $\phi,\psi \in A^*,\, h,g\in A$. The unit
of $A^*$ is $\epsilon$  and the counit is $\phi \mapsto \phi(1)$.

The linear maps defined in (\ref{S epst}) and (\ref{S epss})
are called {\em target\,} and {\em source counital maps}
and are denoted $\eps_t$ and $\eps_s$ respectively :
\begin{equation}
\eps_t(h) = \epsilon(1\1 h)1\2, \qquad
\eps_s(h) = 1\1 \epsilon(h1\2),
\end{equation}
for all $h\in A$. The images of the counital maps
\begin{equation}
A_t = \eps_t(A),  \qquad A_s = \eps_s(A)
\end{equation}
are separable subalgebras of $A$, called target and source {\em bases} 
or {\em counital subalgebras} of $A$.  
These subalgebras commute with each other; moreover
\begin{eqnarray*}
A_t &=& \{(\phi\otimes \id)\Delta(1) \mid \phi\in A^* \} 
= \{ h\in A \mid \Delta(h) = \Delta(1)(h\otimes 1) \}, \\
A_s &=& \{(\id \otimes \phi)\Delta(1) \mid \phi\in A^* \}
= \{ h\in A \mid \Delta(h) = (1\otimes h)\Delta(1) \},
\end{eqnarray*}
i.e., $A_t$ (respectively, $A_s$) is generated 
by the right (respectively, left) tensor factors of $\Delta(1)$
in the shortest possible presentation of $\Delta(1)$ in $A\otimes_k A$. 

The category $\Rep(A)$ of left $A$-modules is a rigid tensor category 
\cite{NTV}. The tensor product of two $A$-modules $V$ and $W$ is 
$V \otimes_{A_t} W$ with the $A$-module structure defined 
via $\Delta$. The unit object $\mathbf{1}$ of $\Rep(A)$ is the target counital 
algebra $A_t$ with the action 
$h\cdot z = \eps_t(hz)$ for all $h\in A,\, z\in A_t$.

For any algebra $B$ we denote 
by $Z(B)$ the center of $B$.
Then the unit object of $\Rep(A)$ is irreducible if and only if $Z(A)\cap A_t = k$.
In this case we will say that $A$ is {\em connected}.
We will say that that $A$ is {\em coconnected }
if $A^*$ is connected, and that $A$ is {\em biconnected} if it is both connected 
and coconnected.

If $p\neq 0$ is an idempotent in $A_t\cap A_s\cap Z(A)$, then
$A$ is the direct sum of weak Hopf algebras 
$pA$ and $(1-p)A$. Consequently, we say that $A$ is {\em indecomposable}
if $A_t\cap A_s\cap Z(A) = k1$.

Every weak Hopf algebra $A$ contains a canonical {\em minimal\,}
weak Hopf subalgebra $A_{{\rm min}}$ generated, as an algebra, 
by $A_t$ and $A_s$ \cite[Section 3]{N}. In other words,
$A_{{\rm min}}$ is the minimal weak Hopf subalgebra of $A$
that contains $1$. Obviously, $A$ is an ordinary Hopf algebra
if and only if  $A_{{\rm min}} = k1$. Minimal weak Hopf algebras over $k$,
i.e., those for which  $A= A_{{\rm min}}$,
were completely classified in \cite[Proposition 3.4]{N}.

The restriction of $S^2$ on $A_{{\rm min}}$ is always an inner automorphism
of $A_{{\rm min}}$, see \cite{N}. 

\begin{remark}
Unless indicated otherwise, in what follows
we will use only weak Hopf algebras satisfying 
the following natural property :
\begin{equation}
\label{regularity property}
S^2|_{A_{{\rm min}}} = \id.
\end{equation}
This property has an easy categorical interpretation. Let $\mathbf{1} = A_t$
be the trivial $A$-module. Then \eqref{regularity property} is satisfied
if and only if the canonical isomorphism $\mathbf{1} \to \mathbf{1}^{**}$ is 
the identity map.\end{remark}

\begin{definition}
We will call a weak Hopf algebra satisfying \eqref{regularity property}
{\em regular}. 
\end{definition}

\begin{remark} It was shown in \cite[6.1]{NV} that every weak Hopf algebra
can be obtained as a twisting of some regular weak Hopf algebra
with the same algebra structure. \end{remark}

\subsection{Integrals}
\label{subsection about integrals}

The following notion of an integral in a weak Hopf algebra is a generalization
of that of an integral in a usual Hopf algebra \cite{M}.

\begin{definition}[\cite{BNSz}]
\label{integral}
A left {\em integral} in $A$ is an element $\ell\in A$ such that 
\begin{equation}
h\ell =\eps_t(h)\ell, \qquad (rh = r\eps_s(h)) \qquad \mbox{ for all } h\in A. 
\end{equation}
\end{definition}
The space of left integrals in $A$ is a left ideal of $A$ of dimension equal to
$\dim_k(A_t)$.

Any left integral $\lambda$ in $A^*$ satisfies the following
invariance property :
\begin{equation}
\label{left invariance}
g\1 \lambda(hg\2) = S(h\1)\lambda(h\2g),
\qquad g,h\in A.
\end{equation}

In what follows we use the Sweedler arrows $\actl$ and $\actr$
for the dual actions :
\begin{equation}
\label{Sweedler arrows}
h\actl\phi (g) = \phi(gh)
\quad \text { and } \quad
\phi\actr h (g)= \phi(hg).
\end{equation}   
for all $g,\,h\in A,\phi\in A^*$.

Recall that a functional $\phi\in A^*$ is {\em non-degenerate}
if $\phi\circ m$ is a non-degenerate bilinear form on $A$. 
Equivalently, $\phi$ is non-degenerate if
the linear map $h\mapsto (h \actl \phi)$
is injective. An integral (left or right) in a weak Hopf algebra
$A$ is called {\em non-degenerate} if
it defines a non-degenerate functional on $A^*$. A left integral $\ell$
is called {\em normalized} if $\eps_t(\ell)=1$. 

It was shown by P.~Vecsernyes \cite{V} that a finite dimensional
weak Hopf algebra always possesses a non-degenerate left integral. 
In particular, a finite dimensional weak Hopf algebra is a
Frobenius algebra (this extends the well-known result of Larson
and Sweedler for usual Hopf algebras).

Maschke's theorem for weak Hopf algebras, proved in \cite[3.13]{BNSz},
states that a weak Hopf algebra $A$ is semisimple if and only if  $A$ is
separable, and if and only if  there exists a normalized left integral in $A$.
In particular, every semisimple weak Hopf algebra is finite dimensional. 

For a finite dimensional $A$ there is a useful notion of duality between 
non-degenerate
left integrals in $A$ and $A^*$ \cite[3.18]{BNSz}. If $\ell$
is a non-degenerate left integral in $A$ then there exists 
a unique $\lambda\in A^*$ such that
$\lambda\actl \ell = 1$. This $\lambda$ 
is a non-degenerate left integral in $A^*$. Moreover,
$\ell\actl \lambda =\epsilon$. Such a pair of non-degenerate integrals 
$(\ell,\,\lambda)$ is called a pair of {\em dual} integrals.

An invertible element $g\in A$ is called {\em group-like}
if $\Delta(g) = (g\otimes g)\Delta(1) = \Delta(1)(g\otimes g)$ 
\cite[Definition 4.1]{N}. Group-like elements of $A$ form a group $G(A)$
under multiplication. This group has a normal subgroup $G_0(A) := G(A_{{\rm min}})$ 
of {\em trivial} group-like elements. If $A$ is finite dimensional,
the quotient group  $\widetilde{G}(A) = G(A)/ G_0(A)$ is finite.
It was shown in \cite{N} that if $A$ is 
finite dimensional
and $\ell\in A$ and $\lambda\in A^*$ is a dual pair of left integrals,
then there exist group-like elements $\alpha\in G(A^*)$ and $a\in G(A)$,
called {\em distinguished} group-like elements,
whose classes in $\widetilde{G}(A^*)$ and $\widetilde{G}(A)$
do not depend on the choice of $\ell$ and $\lambda$, such that
\begin{equation}
S(\ell) = \alpha\actl \ell \qquad \mbox{ and  } \qquad S(\lambda)=a \actl \lambda.
\end{equation} 
(Note that $\alpha$ and $a$ themselves depend on the choice of $\ell$ and $\lambda$).

The following result is crucial for this paper. 
It is an analogue of
Radford's formula for usual Hopf algebras.

\begin{theorem}\label{S4} \cite[Theorem 5.13]{N} One has
\begin{equation}
S^4(h) = a^{-1}(\alpha\actl h \actr \alpha^{-1})a. 
\end{equation}
for all $h\in A$.
\end{theorem}

A weak Hopf algebra $A$ is said to be {\em unimodular} if it contains
a non-degenerate $2$-sided integral; equivalently, $A$ is unimodular
if the coset of $\alpha$ in $G(A^*)$ is trivial. A semisimple weak Hopf 
algebra is automatically unimodular \cite{BNSz}.

\subsection{Proof of Theorem \ref{fourstars}}\label{pffourstars}

Let $\mC$ be a fusion category over $k$.
By Theorem \ref{Schlah1}, $\mC$ is the representation category of a 
semisimple connected weak Hopf algebra $A$. Since $A$ is semisimple, it is 
unimodular, and hence one can choose a left integral $\ell$ such that
the corresponding element $\alpha=\epsilon$. Thus, by Theorem \ref{S4}, 
one has $S^4(h)=a^{-1}ha$, 
for some group-like element $a\in G(H)$, and hence 
for any simple object $V\in \mC=\Rep(A)$, the operator $a^{-1}|_V$
defines an isomorphism $V\to V^{****}$. Since $a$ is grouplike, 
this isomorphism respects the tensor product. 
The theorem is proved. 

\subsection{Canonical integrals}

Using the squared antipode, one can construct integrals 
on a weak Hopf algebra $A$ called the {\em canonical integrals}. 
Namely, we have the following proposition, a particular case of which 
was proved in \cite{BNSz}. 

\begin{proposition}
\label{chi}
Let $A$ be a finite-dimensional weak Hopf algebra
and let $\{q_i\}$ be primitive central idempotents in $A_t$.
Then $\chi_i(h) := \Tr(L_h\circ S^2|_{AS(q_i)})$, where $L_h \in \End_k(A)$
is given by $L_h(x)=hx,\,h, x\in A$, defines a left integral in $A^*$
such that 
\begin{eqnarray}
\label{chi1}
\chi_i(q_j) &=& \Tr(S^2|_{q_jAS(q_i)})  \qquad \mbox{and } \\
\label{q-trace}
\chi_i(gh) &=& \chi(h S^2(g)),
\end{eqnarray} 
for all $g,h\in A$.
\end{proposition}
\begin{proof}
The identities \eqref{chi1} and \eqref{q-trace} follow immediately 
from the definition of $\chi$, since $A$ is regular and 
$S^2\circ L_g = L_{S^2(g)}\circ S^2$ for all $g\in A$.

Next, for all $\phi\in A^*$ we define $R_\phi \in \End_k(A)$
by $R_\phi(x) = x\actr \phi,\, x\in A^*$. 
Then $R_\phi(xS(q_i)) = R_\phi(x)S(q_i)$ so that $AS(q_i)$ is
$R_\phi$-invariant. A direct computation shows
that  $L_{(h \actr \phi)} = R_{\phi\1}\circ L_h \circ R_{S^{-1}(\phi\2)}$
for all $h\in A,\, \phi\in A^*$,
and hence,
\begin{eqnarray*}
(\phi \otimes \chi_i)(\Delta(h))
&=& \Tr(L_{(h \actr \phi)}\circ S^2|_{AS(q_i)}) \\
&=& \Tr( R_{\phi\1}\circ L_h \circ R_{S^{-1}(\phi\2)}\circ S^2|_{AS(q_i)}) \\
&=& \Tr( R_{\eps_t(\phi)}\circ  L_h \circ S^2|_{AS(q_i)}),
\end{eqnarray*}
therefore, $  (\id \otimes \chi_i)\Delta(h) =  
(\eps_t \otimes \chi)\Delta(h)$ for all $h\in A$, i.e., $\chi_i$
is a left integral in $A^*$.
\end{proof}

\subsection{Sufficient conditions for semisimplicity}

The next proposition is a refinement of \cite[Proposition 6.4]{N}
and is analogous to the Larson-Radford theorem \cite{LR1} for usual
Hopf algebras which says that if $A$ is a finite dimensional Hopf algebra
with $\Tr(S^2|_A)\neq 0$ then $A$ is semisimple and cosemisimple.

\begin{proposition}
\label{semisimplicity condition}
Let $A$ be a finite dimensional weak Hopf algebra. If for every
primitive idempotent $p$ in $A_t \cap Z(A)$ there exist 
idempotents $q_p,\, q_p'$ in $pA_t$ such that 
\begin{equation}
\label{Tr S2 neq 0}
\Tr(S^2|_{q_pAS(q_p')}) \neq 0,
\end{equation}
then $A$ is semisimple.
\end{proposition}
\begin{proof}
Let $\ell$ and $\lambda$ be a dual pair of non-degenerate left integrals
in $A$ and $A^*$ respectively. Then $(\ell\2 \actl \lambda) \otimes
q_p S^{-1}(\ell\1) S(q_p')$ is the dual bases tensor for $q_p A S(q_p')$,
cf.\ \cite{N}. Therefore,
\begin{eqnarray*}
0 &\neq& \Tr(S^2|_{q_pAS(q_p')}) \\
&=& \la \ell\2 \actl \lambda,\, q_p S(\ell\1)S(q_p') \ra \\
&=& \la \lambda,\, q_p \eps_s(q_p' \ell) \ra.
\end{eqnarray*}  
Therefore, $y_p:=\eps_s(q_p'\ell) \neq 0$. Next, we compute
\begin{eqnarray*}
y_p \ell
&=& S(\ell\1) S(q_p')\ell\2 \ell \\
&=& S(\ell\1) S(q_p') \eps_t(\ell\2) \ell \\
&=& S(\ell\1) \eps_t(\ell\2) S(q_p') \ell = S(\ell) y_p.
\end{eqnarray*}
Let $\xi_p$ be an element in $A^*_t$ 
(here and below we write $A^*_t:=(A^*)_t$, $A^*_s:=(A^*)_s$)
such that 
$\xi_p\actl 1 =y_p$. Then from the previous computation we have
\begin{equation}
y_p \ell = (\alpha \actl \ell) y_p,
\end{equation}
where $\alpha$ is a distinguished group-like element of $A^*$.
Arguing as in \cite[Proposition 6.4]{N} one gets
\begin{equation}
\xi_p = S(\xi_p) \alpha,
\end{equation}
whence $\alpha$ is trivial, i.e., $\alpha\in A^*_{\rm min}$.
This means that one can choose $\ell$ in such a way that $\alpha =\epsilon$,
i.e., $\ell = S(\ell)$. Hence, $y_p\in A_s \cap Z(A)$ is a non-zero multiple
of $p$. We conclude that $y :=\eps_s(\ell) =\sum_p\, y_p$ is an invertible
central element of $A$. Hence, $\ell' = y^{-1}\ell$ is a normalized
right integral in $A$. By Maschke's theorem, $A$ is semisimple.
\end{proof}

Let $q$ qnd $q'$ be idempotents in $A_t$, and let
$\rho = \epsilon \actr q$ and  $\rho' = \epsilon \actr q'$.
Then $\rho$ and $\rho'$ are idempotents in $A_t ^*$.

\begin{lemma}
\label{identify} 
The space $\rho A^* S(\rho')$ may be
identified, in an $S^2$-invariant way,
with the dual space  of $qAS(q')$.  
\end{lemma}
\begin{proof}
It is easy to check that for all $h\in A$ and $\phi\in A^*$
one has $\la \phi,\, qhS(q') \ra = \la \rho \phi S(\rho'),\, h\ra$,
whence the statement follows.
\end{proof}

The next three corollaries refine \cite[Corollaries 6.5 - 6.7]{N}.

\begin{corollary}
\label{dual semisimplicity condition}
Let $A$ be a finite dimensional weak Hopf algebra. If for every
primitive idempotent $\pi$ in $A^*_t \cap A^*_s$ there exist 
idempotents $\rho_\pi,\, \rho_\pi'$ in $\pi A^*_t$ such that 
\begin{equation}
\label{Tr S2 neq 0 dual}
\Tr(S^2|_{\rho_\pi A^* S(\rho_\pi')}) \neq 0,
\end{equation}
then $A$ is semisimple. 
\end{corollary}
\begin{proof}
This immediately follows from Proposition~\ref{semisimplicity condition}
and Lemma~\ref{identify}, since the latter
implies $\Tr(S^2|_{\rho_\pi A^* S(\rho_\pi')}) = \Tr(S^2|_{q_pAS(q_p')})$.
\end{proof}

\begin{corollary}
\label{Tr S2 connected}
Let $A$ be a connected finite dimensional weak Hopf algebra.
If there exist idempotents $q,\, q'$ in $A_t$
such that $\Tr(S^2|_{qAS(q')})\neq 0$ then $A$ is semisimple.
\end{corollary}

\begin{corollary}
\label{Tr S2 biconnected}
Let $A$ be a biconnected finite dimensional weak Hopf algebra.
If there exist idempotents $q,\, q'$ in $A_t$
such that $\Tr(S^2|_{qAS(q')})\neq 0$ then $A$ and $A^*$
are semisimple.
\end{corollary}

\end{section}

\begin{section}{$\Tr(S^2)$ and the global dimension}

\subsection{Proof of Theorem \ref{traceofS2}}

Let $A$ be a connected 
semisimple weak Hopf algebra. Then
$\mathcal{C} = \Rep(A)$ is a fusion category. 
In particular, $V^{**}$ is isomorphic to $V$ for any simple 
object $V$, which means that $S^2$ is an inner 
automorphism of $A$ (this is also proved in [BNSz, 3.22]).
So let $g\in A$ be such that
$S^2(x) = gxg^{-1}$ for all $x\in A$. Then $g|_V : V \to V^{**} : v \mapsto 
g\cdot v$ is an isomorphism of representations.

The calculation of $\Tr_V(g)$ and $\Tr_V((g^{-1})^*)$ yields
\begin{eqnarray}
\label{Tr no 1}
\Tr_V(g){\rm Id}_{A_t} &=& \Tr(S(g)1\1|_{V^*})1\2,\\ 
\label{Tr no 2}
\Tr_{V^*}((g^{-1})^*){\rm Id}_{A_t}&=&\Tr(g^{-1}1\1|_{V})1\2.
\end{eqnarray}

\begin{remark} We note an important distinction between the quantum trace in 
an object $V$ of a morphism $V\to V^{**}$, which is denoted 
by $\Tr_V$, and the usual trace of a linear operator in the vector space $V$, 
which is denoted just by $\Tr$. \end{remark}

Let $\{q_i\}$ be primitive idempotents in $Z(A_t)$
and for every $i$ let $n_i^2 = \dim_k(q_iA_t)$ (i.e.,
$n_i$ is the dimension of the irreducible representation
of $A_t$ corresponding to $q_i$).
Multiplying equations \eqref{Tr no 1} and \eqref{Tr no 2}
on the right by $q_j$ and $q_i$ respectively, 
applying $\Tr|_{A_t}$, 
and using that for a regular weak Hopf algebra
$\Tr(z|_{A_t}) = \epsilon(z)$ for all $z\in A_t$ \cite{N},
we get 
\begin{equation*}
\Tr_V(g)n_j^2=\Tr(g|_{q_jV}),\quad \Tr_{V^*}((g^{-1})^*)n_i^2=\Tr(g^{-1}|_{S(q_i)V}),
\end{equation*}
and hence,
\begin{equation}
\Tr(g|_{q_jV})\Tr(g^{-1}|_{S(q_i)V})=|V|^2n_i^2n_j^2.
\end{equation} 
Taking the sum over all isomorphism classes of simple $V$, we obtain : 
\begin{equation}
\label{TrS2 vs global dimension refined}
\Tr(S^2|_{q_jAS(q_i)}) = \dim(\Rep(A)) n_i^2 n_j^2.
\end{equation}
This proves the first statement of the theorem. The second one follows from 
Theorem \ref{posi}. Theorem \ref{traceofS2}  is proved. 

\subsection{Proof of Theorem \ref{dualcat}}

Let $\mC$ be a fusion category, and let 
$\mM$ be an indecomposable module category over $\mC$. 
Let us show that $\mC^*_{\mM}$ is semisimple. 

Choose a semisimple algebra $R$ with blocks labeled by simple objects 
of $\mM$. Let $F$ be the fiber functor 
$F:\mC\to R-\text{bimod}$ associated with $\mM,R$ as explained in Section 2. 
Let $A=\End_k(F)$ be the corresponding weak Hopf algebra. 
Consider the weak Hopf algebra $A^{*cop}$. The representation category of 
$A^{*cop}$ is $\mC^*_{\mM}$. Thus, $A^{*cop}$ is connected 
(as $\mM$ is indecomposable). So, $A$ is biconnected. Also, 
by Theorem \ref{traceofS2}, we have 
$\Tr(S^2|_{q_jAS(q_i)})=\dim(\mC)n_i^2n_j^2\ne 0$. 
Thus, by Corollary \ref{Tr S2 biconnected}, $A^{*cop}$, and hence 
$\mC^*_{\mM}$ is semisimple. The equality 
$\dim(\mC)=\dim(\mC^*_{\mM})$ follows from 
Theorem \ref{traceofS2} and Lemma \ref{identify}
(noting that the block sizes corresponding to the idempotents
$q$ and $\epsilon\actr q$ are the same). Theorem \ref{dualcat} is proved. 

\subsection{Proof of Theorem \ref{m1m2}}

Let $\mM$ be a possibly decomposable module category over $\mC$. 
So $\mM=\oplus \mM_j$ is a direct sum of indecomposable subcategories.
Now for each $i$ choose semisimple algebras $R_i$ 
with blocks labeled by simple objects of $\mM_i$,
 and let $R=\oplus R_i$.
Let $F$ be the fiber functor $\mC\to R-\text{bimod}$,
associated to $(\mM,R)$, and let $B=\End_k(F)$ be the corresponding 
weak Hopf algebra. Then by Theorem \ref{traceofS2}
and Corollary~\ref{dual semisimplicity condition} applied to 
$A=B^*$, we find that $B^*$ is semisimple. 
In other words, the the category 
$\mC^*_{\mM}$ is semisimple. This implies Theorem \ref{m1m2}, since 
the category ${\rm Fun}_{\mC}(\mM_1,\mM_2)$ is contained in 
$\mC^*_{\mM_1\oplus \mM_2}$. Theorem \ref{m1m2} is proved. 

\subsection{A refinement of Theorem \ref{traceofS2}}

Let $A$ be a connected, but not necessarily coconnected,
 semisimple weak Hopf algebra.
Let $\pi_i$ be the primitive idempotents in $A^*_s\cap Z(A^*)$, and 
let $A_{ij}\subset A$ be defined by $A_{ij}^*=\pi_iA^*S(\pi_j)$. 
It is easy to see that $A_{ij}=p_jAp_i$, 
where $p_i$ are 
the primitive idempotents of $A_s\cap A_t$. 
In particular, $A_{ii}$ are semisimple algebras. 

In fact, for every $i$ the algebra $A_{ii}$ has a natural 
structure of a biconnected
weak Hopf algebra, and its representation category 
$\Rep(A_{ii})$ is equivalent to $\Rep(A)$ via the usual 
Morita equivalence $V\to p_i V$, $V\in \Rep(A)$. 
Indeed, this map is a tensor functor
since $\Delta(p_i) =(p_i\otimes p_i)\Delta(1)$.
To prove that it is an equivalence it suffices to show
that the central support $\hat p_i$ of $p_i$ equals $1$. 
 From the formula for $\Delta(p_i)$ above we have 
$\Delta(\hat p_iA) \subset \hat p_iA\otimes \hat p_iA$. But $\hat{p_i}A$ 
contains the matrix block $q_0A$ corresponding to the trivial
representation of $A$ (which is irreducible since $A$
is connected). Hence, $\Delta(q_0A)\subset \hat p_iA\otimes \hat p_iA$
and so $\hat p_i=1$. 

Let
$q_{ir}$ be the primitive idempotents of $Z(A_{rr})$.
Let $n_{ir}^2=\dim(q_{ir}(A_{rr})_t)$. 
We have the following refinement of Theorem \ref{traceofS2}.

\begin{theorem}\label{refine} One has
$\Tr(S^2|_{q_{jr}A_{rs}S(q_{is})})=\dim(\Rep(A))n_{jr}^2n_{is}^2$.
\end{theorem}

\begin{proof} The proof is completely analogous to 
the proof of Theorem \ref{traceofS2}. Namely, let 
$g$ be an element of $A$ such that $gxg^{-1}=S^2(x)$. 
Then $p_rgp_r$ realizes $S^2$ in $A_{rr}$. Hence, 
as was shown in the proof of Theorem \ref{traceofS2}
for any $V\in \Rep(A)$, one has
$$
\Tr(p_rgp_r|_{q_{ir}V})=\Tr_V(g)n_{ir}^2, 
$$
and 
$$
\Tr(p_rg^{-1}p_r|_{S(q_{ir})V})=\Tr_V((g^{-1})^*)n_{ir}^2. 
$$
Now, $A_{rs}=\oplus_{V}p_sV\otimes (p_rV)^*$, and 
$S^2|_{A_{rs}}=p_sgp_s\otimes (p_rg^{-1}p_r)^*$. 
Hence, taking the sum over all $V$, we arrive at the result. 
\end{proof}

\subsection{Proof of Proposition \ref{samedim}}
Let $\mC$ be an indecomposable
multi-fusion category. Fix an index 
$r$ and let $\mM=\oplus_i \mC_{ir}$. It is 
clear that $\mM$ is an indecomposable $\mC$-module 
category, and the category $\mC^*_\mM:={\rm Fun}_\mC(\mM,\mM)$
equals $\mC_{rr}^{op}$. 

Choose a semisimple algebra $R$
with blocks labeled by simple objects of $\mM$. 
Let $R_j$ be the two sided ideal in $R$ whose blocks 
are labeled by simple objects of $\oplus_i \mC_{ij}\subset \mM$. Let 
$n_j^2=\dim(R_j)$.  
Let $F$ be the corresponding 
fiber functor, and $B=\End_k(F)$ the corresponding weak
Hopf algebra, and let $A=B^*$. Then $A$ is connected and semisimple.
Therefore, by Theorem \ref{refine}, 
$\Tr(S^2|_{q_{jr}A_{rs}S(q_{is})})=\dim(\Rep(A))n_{is}^2n_{jr}^2$. 
 Thus, using Lemma \ref{identify}, 
we get 
\begin{eqnarray*}
\dim(\mC_{rs})n_{is}^2n_{jr}^2
&=& \Tr(S^2|_{(\epsilon\actr q_{jr})A_{rs}^*S(\epsilon\actr q_{is})})\\
&=& \Tr(S^2|_{q_{jr}A_{rs}S( q_{is})})\\
&=& \dim(\Rep(A))n_{is}^2n_{jr}^2,
\end{eqnarray*}
so $\dim(\mC_{rs})=\dim(\Rep(A))$, as desired.

\subsection{Proof of Theorem \ref{multi-fus}}
 
Without loss of generality it can be assumed that the 
multi-fusion category $\mC$ is indecomposable. 
Let $\mM_1$, $\mM_2$ be module categories over $\mC$. 
Let $\mathbf 1_i$ be 
the simple constituents of the unit object of $\mC$.
Then it is shown in the standard way that 
for each $i$, the restriction functor 
 ${\rm Fun_{\mC}}(\mM_1,\mM_2)\to {\rm Fun_{\mC_i}}(\mathbf 1_i\mM_1,
\mathbf 1_i\mM_2)$ is an equivalence of categories
(one can construct a quasi-inverse functor by extending 
a functor from $\mathbf 1_i\mM_1$ to $\mathbf 1_j\mM_1$ by tensoring with
objects of $\mC_{ij}$ and checking various compatibilities). 
Therefore, the result follows from Theorem \ref{m1m2}.

\subsection{Surjective and injective functors}

Let $\mC,\mD$ be multi-fusion categories, and 
$F:\mC\to \mD$ be a unital tensor functor. 

We will say that $F$ is injective if it is fully faithful
(in this case $F$ identifies $\mC$ with a full tensor subcategory of $\mD$). 
We will say that $F$ is surjective if any simple object $X$ 
of $\mD$ is contained in an object of the form $F(T)$, $T\in \mC$. 
If $F$ is surjective and $\mC$,$\mD$ are indecomposable, 
 then $\mD$ is said to be a quotient category of $\mC$. 

Now let $\mM$ be a faithful module category over $\mD$. Then it is also a 
faithful module category 
over $\mC$. Let $R$ be a semisimple algebra with blocks labeled by 
simple objects of $\mM$, and $G$ the corresponding 
fiber functor $\mD\to R$-bimod. Let $A=\End_k(G)$, $B=\End_k(G\circ F)$ 
be the corresponding weak Hopf algebras. Then the functor $F$ induces 
a morphism of weak Hopf algebras $\phi_F:A\to B$, and vice versa. 
It is easy to check that $F$ is surjective, respectively injective
if and only if $\phi_F$ is injective, respectively surjective. 

Let us now consider the corresponding dual categories $\mC^*,\mD^*$
(we will drop the subscript $\mM$, since the only module 
category we will consider in this subsection will be $\mM$). 
We have a natural dual functor $F^*:\mD^*\to \mC^*$, since any 
$\mD$-linear functor $\Phi:\mM\to \mM$ has a natural $\mC$-linear structure, 
defined using $F$. It is clear that the morphism $\phi_{F^*}:B^*\to A^*$ 
is given by $\phi_{F^*}=\phi_F^*$. 

\begin{proposition}\label{star} 
Let $F:\mC\to \mD$ be a unital tensor functor between multi-fusion categories. 
If $F$ is surjective (injective) then $F^*$ is injective (surjective). 
\end{proposition}
\begin{proof} The claim is equivalent to the statement that 
surjectivity of $\phi_F$ is equivalent to injectivity of $\phi_F^*$, and 
vice versa, which is obvious. 
\end{proof}

\subsection{The induction functor and the class equation for fusion categories}
\label{indfu}

In finite group theory, an important role is played 
by the ``class equation'' $\sum_{s}\frac{1}{|G_s|}=1$, where 
the summation is over conjugacy classes of a finite group $G$, and 
$G_s$ are the corresponding centralizers. There is a similar
formula for semisimple Hopf algebras (due to G.I.Kac and Y.Zhu,
see \cite{Ka}, \cite{Zhu}), 
which is also very useful. 
Here we will prove a class equation for any fusion category,
which generalizes the results of Kac and Zhu. 

Let $\mC$ be a fusion category. 
Let $I:\mC\to Z(\mC)$ be the induction functor, 
defined by the condition $\Hom_{Z(\mC)}(I(V),X)=\Hom_{\mC}(V,X)$,
$V\in \mC$, $X\in Z(\mC)$. 

Let $[X:Y]$ denote the multiplicity of a simple object $Y$ in an
object $X$. For $X\in Z(\mC)$, let $X|_\mC$ denote the
corresponding object of $\mC$. 
Then $I(V)=\oplus_{X\in Irr(Z(\mC))}[X|_{\mC}:V]X$. 

\begin{proposition}\label{I(1)}
One has $I(V)|_\mC=\oplus_{Y\in
  Irr(\mC)}Y\otimes V\otimes Y^*$. 
\end{proposition}

\begin{example} If $G$ is a finite group and $\mC=\Rep(G)$, 
then $I(\mathbf 1)$ is the following representation of the quantum double
$D(G)=\Bbb C[G]\ltimes \mathcal O(G)$: it is the regular
representation of the function algebra
$\mathcal O(G)$, on which $G$ acts by conjugation.
In this case proposition \ref{I(1)} for trivial $V$ is the standard fact that 
for conjugation action, $\mathcal O(G)=\oplus_{Y\in
  Irrep(G)}Y\otimes Y^*$. 
\end{example}

\begin{proof}
Let $\mC=\Rep(A)$, where $A$ is a weak
Hopf algebra. Then $Z(\mC)=\Rep(D(A))$, where $D(A)$ is the
quantum double of $A$ (see \cite{NTV}). As a space, 
$D(A)=A^*\otimes_{A_{min}}A$. So we have 
$I(V)|_\mC=D(A)\otimes_A V=A^*\otimes_{A_{min}}V$, with the
action of $A$ defined by 
$h(\phi\otimes v) = h\3 \actl \phi \actr S(h\1) \otimes h\2 v$
for all $h\in A,\, \phi\in A^*,\, v\in V$.
This yields exactly 
$\oplus_{Y\in Irr(\mC)}Y\otimes V\otimes Y^*$, 
where the tensor product is
taken in the representation category of $A$. The proposition is
proved. 
\end{proof}

\begin{proposition}\label{numobj} The number 
$\sum_{X\in Irr(Z(\mC))}[X|_{\mC}:\mathbf 1]^2$ 
is equal to the number of simple objects of $\mC$.
\end{proposition}
 
\begin{proof} 
This sum equals $\dim{\rm End}(I(\mathbf
1))$, which is the same as $[I(\mathbf 1)|_{\mC}:\mathbf 1]$. But
by Proposition \ref{I(1)} we
know that $I(\mathbf 1)|_{\mC}=\oplus_{Y\in Irr(\mC)}Y\otimes Y^*$, 
so the result follows.
\end{proof}

\begin{proposition} (The class equation)\label{cleq}
Let $\mC$ be a spherical fusion category. Then one has 
$$
\sum_{X\in Irr(Z(\mC)):[X|_{\mC}:\mathbf 1]\ne
  0}[X|_{\mC}:\mathbf 1]\frac{1}{m_X}=1,
$$
where $m_X=\frac{\dim \mC}{\dim X}$ is an algebraic integer. 
\end{proposition}

\begin{proof}
The proposition is immediately obtained 
by computing the dimension of $I(\mathbf 1)$ in two ways. 
On the one hand, this dimension is equal to $\dim \mC$ by 
Proposition \ref{I(1)}. On the other hand, it equals $\sum
[X|_\mC:\mathbf 1]\dim X$. The fact that $m_X$ is an algebraic
integer follows from Proposition \ref{moddivi}.
The proposition is proved. 
\end{proof}

\begin{example} Let $\mC=\Rep(G)$.
In this case $Z(\mC)=\Rep D(G)$, so simple objects 
are parametrized by pairs 
(conjugacy class, irreducible representation of the
centralizer). The simple objects $X\in Z(\mC)$ for which
$X|_\mC$ contains $\mathbf 1$ are exactly those for which 
the representation of the centralizer is trivial. 
Thus, the equation of Proposition \ref{cleq} is exactly 
the class equation for $G$. 
\end{example}

\end{section}

\subsection{Graded categories}

Let $G$ be a finite group. 

\begin{definition} A $G$-grading on a tensor category 
$\mC$ is a decomposition $\mC=\oplus_{g\in G}\mC_g$, such that 
$\otimes: \mC_g\times \mC_h\to \mC_{gh}$, 
$\mathbf 1\in \mC_e$, and $*:\mC_g\to
\mC_{g^{-1}}$. We will say that the grading is faithful if 
$\mC_g\ne 0$ for any $g\in G$.
\end{definition}

It is clear that if $\mC$ is a $G$-graded category, then 
$\mC_e$ is a tensor category, and $\mC_g$ are left 
and right module categories over $\mC_e$. 

Let $\mC$ be a fusion category, $G$ a finite abelian group. 

\begin{proposition}\label{gradi} Suppose that 
$Z(\mC)$ contains, as a full tensor subcategory,
the category of modules $\tilde G$ over 
$\mathcal O(G)$, 
and that all simple objects of this subcategory 
map to the neutral object of $\mC$ under 
the functor $Z(\mC)\to \mC$. 
Then $\mC$ admits a faithful $G^\vee$-grading
(where $G^\vee$ is the dual group to $G$).
\end{proposition}

\begin{proof} Let $g\in \tilde G$ be a simple object. 
Since $g$ is trivial as an object of $\mC$, 
it defines a tensor automorphism 
of the identity functor of $\mC$. 
Thus we have a homomorphism 
$G\to {\rm Aut}_{\otimes}(Id)$, and
hence any simple object of $\mC$ defines 
an element of the dual group $G^\vee$.
This is obviously a faithful grading. 
\end{proof}

\subsection{Deformations of weak Hopf algebras} 

It is well known that a non-semisimple finite dimensional algebra 
can be perturbed into a semisimple one. However, 
for (connected) weak Hopf algebras the situation is different. Namely, we have the 
following result. 

\begin{proposition} Let $H_t$ be a continuous family 
of finite dimensional regular connected weak Hopf algebras over $\Bbb C$ 
depending of a real parameter $t\in (a,b)$
(i.e., $H_t$ is independent of $t$ as a vectors space, and 
the structure maps continuously depend on $t$). 
Then either $H_t$ is semisimple for all $t$ 
or $H_t$ is non-semisimple for all $t$.
\end{proposition}

\begin{proof} 
The set $U\subset (a,b)$ where $H_t$ is semisimple is open.
Assume that this set is nonempty. 
The function $f(t)=\Tr(S^2|_{H_t})$ is a continuous function of $t$. 
By Corollary \ref{Tr S2 connected}, $f$ is zero
on the complement $U^c$ of $U$. Also, by Theorem \ref{posi} 
and Remark \ref{tra}, $f$ is a nonzero 
 algebraic integer for any $t\in U$. 
Since $f$ is continuous, this implies that 
$f$ is constant (and nonzero) on each connected component of $U$. 
Thus, $U^c$ is empty. We are done.  
\end{proof}

\subsection{Sphericity of pivotalization}

In this Section we allow the ground field to have a positive
characteristic, different from $2$.

Let $A$ be a quasitriangular weak Hopf algebra, i.e., such that
$\Rep(A)$ is a braided monoidal category. Let $\R\in (A\otimes A)\Delta(1)$ 
be the $R$-matrix of $A$.
It was shown in \cite{NTV} that the Drinfeld element $u = S(\R\II)\R\I$ is invertible,
satisfies  $S^2(h)=uhu^{-1}$ for all $h$ in $A$, and
\begin{equation}
\label{delta u}
\Delta(u^{-1}) = \R_{21}\R (u^{-1}\otimes u^{-1}).
\end{equation}
Furthermore, the element $uS(u)^{-1}$ is group-like and implements $S^4$ via the
adjoint action.

Observe that for any group-like element $\eta\in G(A^*)$ the element  
$g_\eta  = \R\I\la \eta,\, \R\II\ra$ is group-like in $A$ and the map
$\eta\mapsto g_\eta$ is a group homomorphism from $G(A^*)$ to $G(A^{op})$.

Recall that distingusihed elments of $A$ and $A^*$ were defined in 
Section~\ref{subsection about integrals}. In the next Lemma we extend
a result of Radford \cite{R2} relating the distinguished elements and $uS(u)^{-1}$. 

\begin{lemma}
\label{uS(u)-1 vs a}
Let $a\in A$ and $\alpha\in A^*$ be a pair of distinguished group-like elements
chosen as in Section~\ref{subsection about integrals}. Then $g_\alpha a = S(u)u^{-1}$.
\end{lemma}
\begin{proof}
Let $\ell\in A$ be a non-degenerate integral from which $a$ and $\alpha$ are constructed.
It follows from \cite[Section 5]{N} that 
\begin{eqnarray*}
\ell\2 \otimes \ell\1 &=& \ell\1 \otimes S^2(\ell\2)a^{-1}\\
\ell\1 h \otimes \ell\2 &=& \ell\1 \otimes \ell\2 S(h\actr \alpha),\qquad h\in A.
\end{eqnarray*}
Using these identities and properties of integrals and $R$-matrices we compute
\begin{eqnarray*}
\ell\2 S(\R\II \actr\alpha) \R\I \otimes \ell\1 
&=&  \ell\2  \R\I \otimes \ell\1 \R\II \\
&=&  \R\I \ell\1 \otimes \R\II \ell\2 \\
&=& \R\I  S(\R\II) \ell\1 \otimes \ell\2 \\
&=& \R\I  S(\R\II)  S^2(\ell\2)a^{-1} \otimes \ell\1, 
\end{eqnarray*}
whence $S(\R\II \actr\alpha) \R\I = \R\I  S(\R\II) a^{-1} = S(u)a^{-1}$ 
by  non-degeneracy of $\ell$. On the other hand,
from the axioms of $R$-matrix we have  $S(\R\II \actr\alpha) \R\I= S(\R\II) \R\I g_\alpha
=  u g_\alpha$, which completes the proof. 
\end{proof}

Recall \cite[Section 7]{NTV} that a central $\nu\in A$ is a ribbon
element if $\nu=S(\nu)$ and $\Delta(\nu) =
\R_{21}\R (\nu\otimes \nu)$. A quasitriangular weak Hopf algebra $A$
has a ribbon element if and only if $\Rep(A)$ is a ribbon category.

\begin{corollary}
\label{ribonness}
Suppose that $A$ is a unimodular quasitriangular weak Hopf algebra and that
$\alpha$ is gauged to be $\epsilon$. If $g\in A$ is a group-like element such that
$g^2 = a^{-1}$ and $S^2(h)=ghg^{-1}$ for all $h\in A$ then $\nu = u^{-1}g$ is
a ribbon element of $A$.
\end{corollary}
\begin{proof}
In this case we have $g_\alpha= 1$.
Since both $u$ and $g$ implement $S^2$, the element $\nu$ is central.
It follows from Equation~\eqref{delta u}  that $\Delta(\nu) =
\R_{21}\R (\nu\otimes \nu)$. Finally, $S(\nu) = S(u)^{-1}g^{-1} = u^{-1}g = \nu$
by Lemma~\ref{uS(u)-1 vs a}.
\end{proof}

Recall that in Remark \ref{pivotalization},
it was shown that to every fusion category 
$\mC$ one can associate a twice bigger 
fusion category $\tilde\mC$, which has a canonical pivotal
structure. 

\begin{proposition}\label{sphpiv} If $\dim\mC\ne 0$ then 
the canonical pivotal stucture of the category $\tilde \mC$  
is in fact a spherical structure. 
\end{proposition}
\begin{proof}
Let us realize $\tilde\mC$ as the representation category of a
regular weak Hopf algebra $A$. In this case, 
the pivotal structure of $\tilde\mC$ is defined 
by a grouplike element $g\in A$, such that $g^2=a^{-1}$, 
where $a$ is the distinguished group-like element of $A$.
Assume that the distinguished group-like element $\alpha\in
A^*$ is gauged to be $\epsilon$.
Consider the Drinfeld center $Z(\tilde \mC)$, i.e., the
representation category of the Drinfeld double $D(A)$. 
One can check that $a$ is also the distinguished 
grouplike element for $D(A)$.
Since $g$ is a grouplike element such that $gxg^{-1}=S^2(x)$
for $x\in D(A)$, and $g^2=a^{-1}$, it follows from 
Corollary~\ref{ribonness} that $D(A)$ is ribbon, hence,
$Z(\tilde\mC)$ is ribbon.

Thus, we have $\dim(V)=\dim(V^*)$ in $Z(\tilde\mC)$ \cite{K}. 
Hence, the same holds in $\tilde \mC$
(i.e., $\tilde \mC$ is spherical). 
Indeed, let $I: \tilde \mC\to Z(\tilde \mC)$ be the
induction functor defined in subsection \ref{indfu}. 
Then by Proposition \ref{I(1)}, for any $V\in \tilde\mC$, 
$I(V)=\oplus_{X\in Irr(\tilde\mC)} X\otimes
V\otimes X^*$ as an object of $\tilde\mC$. 
Hence, $\dim I(V)=\dim V\dim \tilde \mC$. On the other hand,
$I(V)^*=I(V^*)$, which implies the statement. 
\end{proof}

\begin{corollary}\label{posroot} In any pivotal fusion category $\mC$
over $\CC$, the dimension of a simple object
is a positive number times a root of unity. 
\end{corollary} 

\begin{proof} 
Consider the category $\tilde\mC$. It has two pivotal structures
$b_1,b_2: V\to V^{**}$. Namely, $b_1$ comes from 
$\mC$ and $b_2$ is the canonical pivotal structure. 
By Proposition \ref{sphpiv}, $b_2$ is spherical. 
On the other hand, $b_2^{-1}\circ b_1$ is a tensor 
endomorphism of the identity functor, 
which obviously acts by roots of unity on simple objects. 
This implies the statement. 
\end{proof}

\begin{section}{The co-Hochschild complex of a cosemisimple 
unimodular weak Hopf algebra and its subcomplex of invariants}

\subsection{The definition}

Let $A$ be a weak Hopf algebra and let $B \subset A$ be a weak Hopf
subalgebra. We will define a complex $C^\bullet(A, B)$
that generalizes the subcomplex of the  co-Hochschild complex 
of a Hopf algebra consisting of invariants of the adjoint 
action, cf.  \cite{S, Sch}.

We set $C^0(A, B) = A_t \cap Z(A)$, the algebra of endomorphisms
of the trivial $A$-module; and for all $n\geq 1$, we set
\begin{equation}
\label{Cn}
C^n(A, B) := \{ f\in A^{\otimes n} \mid f = \Ad(1)f \mbox{ and }
\Ad(h)f = \Ad(\eps_t(h))f,\, \forall h\in B \},
\end{equation}
where $\Ad(h)$ denotes the adjoint action of $h$, i.e.,
$\Ad(h)f = \Delta^n(h\1) f \Delta^n(S(h\2))$ for all $h\in B$
and $f\in A^{\otimes n}$. Here $\Delta^n : A \to A^{\otimes n}$
is the iteration of the comultiplication (we set $\Delta^1 =\id$). 
Note that $A^{\otimes n}$ is a non-unital $B$-module with respect 
to $\Ad$. The conditions of \eqref{Cn} mean that $A$ lies in
$A\otimes_{A_t}\cdots \otimes_{A_t} A$ (where $A$ is an $A_t$-bimodule
via $z_1\cdot h\cdot z_2 = z_1 S(z_2) h$ for  $h\in A,\, z_1,z_2, \in A_t$)
viewed as a subspace of $A^{\otimes n}$ and that $f$ commutes with 
$\Delta ^n(h)$ for all $h\in A$.

Let us define a linear map $\delta^n : C^n(A, B)\to A^{\otimes (n+1)}$ by
$\delta^0(f) = f - S(f)$ and
\begin{eqnarray}
\label{delta}
\delta^n(f)
&=& \Ad(1) \{1\otimes f +\sum_{i=1}^n\, (-1)^i (\id_{i-1}\otimes \Delta
\otimes \id_{n-i})(f) +(-1)^{n+1} f\otimes 1\}
\end{eqnarray}
for $n \geq 1$.

\begin{lemma}
For all $n\geq 0$
we have $\text{Im}(\delta^n) \subset C^{n+1}(A, B)$ and $\delta^{n+1}\circ \delta^n =0$,
i.e., the collection $\{C^n(A, B),\,\delta^n\}_{n\geq 0}$ is a complex.
\end{lemma}
\begin{proof}
That $\delta^n(f),\, f\in C^n(A, B)$ is $\Ad(1)$ stable is a part of definition
\eqref{delta}. For all $h\in B$ we compute
\begin{eqnarray*}
\lefteqn{\Ad(h)\delta^n(f) =} \\
&=& \Ad(h)(1\otimes f) - \Ad(h)(\Delta\otimes \id_{n-1})f + \cdots \\
& & \cdots + (-1)^n\Ad(h) (\id_{n-1}\otimes \Delta)f 
    + (-1)^{n+1}\Ad(h)(f\otimes 1) \\
&=& h\1 S(h\3) \otimes \Ad(h\2)f - (\Delta\otimes\id_{n-1})(\Ad(h)f) +\cdots\\
& & \cdots + (-1)^n (\id_{n-1}\otimes \Delta)\Ad(h)f
    + (-1)^{n+1} \Delta^n(h\1) f \Delta^n(S(h\3))\otimes \eps_t(h\2)\\
&=& 1\1\eps_t(h)\otimes \Ad(1\2)f - (\Delta\otimes\id_{n-1})
    (\Ad(\eps_t(h))f) +\cdots\\
& &  \cdots + (-1)^n (\id_{n-1}\otimes \Delta)\Ad(\eps_t(h))f 
    + (-1)^{n+1} \Delta^n(1\1)\Ad(\eps_t(h))f \otimes 1\2 \\
&=& \Ad(\eps_t(h))\delta^n(f),
\end{eqnarray*}
where we used the identity $(\id\otimes\eps_t)\Delta(x) =
(1\otimes x)\Delta(1)$ for all $x\in A$.
Therefore, $\delta^n(f)\in {C}^{n+1}(A, B)$.
Thus, $\delta^n(C^n(A, B)) \subset C^{n+1}(A, B)$.

That  $\delta^{n+1}\circ \delta^n =0$ can be checked in a direct
standard way, by cancelling the similar terms with opposite signs
among $(n+3)(n+2)$ tensors in $\delta^{n+1}\circ\delta^n(f)$.
\end{proof}

The $n$-th cohomology group of the pair $A,\, B$ is
\begin{equation}
\label{Hn}
H^n (A,\,B) := H^n(C^\bullet(A, B))= \mbox{Ker}(\delta^n)/\mbox{Im}(\delta^{n-1}), \qquad
\mbox{ for } n\geq 0.  
\end{equation}

\subsection{The vanishing theorem for cohomology}
 
 Now let $A$ be a regular finite-dimensional semisimple
(and hence, cosemisimple) weak Hopf algebra.

\begin{lemma}
\label{nondeg of chi}
Suppose that $A$ is indecomposable.
Let $q_i$ be a primitive idempotent in $Z(A_t)$
and let $\chi:=\chi_i$ be the corresponding canonical left integral 
in $A^*$ from Proposition~\ref{chi}.
Then $u = 1\1 \chi(1\2)$ is an invertible scalar.
\end{lemma}
\begin{proof} 
Since $\chi$ is a left integral in $A^*$, we have
\begin{equation}
u = 1\1 \chi(1\2) = S(1\1)\chi(1\2) = S(u) \in A_s\cap A_t.
\end{equation}
Next, the following identities hold in every weak Hopf algebra $A$:
\begin{equation*}
\Delta(1)(h\otimes 1) = (\id \otimes \eps_t)\Delta(h), \qquad
(h\otimes 1)\Delta(1) = (\id \otimes \eps_t^{{\text{op}}})\Delta(h),\qquad h\in A,
\end{equation*}
where  $\eps_t^{\text{op}}$ is the target counital map in $A^{\text{op}}$.

Using Proposition~\ref{chi} we compute for all $h\in A$ :
\begin{eqnarray*}
1\1 h \chi(1\2)
&=& h\1 \chi(\eps_t(h\2)) = h\1 \chi(h\2 S(h\3)) \\
&=& h\1 \chi (S^{-1}(h\3)h\2) = h\1\chi( \eps_t^{{\text{op}}}(h\2)) \\
&=& h 1\1 \chi(1\2),
\end{eqnarray*}
hence, $u\in A_s\cap A_t\cap Z(A) = k1_A$.

For any primitive central idempotent $q_j$ in $A_t$ we have
\begin{equation*}
\chi(q_j) = \Tr(S^2|_{q_jAS(q_i)}) =\dim(\Rep(A)) n_i^2 n_j^2 \neq 0,
\end{equation*}
where $n_i^2 = \dim_k(q_iA_t)$ and $n_j^2 =\dim_k(q_jA_t)$.
Therefore, $\chi|_{A_t}\neq 0$ and $u\neq 0$.
\end{proof}

\begin{theorem}
\label{acyclic}
 Let ${C}$ be the complex associated to $A,\,B$ as above.
Then  $H^0(A, B) = A_t\cap A_s\cap Z(A)$ and $H^n(A, B) =0$ for all $n\geq 1$.
\end{theorem}
\begin{proof}
It suffices to prove this theorem in the case when $A$ is indecomposable. In this
case we have $H^0(A, B) = \text{Ker}(\delta^0) = A_t\cap A_s\cap Z(A) =k$.
Let $f\in {C}^n(A, B),\, n\geq 1$ be a cycle, i.e., $\delta^n(f) =0$. 
We will show by a direct computation that 
\begin{equation}
\tilde{f} =(\id_{n-1}\otimes \chi)(f)
\end{equation}
belongs to  ${C}^{n-1}(A, B)$ and that 
$f$ is a scalar multiple of the differential $\delta^{n-1}(\tilde{f})$.

We will write $f=f^{(1)}\otimes\cdots\otimes f^{(n)}$, where a sum
of some set of tensors is understood.  We have :
\begin{eqnarray*}
\lefteqn{ \Ad(1)\tilde{f} = }\\
&=& 1\1 f^{(1)} S(1_{(2n-2)})\otimes\cdots\otimes   
    1_{(n-1)} f^{(n-1)} S(1_{(n)})\chi(f^{(n)}) \\
&=& \Delta^{n-1}(1\1 1'\1)(f^{(1)}\otimes\cdots\otimes f^{(n-1)}) \Delta^{n-1}(S(1\2 1'\4))
    \chi(1'\2 f^{(n)}S(1'\3)) \\
&=&  \Delta^{n-1}(1'\1)(f^{(1)}\otimes\cdots\otimes f^{(n-1)}) \Delta^{n-1}(S(1'\4))
    \chi(1'\2 f^{(n)}S(1'\3)) = \tilde{f},
\end{eqnarray*}
where $1'$ stands for a second copy of $1$ and we use that in any 
regular weak Hopf algebra one has $(1\1 \otimes 1\otimes 1\otimes 1\2)\Delta^4(1)
= \Delta^4(1)$. Next, we have for all $h\in B$ :
\begin{eqnarray*}
\lefteqn{ \Ad(h)\tilde{f} = }\\
&=& h\1 f^{(1)} S(h_{(2n-2)})\otimes\cdots\otimes   
    h_{(n-1)} f^{(n-1)} S(h_{(n)})\chi(f^{(n)}) \\
&=& h\1 f^{(1)} S(h_{(2n-1)})\otimes\cdots\otimes   
    h_{(n-1)}) f^{(n-1)} \eps_s(h_{(n)}) S(h_{(n+1)}) \chi(f^{(n)}) \\
&=& h\1 f^{(1)} S(h_{(2n-1)})\otimes\cdots\otimes   
    h_{(n-1)} f^{(n-1)} S(h_{(n+1)})\chi(f^{(n)}S(\eps_s(h_{(n)}))) \\
&=& h\1 f^{(1)} S(h_{(2n)})\otimes\cdots\otimes   
    h_{(n-1)} f^{(n-1)} S(h_{(n+2)})\chi(h_{(n)} f^{(n)}S(h_{(n+1)}))  \\
&=& (\id_{n-1}\otimes \chi)(\Ad(h)f)
= (\id_{n-1}\otimes \chi)(\Ad(\eps_t(h))f) = \Ad(\eps_t(h))\tilde{f}.
\end{eqnarray*}
Here we used that 
\begin{equation}
\label{Delta1}
f = (1\otimes\cdots\otimes 1\otimes \Delta(1))f
\end{equation}
for all $f\in {C}^n(A, B)$, the weak Hopf algebra identities
\begin{equation*}
x\1\eps_s(x\2)=x \quad \mbox{ and } \quad 
(y\otimes 1)\Delta(1) = (1\otimes S^{-1}(y))\Delta(1)
\end{equation*}
for all $x\in A$ and $y\in A_s$ (cf.\ \cite{NV}), and Proposition~\ref{chi}.

Hence, $\tilde{f}\in {C}^{n-1}(A, B)$.  Next, since $(\id_n\otimes \chi)(\delta^n(f))=0$, 
we get
\begin{eqnarray*}
0 
&=& (\id_n\otimes \chi)\Ad(1) \{ 
    1\otimes f^{(1)}\otimes\cdots\otimes f^{(n)} 
     - \Delta(f^{(1)})\otimes\cdots\otimes f^{(n)} + \\
& & \cdots + (-1)^n  f^{(1)}\otimes\cdots\otimes 
    \Delta(f^{(n)}) + (-1)^{n+1} f^{(1)}\otimes\cdots\otimes f^{(n)} \otimes 1\} \\
&=& \Ad(1) \{
    1\otimes f^{(1)}\otimes\cdots\otimes f^{(n-1)} \chi(f^{(n)}) 
    - \Delta(f^{(1)})\otimes\cdots\otimes f^{(n-1)} \chi(f^{(n)}) + \\
& & \cdots + (-1)^n  f^{(1)}\otimes\cdots\otimes (\id\otimes\chi)
    \Delta(f^{(n)})\} + (-1)^{n+1} f \Delta^n(1\1) \chi(1\2).
\end{eqnarray*} 
Let $u =1\1 \chi(1\2)$. Then 
\begin{eqnarray*}
f \,\Delta^n(u)
&=& (-1)^n \delta^{n-1}(\tilde{f}) \\ 
& &  + \Ad(1)   \{f^{(1)}\otimes\cdots\otimes (\id\otimes\chi)
    \Delta(f^{(n)}) - f^{(1)}\otimes\cdots\otimes f^{(n-1)}
     \otimes \chi(f^{(n)})1\}.
\end{eqnarray*}
Let us show that 
\begin{equation}
\label{it remains}
\Ad(1) \{ (\id_{n-1}\otimes \chi)(f) \otimes 1 - (\id_n \otimes \chi)
(\id_{n-1}\otimes\Delta)(f) \} =0
\end{equation}
for all $f\in C^n(A, B)$. Observe that both sides of \eqref{it remains}
belong to $C^{n+1}(A, B)$. We compute :
\begin{eqnarray*}
\lefteqn{ \Ad(1)\{f^{(1)}\otimes\cdots\otimes f^{(n-1)}
\otimes {f^{(n)}}\1 \chi({f^{(n)}}\2)\} = }\\
&=& f^{(1)}\otimes\cdots\otimes f^{(n-1)} \otimes S(1\1)\chi(1\2 f^{(n)}) \\ 
&=& f^{(1)}\otimes\cdots\otimes 1'\1 f^{(n-1)} \otimes 1'\2 S(1\1)
    \chi(1\2 f^{(n)})   \\
&=& f^{(1)}\otimes\cdots\otimes 1'\1 1\1 f^{(n-1)} \otimes 1'\2 
    \chi(1\2 f^{(n)})  \\
&=& f^{(1)}\otimes\cdots\otimes 1\1 f^{(n-1)} \otimes 1\2 
    \chi(f^{(n)}\\
&=&  \Ad(1)\{f^{(1)}\otimes\cdots\otimes f^{(n-1)} \otimes  \chi(f^{(n)})1 \}.
\end{eqnarray*}
Here $1'$ stands for another copy of $1$. We used the definition of 
a left integral, several identities in a weak Hopf algebra, and 
equation~\eqref{Delta1} applied to $f$ and $\tilde{f}$. 

Therefore, 
$f\Delta^n(u) = (-1)^n \delta^{n-1}(\tilde{f})$.
Now $u$ is a non-zero scalar by Lemma~\ref{nondeg of chi}, hence
$f = \delta^{n-1}((-1)^n u^{-1}\tilde{f})$,
i.e., $\text{Ker}\, \delta^n \subset \text{Im}\, \delta^{n-1}$
and  $H^n(A, B)=0$ for $n\geq 1$.
\end{proof}

\end{section}

\begin{section}
{The Yetter cohomology of a tensor category and Ocneanu rigidity}

\subsection{The definition}
The following cohomology of a tensor category 
with respect to a tensor functor was defined 
in \cite{Y1,Y2}, motivated by the previous work \cite{CY},
and independently in \cite{Da}. 

Let $\mC,\mC'$ be tensor categories over a field $k$, and 
$F:\mC\to \mC'$ be a (unital) tensor functor.
 
For a nonnegative integer $n$, let $\mC^n$ denote the 
$n$-th Cartesian power of $\mC$. 
In particular, $\mC^0$ has one object (empty tuple $\emptyset$) 
and one morphism (identity). Define the functor 
$T_n: \mC^n\to \mC$ by
$T_n(X_1,...,X_n):=X_1\otimes...\otimes X_n$. 
In particular, $T_0:\mC^0\to \mC$ is defined by 
$T_0(\emptyset)=\mathbf 1$,
and $T_1=\Id$. Let $C^n_F(\mC)=\text{End}(T_n\circ F^{\otimes n})$
(so e.g., $C^0_F(\mC)=\End(\mathbf 1_{\mC'})$). 
We define a differential $d: C^n(\mC)\to C^{n+1}(\mC)$ 
by the formula
$$
df=\Id\otimes f_{2,...,n+1}-f_{12,...,n+1}+f_{1,23,...,n+1}-...+
(-1)^nf_{1,...,nn+1}+(-1)^{n+1}f_{1,...,n}\otimes \Id,
$$
where for instance $f_{12,3,...,n+1}$ is the endomorphism
of the product of $n$ objects $F(X_1)\otimes F(X_2),F(X_3),...,F(X_{n+1})$, 
and we use the identification $F(X_1)\otimes F(X_2)\to F(X_1\otimes X_2)$ 
defined by the tensor structure on $F$. 
 
It is easy to show that $d^2=0$. 
Thus $(C^\bullet(\mC),d)$ is a complex. We will call  
the cohomology groups $H^n_F(\mC),n\ge 0$ of this complex 
the Yetter cohomology groups of $\mC$ with respect to $F$. 
In the important special case $\mC=\mC'$, $F=\Id$,
we will denote this cohomology simply by $H^i(\mC)$ and call it the 
Yetter cohomology of $\mC$. 

As usual, low dimensional Yetter cohomology groups 
have an independent meaning \cite{Y1,Y2,Da}. 
The group $H^1_F(\mC)$ classifies
derivations of $F$ as a tensor functor. 
The group $H^2_F(\mC)$ classifies first order deformations 
of the tensor structure on the functor $F$. 
The group $H^3(\mC)$ classifies first order defomations
of the associativity constraint in $\mC$, i.e. of the structure of $\mC$ 
itself. As usual, obstructions to these deformations live in the cohomology 
groups one degree higher. 

To illustrate this definition, we give a few examples. 

\begin{example} If $\mC$ is 
$\Rep(H)$ for a Hopf algebra $H$ and $F$ the forgetful functor, then 
$H^i_F(\mC)$ is the usual Hochschild cohomology of $H$ as a coalgebra. 
\end{example}

\begin{example}
Let $G$ be a finite group, and  
$\mC={\rm Vec}_G$ be the category 
of $G$-graded vector spaces over an algebraically closed 
field $k$ (of any characteristic). 
Then $H^i(\mC)=H^i(O(G))$, the Hochschild cohomology 
of $O(G)$ as a coalgebra, where $O(G)$ is 
the Hopf algebra of $k$-valued functions on $G$. In other words, 
$H^i(\mC)=H^i(G,k)$, the group 
cohomology of $G$ with coefficients in $k$. 
\end{example}

\begin{example}
Let $G$ be a reductive algebraic group over $k$ 
with Lie algebra $\g$.
Let $\Rep(G)$ denote the category of algebraic representations of 
$G$. Then $H^i(\Rep(G))=(\Lambda^i\g)^G$ for all $i$. 
Indeed, $\Rep(G)=\Rep(O(G)^*)$. 
Therefore, $C^n(\Rep(G))=(O(G^n)^*)^G$, where $G$ acts 
diagonally by conjugation. Since $G$ is reductive, 
the cohomology of this complex is the $G$-invariants 
in the Hochschild cohomology of $O(G)$ with coefficients 
in the trivial representation (corresponding to $1\in G$). 
This cohomology is well known to be $\Lambda\g$, so 
the answer is $(\Lambda\g)^G$, as desired. 
Thus, if $\g$ is a simple Lie algebra, then there is no nontrivial 
derivations or tensor structure deformations
of the identity functor, but there exists 
a unique (up to scaling) first order deformation 
of the associativity constraint, corresponding to a basis element in 
$(\Lambda^3\g)^G$. It is easy to guess that this deformation 
comes from an actual deformation, namely from the deformation of 
$O(G)$ to the quantum group $O_q(G)$. 
\end{example}

\subsection{Comparison of two complexes}

The following proposition shows that 
for categories and functors coming from weak Hopf algebras,  
the Yetter complex defined here is the same as the complex 
defined in the previous section.

\begin{proposition} \label{twocomp}
Let $A$ be a weak Hopf algebra, 
$B$ a weak Hopf subalgebra in $A$, 
$\mC=\Rep(A)$, $\mC'=\Rep(B)$, 
and $F:\mC\to \mC'$ be the restriction functor. 
Then  the complex $C^\bullet_F(\mC)$ coincides with the complex 
$ C^\bullet(A,B)$ defined in Section 6. 
\end{proposition}
\begin{proof} 
Let us show that $\text{End}(T_n\circ F^{\otimes n})= C^n(A,B)$.  
For this, define a map $\xi:  C^n(A,B)\to \text{End}(T_n)$. 
Recall that for any representations $X_1,...,X_n$ of $A$, 
the representation $X_1\otimes_\mC ...\otimes_\mC X_n$ 
is defined to be the image of the projection $\Delta^n(1)$ acting on 
$X_1\otimes...\otimes X_n$ \cite{NTV}. Therefore, any element 
$a\in  C^n(A,B)$ defines an operator on 
the usual tensor product $X_1\otimes...\otimes X_n$ which 
preserves the subspace $X_1\otimes_\mC ...\otimes_\mC X_n$
(since $a=\text{Ad}(1)a$). 
Let us denote by $\xi(a)$ its restriction to this subspace. 
Since $a$ is invariant under the adjoint action, we find that $\xi(a)$ is 
not only a linear map but also a morphism of representations. 
So we have defined the map $\xi$. 

Now we need to show that the map $\xi$ is an isomorphism.  
For this purpose, we construct the inverse map. 
Take $f\in \text{End}(T_n\circ F^{\otimes n})$. Extend it by zero 
to the space $\text{Ker} \Delta^n(1)$ 
in $X_1\otimes...\otimes X_n$. Then we get 
a linear operator $\tilde f$ on $X_1\otimes...\otimes X_n$,
which is functorial in $X_1,...,X_n$. Thus, $\tilde f$ corresponds 
to a unique element $\hat f$ of $A^{\otimes n}$. 
It is clear that $\hat f=\text{Ad}(1)\hat f$, since both sides define 
the same operators in $X_1\otimes...\otimes X_n$. 
It is also clear that $\hat f$ is invariant under adjoint action, as 
$f$ is a morphism in the category. Finally, it is clear that 
$\xi(\hat f)=f$, and $\widehat{\xi(a)}=a$. Thus, $\xi$ is an isomorphism.

Finally, it is easy to see that the 
map $\xi$ identifies the differentials of the two complexes. 
The proposition is proved. 
\end{proof}

\subsection{Proof of Theorem \ref{van} and Theorem \ref{multocnrig}}

Let $\mC,\mC'$ be multifusion categories, and $F:\mC\to \mC'$ a tensor 
functor between them. We can assume without loss of generality that $F$ is 
faithful and surjective. 
Let $R$ be a semisimple algebra, and let $G:\mC'\to R-\text{bimod}$
be a fiber functor on $\mC'$ (we know by now that it necessarily exists for a 
suitable $R$). Let $B=\End(G)$ be the corresponding weak Hopf algebra. 
Let $A=\End(G\circ F)$. Then $B$ is a weak Hopf subalgebra of $A$, 
$\mC'=\Rep(B)$, $\mC=\Rep(A)$, and $F$ is the restriction functor. 
Thus, by Proposition \ref{twocomp}, the complex $C^\bullet_F(\mC)$ is 
isomorphic to the complex $C^\bullet(A,B)$. Now, since 
$A$ is semisimple, it is also cosemisimple by Theorem \ref{sscoss}.
Hence, by Theorem \ref{acyclic}, $C^\bullet_F(\mC)$ is 
acyclic in positive degrees. Theorem \ref{van} is proved.  

Now let us prove Theorem \ref{multocnrig}.
Let $\mC$ be a multifusion category. 
By theorem \ref{van}, we have $H^i(\mC)=0$, $i>0$, in particular 
$H^3(\mC)=0$. 

Now consider the set of all admissible associativity constraints
for a tensor category with the same Grothendieck ring as $\mC$. 
It is an affine algebraic variety. This variety is acted upon 
by the group of twists. 
By a standard argument, one finds that 
the quotient of the tangent space to this variety 
at the point $\mC$ by the tangent space to 
the orbit of $\mC$ is equal to  $H^3(\mC)$.
Since $H^3(\mC)=0$, the group action is locally transitive, 
and hence has finitely many orbits. We are done. 

In conclusion we give the formulation of Ocneanu rigidity for
semisimple weak Hopf algebras. This is a generalization 
to the weak case of Stefan's rigidity theorem for semisimple 
Hopf algebras  \cite{S}, which was conjectured in \cite[3.9]{N}.

\begin{theorem}\label{rigweak} A regular semisimple 
finite dimensional weak Hopf algebra 
does not admit nontrivial 
 deformations in the class of regular weak Hopf algebras. 
In particular, there are finitely many such 
weak Hopf algebras in each dimension. 
\end{theorem}

\begin{proof}
Let $A$ be a regular semisimple weak Hopf algebra, 
$\Rep(A)$ its representation category, 
$R$ its base, and let $F:\Rep(A)\to 
R-\text{bimod}$ be  the forgetful functor.  
Suppose we have a deformation of $A$. 
Then by Theorem \ref{multocnrig}, this deformation 
produces a trivial deformation of 
$\Rep(A)$ as a tensor category. It is also clear that the base 
$R$ of $A$, being a semisimple algebra, is deformed trivially. 
This means that this deformation comes from a deformation
of the forgetful functor $F$ as a tensor functor
(indeed, a regular weak Hopf algebra is completely determined 
by the associated forgetful functor, since the symmetric separability 
idempotent is unique). 
But the deformation of $F$ must also be trivial,
 by Theorem \ref{ocnrigfun}. This implies the first statement of 
the theorem.  The second statement follows from the first one 
as in the proof of Theorem \ref{multocnrig}.  
\end{proof}

\end{section}

\section{Frobenius-Perron dimensions}

In this section we will give some applications of the
Frobenius-Perron theorem to the theory of fusion categories. 
We work with categories over $\CC$. 

\subsection{Definition of the Frobenius-Perron dimensions}

We start by recalling the Frobenius-Perron theorem (see
\cite{Ga}).

\begin{theorem} Let $A$ be a square matrix with nonnegative entries. 
\begin{enumerate}
\item[(1)] $A$ has a nonnegative real eigenvalue. The largest
nonnegative real eigenvalue $\lambda(A)$ of $A$ dominates
absolute values of all other eigenvalues of $A$. 
\item[(2)] If $A$ has strictly positive entries then $\lambda(A)$ is a
simple positive eigenvalue, and the corresponding eigenvector can be
normalized to have strictly positive entries. 
\item[(3)] If $A$ has an eigenvector $\mathbf f$ with strictly positive
entries, then the corresponding eigenvalue is $\lambda(A)$. 
\end{enumerate}
\end{theorem}

Now let $A$ be a unital based ring of finite rank in the sense of
Lusztig (see, 
e.g., \cite{O} for a definition). Let $b_i$ be the basis of $A$ ($b_0=1$), 
and $[b_i]$ be the matrix of multiplication by $b_i\in A$ in the basis
$b_j$. This matrix has nonnegative entries. 
Let $\lambda_i$ be the largest nonnegative eigenvalue
of $[b_i]$. Since $b_i$ is clearly not nilpotent, 
this eigenvalue is actually positive. 

We will call $\lambda_i$
the {\it Frobenius-Perron dimension} of $b_i$. 
We note that Frobenius-Perron dimensions 
for commutative based rings were defined 
and used in the book \cite{FK}.

\begin{theorem}\label{fp0}
The assignment $b_i\to \lambda_i$ extends to a homomorphism of
algebras $A\to \mathbb R$. 
\end{theorem}
\begin{proof}
For any $j$ define 
$z_j=\sum_i b_ib_jb_i^*$. We claim that $z_j$ is central. Indeed,
let $b_ib_j=\sum_k N_{ij}^kb_k$. 
Then $N_{ri}^k=N_{k^*r}^{i^*}$, so
$$
b_rz_j=b_r\sum_i b_ib_jb_i^*=\sum_{k,i}N_{ri}^kb_kb_jb_i^*=
\sum_{k,i}N_{k^*r}^{i^*}b_kb_jb_{i^*}=
\sum_k b_kb_jb_{k^*}b_r=z_jb_r.
$$
Let $z=\sum_j z_j$. Since 
$z_j$ contain $b_j$ as a summand, the matrix $[z]$ of multiplication by $z=
\sum_j z_j$ has strictly positive entries. 
Let $\mathbf f$ be the unique eigenvector of $[z]$ 
with positive entries and $f_0=1$ (its existence and uniqueness
follows from the Frobenius-Perron theorem). 
Then $[b_j]\mathbf f$ must be a multiple of $\mathbf f$, since 
$\mathbf f$ has positive entries and $[z][b_j]=[b_j][z]$. So by the
Frobenius-Perron theorem, $[b_j]\mathbf f = \lambda_j\mathbf f$, 
whence the statement follows.   
\end{proof}

We note that $b_i\to \lambda_i$ is the unique homomorphism 
of $A$ to $\mathbb R$ such that images of the basis 
vectors are positive. This is demonstrated by the following simple and 
probably well known lemma. 

Let $A$ be a $\mathbb Z_+$-ring 
of finite rank with $\mathbb Z_+$-basis $b_i$ (see \cite{O} for definition). 

\begin{lemma}\label{uniq} If there exists a homomorphism $\phi:A\to \mathbb R$ 
such that $\phi(b_i)>0$ for all $i$, then $\phi$ is unique. 
\end{lemma}

\begin{proof} Let $[b_i]$ be the matrix of multiplication by $b_i$ in the basis $b_j$, 
and $\mathbf f$ be the column vector with entries $\phi(b_j)$. 
Then $[b_i]\mathbf f=\phi(b_i)\mathbf f$. But $[b_i]$ is a matrix with nonnegative entries, 
and $\mathbf f$ has positive entries. Thus,
by the Frobenius-Perron theorem, $\phi(b_i)$ is the largest positive 
eigenvalue of $[b_i]$, and hence is uniquely determined. Lemma is proved. 
\end{proof}

\begin{remark} 
M.~M\"uger pointed out to us that for any $i$ one has
$\lambda_i\ge 1$; moreover if $\lambda_i<2$ then $\lambda_i=
2\cos(\pi/n)$ for some integer $n\ge 3$ (this follows from the well known
Kronecker theorem on the eigenvalues of positive integer matrices).
\end{remark}

One can also define the Frobenius-Perron dimensions of 
basis elements of a based module over a unital based ring. 
Indeed, let $M$ be an indecomposable based $A$-module with basis $m_j$ (see
\cite{O} for definition). 

\begin{proposition}\label{frobmod0} There exists a unique up to scaling 
common eigenvector $\mathbf m$ of the matrices $[b_i]|_M$ of
the operators $b_i$ in $M$,
which has strictly positive entries. 
The corresponding eigenvalues are $\lambda_i$. 
\end{proposition} 

\begin{proof} Let $z$ be the element from the proof of Theorem
  \ref{fp0}. 
The matrix $[z]|_M$ has strictly positive entries. 
Let $\mathbf m$ be an eigenvector of $[z]|_M$ with positive
entries. It exists and is unique up to a positive scalar, 
by the Frobenius-Perron theorem. Hence, it is an eigenvector
of $[b_i]|_M$, and the eigenvalues are $\lambda_i$ 
by Lemma \ref{uniq}.
\end{proof}

Now let 
$\mC$ be a fusion category. 
Let $K(\mC)$ be the Grothendieck ring of $\mC$. 
For any object $X\in \mC$, define the Frobenius-Perron dimension
of $X$, ${\rm FPdim}(X)$, to be the largest positive eigenvalue of 
the matrix $[X]$ of multiplication by $X$ in $K(\mC)$. 
Since $K(\mC)$
is a unital based ring of finite rank, 
Theorem \ref{fp0} implies the following result. 

\begin{theorem}\label{fp} The assignment $X\to {\rm FPdim}(X)$
  extends to a 
homomorphism of algebras 
$K(\mC)\to \mathbb R$. 
\end{theorem}

Similarly, one can define Frobenius-Perron dimensions of simple objects 
of a module category over a fusion category. Namely, let $\mM$ be an 
indecomposable module category over a fusion category $\mC$. 
For any $X\in \mC$ let $[X]_\mM$ be the matrix 
of action of $X$ on the Grothendieck group of $\mM$. 
Since the Grothendieck group $K(\mM)$ is an indecomposable  based module over 
$K(\mC)$, Proposition \ref{frobmod0} implies
 
\begin{proposition}\label{frobmod} There exists a unique up to scaling 
common eigenvector $\mathbf m$ of the matrices $[X]_\mM$,
which has strictly positive entries. 
The corresponding eigenvalues are ${\rm FPdim}(X)$. 
\end{proposition} 

The entries $\mu_i$ of $\mathbf m$ are called the Frobenius-Perron dimensions 
of simple objects $M_i$ of $\mM$. Unlike Frobenius-Perron dimensions 
of objects of $\mC$, they are well defined only up to scaling. 

\subsection{Frobenius-Perron dimension of a fusion category}

Let $\mC$ be a fusion category. Define the Frobenius-Perron
dimension ${\rm FPdim}(\mC)$ of
$\mC$ to be the sum of squares of Frobenius-Perron dimensions of simple
objects of $\mC$. 

Also, define the regular representation $R_\mC$ to be the
(virtual) object \linebreak $\oplus_{X\in Irr(\mC)} {\rm FPdim}(X)X$, 
where $Irr(\mC)$ is the set of isomorphism classes of simple
objects of $\mC$. It is clear that ${\rm FPdim}(\mC)={\rm
  FPdim}(R_\mC)$. 

\begin{proposition} \label{freeness}
Let $F:\mC\to \mD$ be a surjective tensor functor between fusion categories.
Then $F(R_\mC)=\frac{{\rm FPdim}(\mC)}{{\rm FPdim}(\mD)}R_\mD$. 
\end{proposition}

\begin{proof} Let $X$ be an object of $\mC$ such that $F(X)$
contains all simple objects of $\mD$. It exists since $F$ is
surjective. The vectors of multiplicities of 
$F(R_\mC)$ and $R_\mD$ are eigenvectors of the matrix of multiplication
by $F(X)$. The matrix and the
vectors have strictly
positive entries. By the Frobenius-Perron theorem, this implies
that $F(R_\mC)$ and $R_\mD$ differ by a positive scalar. 
This scalar is readily found to be the claimed one by taking
${\rm FPdim}$ of $F(R_\mC)$ and $R_\mD$.
\end{proof}

\begin{corollary}\label{freeness1} Let $A\subset B$ be an inclusion of finite
dimensional semisimple quasi-Hopf algebras. Then $B$ is free as a
left $A$-module. 
\end{corollary}

\begin{proof} In this case Frobenius-Perron dimensions coincide
  with the usual ones, so if $\mC=\Rep(B)$, $\mD=\Rep(A)$ then 
$R_\mC=B$, $R_\mD=A$, proving the claim.
\end{proof}

\begin{remark} For Hopf algebras, Corollary \ref{freeness1} is the well known freeness
theorem of Nichols and Zoeller \cite{NZ} (it is true in
the nonsemisimple case as well). In the quasi-Hopf case 
it was independently obtained by P. Schauenburg (without the
semisimplicity assumption), \cite{Scha}.
\end{remark}

\begin{corollary}\label{divis} If $F:\mC\to \mD$ is surjective then ${\rm
    FPdim}(\mC)$ is divisible by ${\rm FPdim(\mD)}$ 
(i.e. the ratio is an algebraic integer). 
\end{corollary}

\begin{proof} Taking the multiplicity of the neutral object of
  $\mD$ in Proposition \ref{freeness}, 
we find $\sum_{X\in Irr(\mC)}{\rm FPdim}(X)[F(X):\mathbf 1]=
\frac{{\rm FPdim}(\mC)}{{\rm FPdim}(\mD)}$, which yields the
statement. 
\end{proof} 

\begin{proposition}\label{squar} For any fusion category $\mC$, one has 
${\rm FPdim}(Z(\mC))={\rm FPdim}(\mC)^2$. 
\end{proposition}

\begin{remark} If $\mC$ is a representation category of a 
quasi-Hopf algebra $A$, then Proposition \ref{squar} is 
saying that the Frobenius-Perron dimension of $Z(\mC)$ is
$\dim(A)^2$. This is a simple consequence of the result of 
Hausser and Nill \cite{HN}, who constructed a quasi-Hopf algebra 
$D(A)$ (the double of $A$), isomorphic to $A\otimes A^*$ as a
vector space, such that 
$Z(\mC)=\Rep(D(A))$. \end{remark}

\begin{proof} Recall from Section 5 
that $V\to I(V)$ denotes the induction
functor. By Proposition \ref{I(1)}, we have 
\begin{eqnarray*}
{\rm FPdim}(I(V))
&=& \oplus_{X\in Irr(Z(\mC))}{\rm FPdim}(X)[I(V):X]\\
&=& \oplus_{X\in Irr(Z(\mC))}{\rm FPdim}(X)[X|_\mC:V] \\
&=& [R_{Z(\mC)}:V]=\frac{{\rm FPdim}Z(\mC)}{{\rm FPdim}\mC}
[R_\mC:V]=\frac{{\rm FPdim}Z(\mC)}{{\rm FPdim}\mC}{\rm FPdim}(V).
\end{eqnarray*}
In particular, ${\rm FPdim}(I(\mathbf 1))=\frac{{\rm
    FPdim}Z(\mC)}{{\rm FPdim}\mC}$.

On the other hand, by Proposition \ref{I(1)}, $I(\mathbf 1)=\oplus_{Y\in
  Irr(\mC)}Y\otimes Y^*$, which implies the desired statement
(since the dimension of the right hand side is ${\rm
  FPdim}(\mC)$). The proposition is
proved. 
\end{proof}

\begin{corollary}\label{fpdual} Let $\mC$ be a fusion category and $\mC^*$ be
  its dual with respect to an indecomposable module category. 
Then ${\rm FPdim}(\mC)={\rm FPdim}(\mC^*)$. 
\end{corollary}

\begin{proof} This follows from Proposition \ref{squar} and the
  fact that
$Z(\mC)$ is equivalent to $Z(\mC^*)$ \cite{O1,Mu2}.  
\end{proof} 

\begin{proposition}\label{subdivi}
Let $\mC$ be a full tensor subcategory of a fusion category
$\mD$. Then ${\rm FPdim}(\mD)/{\rm FPdim}(\mC)$ is an algebraic
integer. 
\end{proposition}

\begin{proof}
Consider the decomposition of $\mD$, as a left module category 
over $\mC$, in a direct sum of indecomposables. 
We have $\mD=\oplus_{i=1}^r \mM_r$. 
Let us consider the duals of $\mC$ and $\mD$
with respect to this module category. We have $\mD^*=\mD^{op}$,
and $\mC^*$ is a multi-fusion category $\mC^*=\oplus_{i,j=1}^r \mC^*_{ij}$, 
where 
$\mC^*_{ij}:={\rm Fun}(\mM_j,\mM_i)$. 
Note that the Frobenius-Perron dimensions 
of $\mC_{ii}^*$ are equal to ${\rm FPdim}(\mC)$
by Corollary \ref{fpdual}. 

Let $E_{ij}$ denote nonzero virtual objects of $\mC^*_{ij}$, 
which satisfy the condition \linebreak $X\otimes E_{ij}\otimes Y=
{\rm FPdim}(X){\rm FPdim}(Y)E_{ij}$, $X\in \mC^*_{ii},Y\in
\mC^*_{jj}$. Since the category $\mC_{ij}^*$ is an indecomposable
bimodule over $\mC_{ii}^*$ and $\mC_{jj}^*$, by the
Frobenius-Perron 
theorem, the object $E_{ij}$ is unique up to scaling. 
We require that the coefficients in $E_{ij}$ be positive, which 
determines it uniquely up to a positive scalar. 

To fix the scaling completely, we will require that 
$E_{ii}$, which is a multiple of the regular representation
$R_{ii}$ in
$\mC^*_{ii}$, be actually equal to ${\rm FPdim}(\mC)^{-1}R_{ii}$.  
Further, we normalize $E_{1j}$ arbitrarily, then normalize 
$E_{j1}$ by the condition $E_{j1}\otimes E_{1j}=E_{jj}$, and
finally define $E_{jk}$ as $E_{j1}\otimes E_{1k}$. 
It is easy to see that $E_{ij}$ satisfy the matrix unit relations
$E_{ij}\otimes E_{kl}=\delta_{jk}E_{il}$. 

Now let $F:\mD^*\to \mC^*$ be the natural functor, 
$R$ be the regular representation of $\mD^*$, and consider
$F(R)\in \mC^*$. We claim that $F(R)=\oplus_{i,j} b_{ij}E_{ij}$, where
$b_{ij}>0$. Indeed, $F(R)^{\otimes 2}={\rm FPdim}(\mD)F(R)$, 
so for any $i$ the virtual object $L_i:=\oplus_j b_{ij}E_{ij}$ is
an
eigen-object of the functor of right multiplication by $F(R)$.
Since the category $\oplus_j \mC_{ij}^*$ is an indecomoposable
right module over $\mC^*$, 
by the Frobenius-Perron theorem such eigen-object with positive
coefficients is unique up to scaling. Thus, for any 
$X\in \mC_{ii}^*$, $X\otimes L_i=c(X)L_i$, where $c(X)$
is a scalar. By the Frobenius-Perron theorem $c(X)={\rm
  FPdim}(X)$. Similarly, $P_j:=\oplus_i b_{ij}E_{ij}$ 
satisfies $P_j\otimes X={\rm FPdim}(X)P_j$, $X\in \mC^*_{jj}$. 
This implies the claim. 

Now let $B$ be the matrix $(b_{ij})$. Since 
$F(R)^{\otimes 2}={\rm FPdim}(\mD)F(R)$, 
and $E_{ij}$ behave like matrix units, we have $B^2={\rm
  FPdim}(\mD)B$. Thus eigenvalues of $B$ are ${\rm FPdim}(\mD)$
and zero, and $B$ is diagonalizable. 
Since $b_{ij}>0$, by the Frobenius-Perron theorem
the multiplicity of the eigenvalue ${\rm FPdim}(\mD)$ is $1$, so
the rank of $B$ is $1$. Hence $\Tr(B)={\rm FPdim}(\mD)$. 

On the other hand, let us compute $b_{ii}$. 
We have 
$$
\sum_{X\in Irr(\mD)}{\rm FPdim}(X)[F(X):\mathbf 1_{ii}]=
b_{ii}/{\rm FPdim}(\mC)
$$ 
(as $E_{ii}=R_{ii}/{\rm
  FPdim}(\mC)$). Thus, $b_{ii}/{\rm FPdim}(\mC)$ is an algebraic
integer. Summing over all $i$, we find that 
$\Tr(B)/{\rm FPdim}(\mC)={\rm FPdim}(\mD)/{\rm FPdim}(\mC)$ is an
algebraic integer, as desired. 
\end{proof}

\begin{remark} Since $b_{ii}\ge 1$, we see that 
${\rm FPdim}(\mD)/{\rm FPdim}(\mC)$ dominates the number 
of indecomposable blocks in 
$\mD$ as a left $\mC$-module category.\end{remark}

\begin{remark} Let $\mM_1=\mC$. Then $\mC_{i1}^*=\mM_i$. 
It is easy to see that $F(R)_{i1}=\sum_{X\in Irr(\mM_i)}
{\rm FPdim}(X)X$. By uniqueness of eigen-object of 
$F(R)$ in $\oplus_j\mC^*_{ji}$, we have $F(R)_{i1}\otimes F(R)_{1i}=
{\rm FPdim}(\mC)F(R)_{ii}$. Thus, 
$$
[F(R)_{ii}:\mathbf 1_{ii}]=\dim{\rm End}(F(R)_{i1},F(R)_{i1})/{\rm
FPdim}(\mC)={\rm FPdim}(\mM_i)/{\rm FPdim}(\mC).
$$ 
This implies the refinement of Theorem \ref{subdivi}
stating that ${\rm FPdim}(\mM_i)/{\rm FPdim}(\mC)$ is an
algebraic integer. \end{remark}

\begin{remark} We do not know if the global dimension 
of a fusion category is divisible by the global dimension 
of its full tensor subcategory. \end{remark}

\begin{example} The following example shows that 
Proposition \ref{subdivi} is a property of tensor categories
and not just of based rings.
 
Let $G$ be a finite group and $n$ a nonnegative integer. 
Consider the based ring $A_n(G)$ with basis $G\cup \lbrace
X\rbrace$, with multiplication given by the usual multiplication
in $G$ and the relations $gX=Xg=X,g\in G$, and $X^2=nX+\sum_{g\in
  G}g$. Such rings for special $n,G$ can arise as Grothendieck 
rings of fusion categories (for example, 
the categories of \cite{TY}, or representation 
categories of affine groups ${\rm Aff}(\mathbb F_q)$). 
We claim that if $A_n(G)$ is categorifiable, 
then the number $(n+\sqrt{n^2+4|G|})^2/4|G|$
is an algebraic integer. In particular, 
for $n\ne 0$ one has $|G|\le 2n^2$. (For $n=1$ we find 
$|G|=1,2$, so $G=1$ or $G=\mathbb Z/2\mathbb Z$, 
which was proved by T.~Chmutova even for non-rigid categories). 
Indeed, let $\mC$ be a fusion category with Grothendieck ring 
$A_n(G)$. Then the Frobenius-Perron dimensions of elements of $G$
are $1$, while the Frobenius-Perron dimension of $X$ is the 
largest root of the equation $x^2-nx-m=0$, where $m=|G|$. 
Thus $x=\frac{1}{2}(n+\sqrt{n^2+4m})$. 
By Proposition \ref{subdivi}, 
$x^2/m$ should be an algebraic integer.  
\end{example} 

\begin{proposition}\label{dimgradi}
If a fusion category $\mC$ is 
faithfully graded by a finite group $G$, 
then the Frobenius-Perron dimensions of $\mC_g$ are equal for all $g\in G$. 
In particular, $|G|$ divides ${\rm FPdim}(\mC)$. 
\end{proposition}

\begin{proof}
Let $R$ be the regular object of $\mC$. Then 
$R=\oplus R_g$, where $R_g$ is a virtual object of 
$\mC_g$. We have $R\otimes R_h={\rm FPdim}(R_h)R$.
Taking the $gh$-th component of this equation, we find 
$R_g\otimes R_h={\rm FPdim}(R_h)R_{gh}$. 
Similarly, $R_g\otimes R={\rm FPdim}(R_g)R$, which yields 
$R_g\otimes R_h={\rm FPdim}(R_g)R_{gh}$. 
Since $R_{gh}\ne 0$, we get ${\rm FPdim}(R_g)={\rm FPdim}(R_h)$, i.e.
${\rm FPdim}(R_g)={\rm FPdim}(X)/|G|$. 
\end{proof}

\subsection{Comparison of the global and Frobenius-Perron
dimensions}

\begin{proposition} \label{domi} Let $\mC$ be a fusion category. Then 
for any simple object $V$ of $\mC$ one has $|V|^2\le {\rm FPdim}(V)^2$, 
and hence $\dim(\mC)\le {\rm FPdim}(\mC)$. Hence, if $\dim(\mC)={\rm FPdim}(\mC)$ 
then $|V|^2={\rm FPdim}(V)^2$ for all simple objects $V$. 
\end{proposition}

\begin{proof} Let $X_i$ be the $i$-th simple object of $\mC$,
and $d_i=d(X_i)$ be its dimension introduced in the proof of Theorem \ref{posi}.
Then $|d_i|^2=|V|^2$, and  $d_i$ is an eigenvalue of the matrix $T_i$. 
It follows from Section 3
that the absolute values of the entries of $T_i$ are dominated by the corresponding 
entries of the matrix $N_i$ of multiplication by $X_i$.
Let $f_i$ be the Frobenius-Perron dimension of $X_i$,  
and consider a norm on the Grothendieck group $Gr_{\CC}(\mC)$ given by 
$||\sum a_iX_i||=\sum f_i|a_i|$. Let us estimate the norm of the operator $T_i$. 
We have 
$$
||T_i\mathbf a||=\sum_jf_j|\sum_k(T_i)_{jk}a_k|\le 
\sum_{j,k}f_j(N_i)_{jk}|a_k|=
f_i\sum_kf_k|a_k|=f_i||\mathbf a||
$$
Thus, $||T_i||\le f_i$. This implies 
that $|V|^2\le f_i^2$, which proves the proposition.
\end{proof}

\begin{proposition}\label{divi}
Let $\mC$ be a fusion category of global dimension $D$ and Frobenius-Perron dimension 
$\Delta$. Then $D/\Delta$ is an algebraic integer, which is $\le 1$. 
\end{proposition}

\begin{proof} By Proposition \ref{domi}, we only need to prove that 
$D/\Delta$ is an algebraic integer. 

First of all, we can assume without loss of generality that 
$\mC$ is a spherical category. Indeed, if $\mC$ is not spherical, then, 
as explained in Proposition \ref{sphpiv}, we can construct a 
spherical fusion category 
$\tilde \mC$ which projects onto $\mC$. It is easy to see that 
this category has global dimension $2D$ and Frobenius-Perron dimension $2\Delta$,
so the ratio of dimensions is the same. 

Now assume that $\mC$ is spherical. Consider the Drinfeld center $Z(\mC)$.
This is a modular category
of global dimension $D^2$, and by Proposition \ref{squar}, 
its Frobenius-Perron dimension is $\Delta^2$. Let $S=(s_{ij})$ be its S-matrix. 
It follows from Verlinde's formula (see \cite{BaKi}) that 
the commuting matrices $N_i$ of multiplication by simple objects $X_i$ of $Z(\mC)$
have common eigenvectors $\mathbf v_j$ with eigenvalues $s_{ij}/s_{0j}$. 
In particular, there exists a unique (by nondegeneracy of $S$) distinguished label 
$j=r$ such that $s_{ir}/s_{0r}=f_i$, the Frobenius-Perron dimension of $X_i$. 
Since $S$ is symmetric, and $S^2$ is the charge conjugation matrix, 
we have $\Delta^2=\sum_i f_i^2=\delta_{rr^*}/s_{0r}^2$. Thus $r=r^*$ and 
$\Delta^2=1/s_{0r}^2=D^2/d_r^2$, where $d_r$ is the categorical dimension of 
the $r$-th simple object. So $D^2/\Delta^2=d_r^2$ is an algebraic integer. 
\end{proof}

\subsection{Pseudo-unitary categories}

Let us say that a fusion category $\mC$ is {\em pseudo-unitary} if 
$\dim(\mC)={\rm FPdim}(\mC)$. This property is automatically satisfied for 
unitary categories, which occur in the theory of operator algebras,
but it fails, for example, for the Yang-Lee category. 

\begin{proposition}\label{uniqpivo} 
Any pseudo-unitary 
fusion category $\mC$ admits a unique pivotal 
(in fact, spherical) structure, 
with respect to which the categorical dimensions of all simple objects are positive, and 
coincide with their Frobenius-Perron dimensions. 
\end{proposition}

\begin{proof} Let 
$g: Id\to ****$ be an isomorphism of tensor functors, 
$a:Id\to **$ an isomorphism of additive functors, such that $a^2=g$,  
$d_i$ be the dimensions of simple 
objects $X_i$ associated to $a$ as in Section 3.1,
and $\mathbf d$ the vector with components $d_i$. Then $T_i\mathbf d=d_i\mathbf d$,
and $|d_i|=f_i$ is the Frobenius-Perron dimension of $X_i$. 
But $|(T_i)_{jk}|\le (N_i)_{jk}$. Thus, 
$$
f_if_j=|d_id_j|=|\sum (T_i)_{jk}d_k|\le \sum (N_i)_{jk}f_k=f_if_j
$$
This means that the inequality in this chain is an equality. In 
particular $(T_i)_{jk}=\pm (N_i)_{jk}$, and 
the argument of $d_id_j$ equals the argument 
of $(T_i)_{jk}d_k$ whenever $(N_i)_{jk}>0$.
This implies that whenever $X_k$ occurs in the tensor product 
$X_i\otimes X_j$, the ratio $d_i^2d_j^2/d_k^2$ is positive. 
Thus, the automorphism of the identity functor 
$\sigma$ defined by $\sigma|_{X_i}=d_i^2/|d_i|^2$ 
is a tensor automorphism. Let us twist $g$ by this automorphism. 
After this twisting, the new dimensions $d_i$ will be real. 
Thus, we can assume without loss of generality that 
$d_i$ were real from the beginning. 

It remains to twist the square root $a$ of $g$ by 
the automorphism of the identity functor $\tau$ given by 
$\tau|_{X_i}=d_i/|d_i|$. After this twisting, 
new $T_i$ is $N_i$ and new $d_k$ is $f_k$. This means 
that $\beta_i^{jk}=1$ and $a$ is a pivotal structure
with positive dimensions. It is obvious that such a structure is unique. 
We are done.  
\end{proof}

\subsection{Fusion categories with integer Frobenius-Perron dimension}
 
\begin{proposition}\label{intFP} Let $\mC$ be a fusion category, such that
  $\Delta:={\rm FPdim}(\mC)$ is an integer. Then $\mC$ is
  pseudo-unitary. In particular, the category of representations
  of a semisimple finite dimensional quasi-Hopf algebra is 
pseudo-unitary. 
\end{proposition}

\begin{proof} Let $D$ be the global dimension of $\mC$.
Let $D_1=D,D_2,...,D_N$ be all the algebraic conjugates of $D$, and 
$g_1,...,g_N$ the elements of $Gal(\bar{\mathbb Q}/{\mathbb Q})$ such that 
$D_i=g_i(D)$. Applying Propostion \ref{divi} to the category
$g_i(\mC)$, we find that $D_i/\Delta$ is an algebraic integer
$\le 1$. Then $\prod_i (D_i/\Delta)$ is an algebraic integer $\le 1$. 
But this product is a rational number. So it equals $1$ and hence all
factors equal $1$. Hence $D=\Delta$, as desired. 
\end{proof}

\begin{remark}
We note that Proposition \ref{intFP} may be regarded as a
categorical analog of the well known 
theorem of Larson and Radford \cite{LR2}, saying that in a
semisimple Hopf algebra, the square of the antipode is the
identity. More precisely, using Proposition \ref{intFP}, 
one can obtain the following proof of the Larson-Radford 
theorem (which is close to the original proof). 
Let $A$ be a finite dimensional semisimple Hopf algebra with
antipode $S$. Since $S$ has finite order, $S^2$ is semisimple,
and its eigenvalues are roots of $1$. 
Now, for every simple $A$-module $V$, $S^2$ preserves ${\rm
  End}(V)\subset A$ (as $V^{**}$ is isomorphic to $V$), and 
$|V|^2=\Tr(S^2|_{{\rm End}(V)})$, while ${\rm FPdim}(V)^2=
\dim(V)^2=\dim{\rm End}(V)$. Thus, by Proposition \ref{intFP}, 
$\Tr(S^2|_{{\rm End}(V)})=\dim{\rm End}(V)$. Since $S^2$ is
semisimple and its eigenvalues are roots of unity, 
we have $S^2=1$. \end{remark}

\begin{remark} In general, 
the Frobenius-Perron dimension of a fusion category
may not only be different from its global
dimension, but may in
fact be not conjugate to the global dimension 
by the Galois group action. The same is true for 
the numbers ${\rm FPdim}(V)^2$ and $|V|^2$ for a simple 
object $V$. An example is a tensor product of the Yang-Lee
category and its Galois conjugate. It has global dimension $5$,
but Frobenius-Perron dimension $5(3+\sqrt{5})/2$. 
Also, this category has a simple object $V$ with $|V|^2=1$ but 
${\rm FPdim}(V)^2=(7+3\sqrt{5})/2$. \end{remark}

One may wonder if a fusion category $\mC$ with integer
${\rm FPdim}(\mC)$ must have integer Frobenius-Perron dimensions 
of simple objects. Unfortunately, the answer is ``no'' (see 
\cite{TY}, where dimensions of some simple objects are square
roots of integers). However, the following weaker result is true. 
Let $\mC_{ad}$ be the full tensor subcategory of 
$\mC$ generated by objects occuring in $X\otimes X^*$, where 
$X\in \mC$ is simple (if $\mC=\Rep(G)$ for a group $G$ then 
$\mC_{ad}$ is the representation category of the adjoint group 
$G/Z(G)$ -- the quotient by the center). 

\begin{proposition}
\label{roots} 
Suppose that ${\rm FPdim}(\mC)$ is an integer.
The  Frobenius-Perron dimensions
of simple objects in $\mC_{ad}$ are integers. 
Furthermore, for any simple object $X\in \mC$, one has ${\rm
  FPdim}(X)=\sqrt{N}$, where $N$ is an integer. 
\end{proposition}

\begin{proof} The second statement follows from the first one, since
for any simple object $X\in \mC$ the object $X\otimes X^*$ lies
in $\mC_{ad}$. To prove the first statement, consider 
the object $B=\oplus_i X_i\otimes X_i^*\in \mC_{ad}$, 
where $X_i$ are the simple objects of $\mC$. 
For some $n$ the matrix of multiplication by the object 
$B^{\otimes n}$ has strictly positive entries. 
On the other hand, the vector of dimensions 
of objects of $\mC_{ad}$ is an eigenvector 
of this matrix with integer eigenvalue ${\rm FPdim}(\mC)^n$. 
Thus, by the Frobenius-Perron theorem 
its entries are rational (as $d_0=1$), and hence integer
(as they are algebraic integers).
\end{proof}

\subsection{Fusion categories of Frobenius-Perron dimension
  $p^n$}

\begin{theorem}\label{filtrat}
Let $\mC$ be a fusion category of Frobenius-Perron dimension 
$p^n$, where $n\ge 1$ is an integer, and $p$ is a prime. 
Then: 

(i) There exists a nontrivial automorphism 
$g\in {\rm Aut}_\otimes (Id_\mC)$ of order $p$, 
so that $\mC$ is $\Bbb Z/p\Bbb Z$ graded:
$\mC=\oplus_{j=0}^{p-1}\mC_j$, where $\mC_j$ 
is spanned by simple objects $X$ such that $g|_X=
e^{2\pi ij/p}$. The Frobenius-Perron dimension 
of the fusion category $\mC_0$ and module categories
 $\mC_j$ over it are equal to $p^{n-1}$. 

(ii) $\mC$ admits a filtration by fusion subcategories
$\mC\supset
\mC^{(1)}\supset...\supset \mC^{(n)}=<\mathbf 1>$, 
where $\dim \mC^{(i)}=p^{n-i}$.   
\end{theorem}

\begin{remark}
For representation categories of semisimple Hopf algebras, this 
theorem says that a semisimple Hopf algebra of order $p^n$ has a
central group-like element of order $p$, so one can take a
quotient and obtain a Hopf algebra of dimension $p^{n-1}$. 
This result is due to Masuoka \cite{Ma}. 
\end{remark}

\begin{proof} It is clear that (ii) follows from (i), so
it suffices to prove (i). 
Since the Frobenius-Perron dimension of $\mC$ 
is an integer, it is pseudo-unitary. In particular, 
it is spherical. By Proposition \ref{moddivi},
the dimensions of all objects $X\in Z(\mC)$ which occur in
$I(\mathbf 1)$, are divisors of $p^n$. 
In addition, they are all integers, since 
the image of $X$ in $\mC$ belongs to $\mC_{ad}$.  
By the class equation (Proposition \ref{cleq}),
there have to be at least $p$
such objects $X$, which are invertible. 
This means, the group of invertible objects 
of $Z(\mC)$ is nontrivial. By Proposition \ref{subdivi},
the order of this group must divide $p^{2n}$.
In particular, this group contains a cyclic subgroup 
$G$ of order $p$. The group $G$ maps 
trivially to $\mC$, so by Proposition \ref{gradi}
we get a 
faithful $G^\vee$-grading on $\mC$. 
The rest follows from Proposition \ref{dimgradi}.
The theorem is proved.   
\end{proof}

The following corollary of this result is a categorical generalization 
of (the semisimple case of) the 
well known result of Y.~Zhu \cite{Zhu} saying that a Hopf algebra 
of prime dimension is a group algebra. 

\begin{corollary}\label{zhu}
Let $\mC$ be a fusion category of Frobenius-Perron dimension $p$, where 
$p$ is a prime number. Then $\mC$ is the category of representations 
of the group $\mathbb Z/p\mathbb Z$, in which associativity is defined by a 
cohomology class in $H^3(\mathbb Z/p\mathbb Z,\CC^*)=\mathbb Z/p\mathbb Z$.  
\end{corollary}

\begin{proof}
The corollary is immediate from Theorem \ref{filtrat}. 
Indeed, this theorem shows that $\mC$ is faithfully 
graded by $\Bbb Z/p\Bbb Z$, so by Proposition \ref{dimgradi}, 
all $\mC_g$ are 1-dimensional, i.e. $\mC_g$ contains a single simple object, 
which is invertible. Thus, $\mC$ has $p$ simple objects, 
which are invertible, and we are done.  
\end{proof}

\begin{corollary} \label{quasih}
A semisimple quasi-Hopf algebra of prime dimension $p$ 
is twist equivalent to the group algebra of a cyclic group of order $p$, with
associator defined by a 3-cocycle. 
\end{corollary}

The following proposition 
classifies fusion categories of dimension $p^2$.

\begin{proposition} Let $\mC$ be a fusion category of 
Frobenius-Perron dimension $p^2$, where $p$ is a prime number. 
If $p$ is odd, then simple objects of $\mC$  are invertible, 
so $\mC$ is the category of representations of a group of order 
$p^2$ with associativity defined by a 3-cocycle. 
If $p=2$, then either $\mC$ is the representation category 
of a group of order $4$ with a 3-cocycle, or it has 3 objects 
of dimensions $1,1$ and $\sqrt{2}$, and thus belongs to the list of 
\cite{TY}. 
\end{proposition}

\begin{proof}
Assume that not all simple objects of $\mC$ are invertible. 
Our job is to show that $p=2$. 
Indeed, in this case the dimension is 4,
and we have a subcategory of dimension 2
by Theorem \ref{filtrat}. Since dimensions of all simple objects are 
square roots of integers, there is room for only one new object 
of dimension $\sqrt{2}$. All such categories are listed in \cite{TY}
(there are two of them, and they correspond 
in physics to ``Ising model''). 

Assume the contrary, i.e. that $p$ is odd. 
We have a grading of $\mC$ by $\Bbb Z/p\Bbb Z$, and for all 
$j\in \Bbb Z/p\Bbb Z$ the dimension of $\mC_j$ is $p$ 
by Proposition \ref{dimgradi}. Clearly, 
$\mC_0$ holds $p$ invertible simple objects, and no other simple objects. 
It is also clear that all invertible 
objects of $\mC$ are contained in $\mC_0$ (as by Proposition \ref{subdivi}, 
the order of the group of invertibles is a divisor of $p^2$).
Further, any object $X$ in $\mC_j$, $j\ne 0$, must be invariant
under tensoring (on either side) with invertibles from $\mC_0$, 
since otherwise the dimension of $\mC_j$ would be greater than $p$ 
(as the dimension of $X$ is greater than $1$). This means 
that $X\otimes X^*$ contains all the invertibles from $\mC_0$, and hence 
the dimension of $X$ is at least $\sqrt{p}$. Since the dimension 
of $\mC_j$ is $p$, we find that the dimension of $X$ is exactly 
$\sqrt{p}$, and $X$ is the only simple object of $\mC_j$. 
Let $j,k\in \Bbb Z/p\Bbb Z$ be nonzero and distinct
(they exist since $p>2$). Let $X_j,X_k$ be the only 
simple objects of $\mC_j,\mC_k$. Then 
$X_j\otimes X_k^*$ is a $p$-dimensional object, which 
is an integer multiple of a $\sqrt{p}$-dimensional object, which spans
$\mC_{j-k}$. Contradicton. The proposition is proved. 
\end{proof} 

\subsection{Categories with integer Frobenius-Perron dimensions
  of simple objects}

\begin{theorem}\label{quasihopf} Let $\mC$ be a fusion category.
The following two conditions are equivalent.
\begin{enumerate}
\item[(i)]  for all $X\in \mC$, the number ${\rm FPdim}(X)$ is an integer. 
\item[(ii)] $\mC$ is a representation category of a finite dimensional quasi-Hopf algebra
$A$. 
\end{enumerate}
\end{theorem} 

\begin{remark} It is shown in \cite{EG2} that $A$ is determined uniquely up to 
twisting. \end{remark}

\begin{proof} (i)$\to$ (ii).
Let $F:\mC\to {\rm Vect}$ be the additive functor that maps 
simple objects $X_i$ of $\mC$ to vector spaces $F(X_i)$ of dimension 
${\rm FPdim}(X_i)$. For any $i,j$, choose any isomorphism 
$J_{ij}: F(X_i)\otimes F(X_j)\to F(X_i\otimes X_j)$ (this is possible by 
Theorem \ref{fp}). Then by a standard reconstruction 
argument the algebra $\End_k(F)$ has a natural structure 
of a quasi-Hopf algebra (changing $J$ corresponds to twisting the quasi-Hopf algebra). 
We note that $A$ is not, in general, a Hopf algebra, since $J$ is not required to be a 
tensor structure on $F$. 

(ii)$\to$ (i). This is clear since for $\mC=\Rep(A)$, ${\rm FPdim}(X)={\rm Dim}(X)$, 
where ${\rm Dim}(X)$ denotes the dimension of the underlying vector space.  
\end{proof}

It is obvious that the class of fusion categories with integer 
Frobenius-Perron dimensions is closed under the operations of tensor product,
passing to the opposite category, and taking a fusion subcategory.
The following theorem shows that 
it is also closed under taking duals (i.e., in the language of \cite{Mu1},\cite{O} 
the property of integer Frobenius-Perron dimensions is weak Morita invariant). 

\begin{theorem}\label{dualfrob} Let $\mC$ be a fusion category with integer 
Frobenius-Perron dimensions, $\mM$ be an indecomposable module category over 
$\mC$. Then the fusion category $\mC_{\mM}^*$ has integer Frobenius-Perron dimensions. 
\end{theorem}

\begin{proof} 
We know from \cite{Mu1},\cite{O1} that $Z(\mC)$ is canonically equivalent to 
$Z(\mC_\mM^*)$. So let $V$ be an object 
of $Z(\mC)$ such that the matrix of multiplication 
by $V$ in both $\mC$ and $\mC^*_\mM$ has strictly positive entries
(it exists, since both 
$\mC$ and $\mC^*_\mM$ are fusion categories, and so the natural functors
$Z(\mC)\to \mC$, $Z(\mC^*)\to \mC^*_\mM$ are surjective).

Let $[X], [X]_*$ be the matrices of multiplication by $X\in Z(\mC)$ 
in the Grothendieck 
groups $K(\mC),K(\mC^*_\mM)$, respectively. 
For a square matrix $A$ with 
nonnegative entries, let $t(A)$ be the largest real eigenvalue of $A$.
By the assumptions of the theorem, $t([V])$ is an integer
(it is the Frobenius-Perron dimension of $V$ as an object of $\mC$). 
Since the asssignments $X\to t([X])$, $X\to t([X]_*)$, $X\in Z(\mC)$, are both 
homomorphisms of algebras, Lemma \ref{uniq} implies 
that $t([V])=t([V]_*)$. Hence, $t([V]_*)$ is an integer. 

The vector of Frobenius-Perron dimensions of simple objects 
in $\mC_\mM^*$ is the unique, up to scaling, eigenvector
of $[V]_*$ with eigenvalue $t([V]_*)$. 
Since the matrix $[V]_*$ has integer entries, the entries 
of this vector are rational numbers (as the zeroth entry is $1$). 
This means that they are integers (as they are eigenvalues of integer 
matrices and hence algebraic integers).
The theorem is proved.
\end{proof}

\begin{corollary}\label{quotcat} 
The class of fusion categories with integer Frobenius-Perron dimensions
is closed under the operation of taking a component of a quotient category. 
\end{corollary}

\begin{proof}
Let $\mC$ be a fusion category with integer Frobenius-Perron dimensions, 
$\mD$ an indecomposable multi-fusion category, and 
and $F:\mC\to \mD$ be a surjective tensor functor (so $\mD$ is a quotient of $\mC$). 
Let $\mM$ be an indecomposable module category over $\mD$. Then 
$F^*:\mD^*\to \mC^*$ is an injective tensor functor between the
corresponding dual fusion categories
($\mC^*$ is fusion since $\mM$ is indecomposable over $\mC$, by surjectivity of $F$). 
By Theorem \ref{dualfrob}, $\mC^*$ has integer Frobenius-Perron dimensions. 
Hence, so does $\mD^*$ as its full subcategory. But any component
category $\mD_{ii}$ of $\mD$ is dual 
to $\mD^*$ with respect to the $i$-th part of $\mM$ as a 
$\mD^*$-module. Thus, $\mD_{ii}$ has integer Frobenius-Perron dimensions 
by Theorem \ref{dualfrob}. 
\end{proof}

\begin{proposition}\label{frobmod1} 
Let $\mC$ be a fusion category with integer Frobenius-Perron dimensions.
Then for any indecomposable module category $\mM$ over $\mC$, the
Frobenius-Perron dimensions 
of objects in $\mM$ are integers (under a suitable normalization).
\end{proposition}

\begin{proof} Since the matrix $[V]_\mM$ from Proposition \ref{frobmod}
has integer entries, and its largest real eigenvalue ${\rm FPdim}(V)$ is an integer
and has multiplicity $1$, the 
entries of $\mathbf m$ are rational.
\end{proof}

In \cite{O2}, it is asked whether 
there are finitely many fusion categories with a given number 
of simple objects, and proved that the answer is ``yes'' 
for representation categories of Hopf algebras. Using the above
techniques, one can show that the answer is ``yes'' 
for fusion categories with integer Frobenius-Perron
dimension (which includes Hopf and quasi-Hopf algebras). 
Namely, we have

\begin{proposition}
There are finitely many fusion categories of integer
Frobenius-Perron dimension with a given number of simple
objects. In particular, there are finitely many equivalence
classes of semisimple quasi-Hopf algebras with a given number 
of matrix blocks. 
\end{proposition}

\begin{proof} Let $\mC$ be such a category with $N$ simple objects. 
Since the Frobenius-Perron dimension of $\mC$ is an integer, it is 
pseudo-unitary, in particular has a canonical spherical structure. 
Let $d_i$ be the categorical (=Frobenius-Perron)
dimensions of simple objects of $Z(\mC)$
which occur in $I(\mathbf 1)$, and 
$D$ be its global (=Frobenius-Perron) dimension.
Then $d_i$ are integers (by Proposition \ref{roots}, 
since for $X_i$ occuring in $I(\mathbf 1)$, $X_i|_\mC$ 
lies in $\mC_{ad}$ by Proposition \ref{I(1)}). 
Thus, $D/d_i=m_i$ are (usual)
integers by Proposition \ref{moddivi}, 
and $\sum_i[X_i|_\mC:1]\frac{1}{m_i}=1$. 
The number of summands $1/m_i$ here 
is $\sum [X_i|_{\mC}:\mathbf 1]$, which is at most
$N$ by Proposition \ref{numobj}. By a classical argument 
of Landau (see \cite{O2} and references therein), 
the number of such collections $m_i$ 
is finite for any given $N$. 
On the other hand, $m_0=D$ (as $X_0$ is the neutral object). Thus, $D$ 
is bounded by some function of $N$. 
But for a given $D$, 
there are only finitely many fusion rings, and hence, by Ocneanu
rigidity, finitely many fusion categories. 
\end{proof}

\subsection{Group-theoretical fusion categories}

An example of a fusion category with integer Frobenius-Perron dimensions 
is the category $\mC(G,H,\omega,\psi)$
associated to quadruples $(G,H,\omega,\psi)$, where 
$G$ is a finite group, $H$ is a subgroup, $\omega\in Z^3(G,\CC^*)$ is  a 3-cocycle,  
$\psi\in C^2(H,\CC^*)$ a 2-cochain such that $d\psi=\omega|_H$ (see \cite{O1}). 
Namely, let ${\rm Vec}_{G,\omega}$ be the category of 
finite dimensional $G$-graded vector spaces with associativity defined by $\omega$.
Let ${\rm Vec}_{G,\omega}(H)$ be the subcategory of ${\rm Vec}_{G,\omega}$ of 
objects graded by $H$. 
Consider the twisted group algebra $A=\CC_\psi[H]$. Since $d\psi=\omega|_H$, it is an  
associative algebra in this category. 
Then $\mC(G,H,\omega,\psi)$ is defined to be 
the category of $A$-bimodules in ${\rm Vec}_{G,\omega}$. 

\begin{remark} If $(\omega',\psi')$ is another pair of elements such that 
$\omega'=\omega+d\eta,\psi'=\psi+\eta|_H+d\chi$, 
$\eta\in C^2(G,\CC^*)$, $\chi\in C^1(H,\CC^*)$, then 
$\mC(G,H,\omega,\psi)$ is equivalent to $\mC(G,H,\omega',\psi')$. 
Thus, the categories $\mC(G,H,\omega,\psi)$ 
are parametrized by the set $S(G,H)$ of equivalence classes of
pairs $(\omega,\psi)$ 
under 
the above equivalence, which is a fibration over ${\rm Ker}(H^3(G,\CC^*)\to 
H^3(H,\CC^*))$ with fiber being a torsor over ${\rm Coker}(H^2(G,\CC^*)\to H^2(H,\CC^*))$. \end{remark} 

\begin{definition} A fusion category is said to be group-theoretical if it is 
of the form $\mC(G,H,\omega,\psi)$.
\end{definition}

\begin{remark} We note that there may be more than one way
to represent a given fusion category in the form
$\mC(G,H,\omega,\psi)$. In particular, $G$ is not uniquely determined. 
However, $|G|$ is uniquely determined, since it equals the global dimension 
of $\mC(G,H,\omega,\psi)$. \end{remark}

Group-theoretical fusion categories have the following simple characterization. 
Let us say that a fusion category $\mD$ is pointed if all its simple objects are invertible. 
Any pointed fusion category has the form ${\rm Vec}_{G,\omega}$.

\begin{proposition}\label{1dim}\cite{O1}. A fusion category $\mC$ is group-theoretical 
if and only if it is dual to a pointed category with respect to some 
indecomposable module category.
\end{proposition}

\begin{corollary}\label{intdim} Group-theoretical fusion categories have integer 
Frobenius-Perron dimensions of simple objects. In particular,
they are pseudo-unitary.
\end{corollary}

\begin{proof} This follows from Proposition \ref{1dim}, Theorem \ref{dualfrob}, 
and the fact that any invertible object has Frobenius-Perron dimension $1$.  
\end{proof} 

It is clear from the results of 
\cite{O1} that the class of group-theoretical categories is closed 
under tensor product, taking the opposite category, and 
taking the dual category with respect to an indecomposable module category. 
Moreover, the following result shows that 
this class is also closed under taking a subcategory or 
a component of a quotient category. 

\begin{proposition}\label{subquot} 

(i) A full fusion subcategory of a group-theoretical category is group-theoretical.

(ii) A component in a quotient category of a group-theoretical category is group-theoretical.  
\end{proposition}

\begin{proof} (i) Let $\mC\subset \mD$ be fusion categories ($\mC$ is a full 
subcategory), and let $\mD$ be group-theoretical.
Let $\mM$ be an indecomposable module category 
over $\mD$ such that $\mD^*$ is ${\rm Vec}_{G,\omega}$. Then by 
Proposition \ref{star}, $\mC^*$ is a quotient of ${\rm Vec}_{G,\omega}$. 
Let $I$ label the component categories of $\mC^*$.
 It is clear that for each $g\in G$ and $i\in I$ there exists 
a unique $g(i)\in I$ such that the functor $g\otimes$ is an equivalence
$\mC_{ir}\to \mC_{g(i)r}$ for each $r\in I$. So we have an action of $G$ on $I$. 
Let $j\in I$ and $H$ be the stabilizer of $j$ in $G$. 
For any $g\in H$, denote by $\bar g_{jj}$ the projection of $g$
to the fusion category $\mC^*_{jj}$.
Then the assignment $g\to \bar g_{jj}$ is a surjective tensor functor 
${\rm Vec}_{H,\omega}\to \mC^*_{jj}$. This functor 
must map invertible objects to invertible objects. Hence, 
all simple objects of $\mC^*_{jj}$ are invertible, so $\mC^*_{jj}$ is pointed, and hence $\mC$ 
is group-theoretical (as it is dual to $\mC^*_{jj}$ with respect to a module category).

(ii) Let $\mD$ be group-theoretical, $F:\mD\to \mC$ be surjective
($\mC$ is indecomposable), and $\mC_{ii}$ a component of $\mC$. 
We need to show that $\mC_{ii}$ is group-theoretical. 
Let $\mM$ be an indecomposable module category over $\mC$. 
It suffices to show that $\mC^*$ is group-theoretical, as $\mC^*$ is dual to $\mC_{ii}$
with respect to the i-th part of $\mM$. 
But by Proposition \ref{star}, $\mC^*$ is embedded into $\mD^*$, so it
suffices to know that $\mD^*$ is group-theoretical, which follows by duality from the fact 
that $\mD$ is group-theoretical.  
\end{proof}

\subsection{A question}

We think that the following question is interesting. 

\begin{question} \label{quest} Does there exist a finite
  dimensional semisimple Hopf algebra $H$ whose representation category 
is not group-theoretical?  
\end{question}

\begin{remark} Hopf algebras with group theoretical category of 
representations 
can be completely classified in group theoretical terms 
(see \cite{O1}). Therefore, a negative answer to 
Question \ref{quest}  
would provide a full classification of semisimple Hopf algebras. 
\end{remark}

\begin{remark} The answer to Question \ref{quest}
is ``no'' for triangular Hopf algebras, as 
follows from \cite{EG3}. However, if the answer is ``yes'' 
in general, then it is ``yes'' already for quasitriangular Hopf algebras. 
Indeed, if $\Rep(D(H))$ is group-theoretical 
then so is $\Rep(H)\otimes \Rep(H^{*cop})$ (as a dual category to
$\Rep D(H)$), from which it follows by Proposition \ref{subquot} that 
so is $\Rep(H)$. \end{remark}

\begin{remark} It is interesting that there exists a quasi-Hopf
algebra $H$ whose category of representations is not
group-theoretical. Namely, recall that in \cite{TY}, 
there is a classification of fusion categories whose
simple objects are elements of a group $G$ and 
an additional object $X=X^*$, with fusion rules 
$g\otimes X=X\otimes g=X, X\otimes X=\sum_{g\in G}g$
(and the usual multiplication rule in $G$). The result is 
that $G$ must be abelian, and for each $G$ the categories are 
parametrized by pairs $(q,r)$ where $q$ is a nondegenerate 
symmetric bilinear form on $G$ with coefficients 
in $\CC^*$, and $r$ is a choice of $\sqrt{|G|}$. 
Let us denote such category by $TY(G,q,r)$. 
It is clear that $TY(G,q,r)$ is a representation category of a
quasi-Hopf algebra (i.e. has integer Frobenius-Perron dimensions)
iff $|G|$ is a square. So let us take $G=(\mathbb Z/p)^2$, where $p>2$
is a prime, and let $q$ be given by $q((x,y),(x,y))=ax^2-by^2$, 
where $a,b\in \mathbb F_p$ and $a/b$ is a quadratic non-residue
(i.e. the quadratic form $ax^2-by^2$ is elliptic). 
Then one can check by direct inspection 
that $TY(G,q,r)$ is not group-theoretical. 
(We note that $TY(G,q,r)$ is not a representation category of a
Hopf algebra; on the other hand, if $q$ is a hyperbolic form, 
then $TY(G,q,+\sqrt{|G|})$ is a representation category of a
Hopf algebra, but is group-theoretical). \end{remark}

\subsection{Cyclotomicity of dimension functions}\label{cy}
The forgetful functor $Z(\mC)\to \mC$
induces a homomorphism of rings $F:K(Z(\mC))\to K(\mC)$. It is clear that 
the image of this map is contained in the center $Z(K(\mC))$ of the ring 
$K(\mC)$.

\begin{lemma} \label{centers}
The map $K(Z(\mC))\otimes {\Bbb Q} \to Z(K(\mC))\otimes {\Bbb Q}$ is
surjective. \end{lemma}

\begin{proof} Recall that $I: \mC \to Z(\mC)$ denotes the induction functor,
that is the left adjoint to the forgetful functor $Z(\mC)\to \mC$. We will
denote the induced map $K(\mC)\to K(Z(\mC))$ by the same letter. Let $b_i$ be 
the basis of $K(\mC)$. We have by Proposition \ref{I(1)} 
$F(I(x))=\sum_ib_ixb_i^*$ for any $x\in K(\mC)$. In particular for
$x\in Z(K(\mC))$ we have $F(I(x))=x\sum_ib_ib_i^*$. Finally note that
the operator of multiplication by the element $\sum_ib_ib_i^*\in Z(K(\mC))$ is
a self-adjoint positive definite operator $K(\mC)\to K(\mC)$ (with respect to
the usual scalar product defined by $(b_i,b_j^*)=\delta_{ij}$) and hence
invertible. The Lemma is proved.
\end{proof}
 
\begin{remark} In the special case $\mC =\Rep(H)$ for some Hopf algebra $H$
Lemma \ref{centers} was proved in \cite{KSZ} 6.3.
\end{remark}

\begin{theorem}\label{cyclotomic}
Let $\mC$ be a fusion category over $\Bbb C$ and let $L$ be an irreducible
representation of $K(\mC)$. There exists a root of unity $\xi$ such that
for any object $X$ of $\mC$ one has $\Tr([X],L)\in {\mathbb Z}[\xi]$.
\end{theorem}

\begin{proof} First of all we can assume without loss of generality that
the category $\mC$ is spherical. Indeed as explained in Proposition
\ref{sphpiv} we can construct a spherical fusion category $\tilde \mC$ which 
projects onto $\mC$; moreover the map $K(\tilde \mC)\to K(\mC)$ is surjective
and the simple objects of $\tilde \mC$ map to the simple objects of $\mC$.
 
Now let $\mC$ be spherical and let $b_i$ be the basis of $K(\mC)$. The element 
$e_L=\sum_i\Tr(b_i,L)b_i^*$ is proportional to a primitive central idempotent 
in $K(\mC)$, see e.g. \cite{Lu} 19.2 (b). By Lemma \ref{centers} there exists 
a primitive central idempotent $\tilde e_L\in K(Z(\mC))$ such that $e_L$ is
proportional to $F(\tilde e_L)$. The element $\tilde e_L$ can expressed in
terms of $S-$matrix of the category $Z(\mC)$, see e.g. \cite{BaKi}. 
By Theorem \ref{app} there exists a root of unity $\xi$ such that the entries
of $S-$matrix lie in ${\Bbb Q} (\xi)$ and thus $\tilde e_L\in K(Z(\mC))\otimes
{\Bbb Q} (\xi)$. Hence $e_L$ is proportional to some element of $K(\mC)\otimes
{\Bbb Q} (\xi)$. Now the coefficient of $b_0=1\in K(\mC)$ in $e_L$ equals to
$\dim (L)\in {\Bbb Q} (\xi)$ and is nonzero. Hence $e_L\in K(\mC)\otimes {\Bbb Q} (\xi)$
and $\Tr(b_i,L)\in {\Bbb Q} (\xi)$. Obviously the number $\Tr (b_i,L)$ is an
algebraic integer and the ring of integers in $\Bbb Q(\xi)$ is $\Bbb Z[\xi]$.
The Theorem is proved.
\end{proof}

\begin{remark} Again in the case when $\mC =\Rep(H)$ for some Hopf algebra
$H$ Theorem \ref{cyclotomic} is proved in \cite{KSZ} 6.3.
\end{remark}

\begin{corollary} Any irreducible representation of $K(\mC)$ is defined over
some cyclotomic field. In particular for any homomorphism $\phi: K(\mC)\to \CC$
we have $\phi (b_i)\in {\Bbb Q} [\xi]$ for some root of unity $\xi$.
\end{corollary}

\begin{proof} Let ${\Bbb Q}^{ab}$ be the field of all cyclotomic numbers. We 
already proved that $K(\mC)\otimes {\Bbb Q}^{ab}$ decomposes into direct sum of
simple algebras (such decomposition is the same as the decomposition of 
$1\in K(\mC)$ into the sum of primitive central idempotents). Now it is
well known that the Brauer group of ${\Bbb Q}^{ab}$ is trivial. The result follows.
\end{proof}

\begin{corollary}
Let $\mC$ be a fusion category over $\Bbb C$. 
There exists a root of unity $\xi$ such that
for any object $X$ of $\mC$ one has ${\rm FPdim}(X)\in {\mathbb Z}[\xi]$. 
\end{corollary}

\subsection{Module categories over products} 

Let $a,b$ be two algebraic integers.
We say that $a,b$ are coprime if there exist algebraic integers $p,q$ such that $pa+qb=1$. 

\begin{proposition} Let $\mC,\mD$ be nondegenerate fusion categories 
with coprime Frobenius-Perron dimensions. Then any 
indecomposable module category over $\mC\otimes \mD$ is $\mM\otimes \mN$,
where $\mM$ is an indecomposable module category over $\mC$, 
and $\mN$ over $\mD$. 
\end{proposition}

\begin{proof} Let $\mB$ be an indecomposable module category 
over $\mC\otimes \mD$. Let us restrict it to $\mC$. 
We get $\mB|_{\mC}=\oplus_i\mM_i$, where $\mM_i$ are indecomposable. 
Thus the action of $\mD$ on $\mB$ can be regarded as a tensor functor
$F: \mD\to \oplus_{i,j}{\rm Fun}_{\mC}(\mM_i,\mM_j)$. 
Let $\mD'$ be the image of this functor, i.e. 
the category generated by subquotients of $F(X)$. 
Then $\mD'$ is a multifusion category, which is indecomposable 
(as $\mB$ is indecomposable over $\mC\otimes \mD$). 
Since the functor $F: \mD\to \mD'$ is surjective, the 
Frobenius-Perron dimension $d$ of any component category of $\mD'$ 
divides ${\rm FPdim}(\mD)$. On the other hand, each component category 
of $\mD'$ is a tensor subcategory in $Fun_{\mC}(\mM_i,\mM_i)$ 
for some $i$. The latter is a dual category to $\mC$, so 
it has the same Frobenius-Perron dimension as $\mC$. 
Hence, $d$ must divide ${\rm FPdim}(\mC)$. Thus, $d=1$ by the coprimeness
condition. Hence $\mD'=A-bimod$, where $A$ is a semisimple algebra, 
and 
furthermore, all $\mM_i$ are equivalent to a single module category $\mM$.  
The functor $F:\mD\to \mD'$ defines a module category $\mN$ over $\mD$, 
and it is clear that $\mB=\mM\otimes \mN$. We are done.    
\end{proof}

{\bf Example.} 
Let $\mC_h$ be the fusion category associated to the quantum group 
$U_q(so_3)$ for $q=e^{\pi i/h}$. The Frobenius-Perron dimension 
of this category is $d(h)=\frac{h}{4\sin^2(\pi/h)}$ (this follows from 
[BaKi], (3.3.9)). The denominator of this fraction is an algebraic integer.
Thus, our result implies that if $h_i$ are pairwise coprime then 
any indecomposable module category over $\otimes_i \mC_{h_i}$ is the 
product of those over $\mC_{h_i}$. 

\section{Fusion categories in positive characteristic}

In this section we will consider 
fusion and multi-fusion categories over an algebraically closed 
field $k$ of positive characteristic $p$. 
We will show that many of the results of the previous sections 
continue to hold in this case after very minor modifications. 

\subsection{List of modifications}

Fusion and multi-fusion categories in characteristic $p$ are defined
in the same way as in charactersitic zero. Below we will describe the modifications
which need to be made to validate the statements from the
previous sections in this case. 

{\bf Section 2.} The statements of subsections 2.1, 2.2 which make sense 
generalize without changes, except the last statement of Theorem
\ref{posi}, which claims 
that the global dimension of a fusion category is nonzero. 
This is in fact false 
in positive characteristic. For example, if $G$ is a finite group 
whose order is divisible by $p$, then the category of modules over 
the function algebra ${\rm Fun}(G,k)$ is a fusion category of global dimension $|G|=0$. 
This is the main difference between the zero and positive characteristic case. 
To deal with this difference, we introduce the following definition. 

\begin{definition}
A fusion category over $k$ is said to be nondegenerate if 
its global dimension is nonzero. An indecomposable multi-fusion category is 
nondegenerate if at least one of its component categories is nondegenerate. 
A multi-fusion category is nondegenerate if all its indecomposable parts are 
nondegenerate. A regular semisimple 
weak Hopf algebra over $k$ is nondegenerate if its category of
representations is a nondegenrate multi-fusion category.   
\end{definition}    

Theorems \ref{dualcat}, \ref{m1m2}, 
\ref{multi-fus} remain valid if the category $\mC$ is nondegenerate. 
Proposition \ref{samedim} remains true without any additional conditions. 

In Section 2.5 (and at other places throughout), the algebra $R$ has 
to be chosen in such a way that 
its block sizes are not divisible by $p$ (since only such algebras 
admit a symmetric separability idempotent). Theorem \ref{Schlah} is valid
for such $R$. Corollary \ref{Schlah1} and Theorem \ref{traceofS2} also remain valid. 
Theorem \ref{sscoss} is valid for nondegenerate weak Hopf algebras. 

Theorems \ref{van}, \ref{ocnrig}, \ref{multocnrig}, \ref{ocnrigfun},
and Corollary \ref{ocnrigmod}  remain valid for nondegenerate categories. 
More precisely, in Theorems \ref{van}, 
\ref{ocnrigfun} one only needs nondegeneracy of the category $\mC$.     

\begin{remark} We note that the nondegeneracy condition is essential.
For example, if $A$ is the Hopf algebra  
${\rm Fun}(G,k)$, where the order of $G$ is divisible by $p$, 
then Theorem \ref{dualcat} is false for $\Rep(A)$, and Theorem 
\ref{sscoss} is false for $A$. Similarly, Theorem \ref{van} 
is false for $\Rep(A)$ and the identity functor, since the
corresponding cohomology is $H^i(G,k)$, which may be nonzero
if the order of $G$ is divisible by $p$. However, 
we do not know a counterexample to Theorem \ref{ocnrig}
in positive characteristic. \end{remark} 
 
{\bf Sections 3-7.} While Section 3 makes no sense in positive
characteristic, the material of Section 4 generalizes without
changes, as long as one considers regular weak Hopf algebras. 
The only change in Section 5 is that one should always choose
nonzero block sizes for the algebra $R$
(except subsections 5.8,5.9 where the characteristic should be zero). 
Theorems \ref{acyclic} 
and \ref{rigweak} are valid for nondegenerate weak Hopf algebras.

{\bf Section 8.} The results of Section 8.1 generalize without
changes, since they have nothing to do with the characteristic of
the ground field. So do the results of Section 8.2, except 
Proposition \ref{squar}, Corollary \ref{fpdual}, and 
Proposition \ref{subdivi}, which are valid for nondegenerate categories.
Finally, let us point out that Theorem \ref{quasihopf} is valid 
in positive characteristic. 

{\bf Example.} A semisimple Hopf algebra $H$ is nondegenerate if and
only if it is cosemisimple. Indeed, if $H$ is nondegenerate, then 
so is $\Rep(H)$, and hence $H^*$ is semisimple, since 
$\Rep(H^{*op})$ is a dual category to $\Rep(H)$. Conversely, 
if $H$ is cosemisimple, then $\Tr(S^2)\ne 0$ (\cite{LR1}), so 
the global dimension of $\Rep(H)$ is nonzero, and $H$ is nondegenerate. 

\subsection{Lifting theorems}

In this subsection we will show that a nondegenerate multi-fusion
category over a field of characteristic $p$, and a tensor 
functor between such categories, can be lifted to
characteristic zero. It is an analog 
of the result of \cite{EG1}, where this was shown for semisimple
cosemisimple Hopf algebras and morphisms between them. 

First of all, we need to define fusion categories over any commutative ring. 
To do this, recall that a fusion category over a field $k$ can be regarded 
as a collection of finite dimensional vector spaces
$H_{ij}^k={\rm Hom}(X_i\otimes X_j,X_k)$
(multiplicity spaces), together with a collection of linear maps 
between tensor products of these spaces
(associativity morphism,
evaluation, coevaluation morphisms), satisfying some equations
(axioms of a rigid tensor category). 
Now let $R$ be a commutative ring with unit, and define 
a fusion category over $R$ to be a collection of free 
finite rank $R$-modules $H_{ij}^k$ together with a collection of
module homomorphisms satisfying the same equations.
 
Tensor functors between fusion categories over $k$ can be
defined in similar terms, as collections of linear maps
satisifying algebraic equations; this allows one to define 
tensor functors over $R$ in an obvious way.    

Now let $W(k)$ be the ring of Witt vectors of $k$.
Let $I$ be the maximal ideal in $W(k)$ generated by $p$. 
If $\mC$ is a fusion category over 
$k$, then a lifting of $\mC$ to $W(k)$ is a fusion category 
$\tilde \mC$ over $W(k)$, 
together with an equivalence $\tilde \mC/I\to \mC$,
where $\mC/I$ is the reduction of $\mC$ modulo $p$. 
In a similar way one defines a lifting of a module category 
and, more generally, a lifting of a tensor functor. 

\begin{theorem}\label{lifting}
A nondegenerate multi-fusion category $\mC$ over $k$ 
admits a unique lifting to 
$W(k)$. Moreover, a tensor functor $F$ between 
nondegenerate fusion categories $\mC,\mC'$ 
over $k$ admits a unique lifting to $W(k)$.
In particular, a module category over 
a nondegenerate multi-fusion category $\mC$ has
a unique lifting to $W(k)$.  
\end{theorem}

\begin{proof}
The theorem follows from the facts 
that $H^3(\mC)=0,H^4(\mC)=0$, $H^2_F(\mC)=0$, 
$H^3_F(\mC)=0$, proved in Section 7, in the same way as
in \cite{EG1}.  
\end{proof}

\begin{corollary}\label{brsym}
A nondegenerate braided (symmetric) fusion category over $k$ 
is uniquely liftable to a braided (resp. symmetric) 
fusion category over $W(k)$. 
\end{corollary}

\begin{proof} A braiding on $\mC$ is the same as a splitting 
$\mC\to Z(\mC)$ of the tensor functor $Z(\mC)\to \mC$. 
By Theorem \ref{lifting}, such splitting is uniquely liftable. 
Hence, a braiding on $\mC$ is uniquely liftable. 

Now let us prove the result in the symmetric case. 
A braiding gives rise to a categorical equivalence 
$B: \mC\to \mC^{op}$. A braiding is symmetric iff 
the composition of two categorical equivalences $B$ and $B^{21}$ 
is the identity. Thus, Corollary follows from 
Theorem \ref{lifting}. 
\end{proof}

\begin{remark} Corollary \ref{brsym} is false for degenerate 
categories. For example, let $\tilde \mC$ be the fusion category 
corresponding to the quantum group $U_q(sl_2)$ at $q=e^{\pi
  i/p}$, where $p$ is an odd prime. This category is 
known to have semisimple reduction to $\bar{\mathbb F}_p$, 
which is symmetric. However, this reduction does not lift to a symmetic
category, since it has non-integer Frobenius-Perron dimensions. 
This example also shows that there exist semisimple and
cosemisimple triangular weak Hopf algebras over $k$ which do not
come from group algebras. For usual Hopf algebras, 
this is impossible (see \cite{EG3}). \end{remark}

\subsection{Faithfulness of lifting}

In this section we will show that the lifting procedure
of nondegenerate fusion (braided, symmetric) categories is faithful, i.e. 
if two such categories over $k$ have equivalent liftings, then they
are equivalent. We will also show that if two tensor functors 
from a nondegenerate fusion category over $k$ 
to another fusion category have isomorphic liftings, 
then they are isomorphic. 

Let $R$ be a complete local ring, $\Bbb F$ the fraction field of
$R$, $I\subset R$ the maximal ideal, and $k=R/I$ the
residue field. 

\begin{theorem}\label{faithlif}
(i) Let $\mC_1,\mC_2$ be fusion categories defined
over $R$, such that the fusion categories 
$\mC_1/I$, $\mC_2/I$ are nondegenerate, and the categories 
$\mC_1\otimes_R\Bbb F,\mC_2\otimes_R\Bbb F$
are tensor equivalent. Then $\mC_1,\mC_2$  
(and hence $\mC_1/I,\mC_2/I$) are tensor equivalent. 

(ii) The same statement is valid for braided and symmetric fusion
categories. 
\end{theorem} 

\begin{proof} (i) Pick 
a tensor equivalence $E: \mC_1\otimes_R \Bbb F\to
\mC_2\otimes_R\Bbb F$. 
Now pick an equivalence of $\mC_1$, $\mC_2$ as {\it additive}
categories, which is the same as $E$ at the level of isomorphism
classes of objects. By doing this,
we may assume without loss of generality that $\mC_1$ and $\mC_2$ are 
really one category $(\mC,\otimes)$ over $R$ with two different
associativity isomorphisms $\Phi_1,\Phi_2$
and evaluation maps
${\rm ev}_1,{\rm ev}_2$ (the coevaluation maps
of simple objects of both categories
may be identified with any fixed collection of such maps 
by renormalization). 

We know that there exists an invertible element $J$ of 
$C^2(\mC)\otimes_R \Bbb F$ such that $(\Phi_1,{\rm ev}_1)^J=
(\Phi_2,{\rm ev}_2)$. 
(where ${}^J$ denotes the twisting by $J$). Our
job is to show that such $J$ can be found in $C^2(\mC)$, 
without localization to $\Bbb F$.

Consider the affine group scheme $T$ (defined over $R$)
of all invertible elements of $C^2(\mC)$. Let $T_0$ be the
subscheme of $T$ consisting of ``trivial'' twists, 
i.e. those of the form $\Delta(x)(x^{-1}\otimes x^{-1})$,
where $x$ is an invertible element of $C^1(\mC)$ 
(i.e. the center of $\mC$). It is clear that $T_0$ is central in
$T$. 

{\bf Remark.} Since the algebra $C^2(\mC)$ is a direct sum of
matrix algebras (split over $R$), $T$ is a direct product of general linear
groups. 

Let $Z$ be the affine scheme of all invertible associators 
$\Phi$ and collections of evaluation morphisms ${\rm ev}$
satisfying (together with the fixed coevaluation morphisms) 
the axioms for a rigid tensor category
(i.e. the pentagon relation, the unit object axiom, and the duality axioms).
The group scheme $T$ acts on $Z$ by twisting
(and then renormalizing the coevaluations to be the fixed ones): 
$(\Phi,{\rm ev})\to (\Phi,
{\rm ev})^J$. 
The restriction of this action to $T_0$ is trivial, so 
we get an action of $T/T_0$ on $Z$. Consider the regular mapping 
$\phi: T/T_0\to Z$, given by $J\to (\Phi_1,{\rm ev_1})^J$. 

\begin{lemma} 
\label{fin} The mapping $\phi$ is finite.
\end{lemma}

\begin{proof}
 By Nakayama's lemma, it is sufficient to check the finiteness 
of the reduction $\phi|_k$ of $\phi$ modulo $I$. 

First of all, we have 
$H^2(\mC_1/I)=0$, which implies that the stabilizer 
of $(\Phi_1,{\rm ev}_1)$ in $T/T_0$ is a finite group. Thus the map 
$\phi|_k: T/T_0\to \phi|_k(T/T_0)$ is a finite covering.

It remains to show that $\phi|_k(T/T_0)$ is closed in $Z$. 
To do this, consider its closure $Z'$. Points of $Z'$ 
represent fusion categories over $k$, all of which 
are nondegenerate Indeed, the global dimension is a regular 
(in fact, locally constant)
function on $Z$, and it is nonzero on $\mC_1/I$, so it must be nonzero 
for all points of $Z'$. But for a nondegenerate fusion category 
$\mD$ over $k$ we have  $H^3(\mD)=0$, which implies that 
the action of $T/T_0$ on $Z'$ is locally transitive. 
This means that $Z'$ 
is in fact a single orbit of $T/T_0$, i.e. $\phi|_k(T/T_0)$ is closed. 
The lemma is proved. 
\end{proof}

Thus we have a homomorphism $\phi^*: R[Z]\to R[T/T_0]$, which
makes $R[T/T_0]$ into a finitely generated $R[Z]$-module. 
 
Now let $J\in (T/T_0)(\Bbb F)$ be such that $(\Phi_1,{\rm ev}_1)^J=
(\Phi_2,{\rm ev}_2)$. 
Thus, $J$ can be regarded as a homomorphism 
$J: R[T/T_0]\to \Bbb F$. This homomorphism 
satisfies the condition $J(af)=\gamma(a)J(f)$, 
$f\in R[T/T_0]$, $a\in R[Z]$. 
where $\gamma: R[Z]\to R$ corresponds to the point 
 $(\Phi_2,{\rm ev}_2)$. Now, by Lemma \ref{fin}, 
for any $f\in R[T/T_0]$ there exist $b_1,...,b_n\in R[Z]$ such
that $f^n+b_1f^{n-1}+...+b_n=0$. But then 
$J(f)^n+\gamma(b_1)J(f)^{n-1}+...+\gamma(b_n)=0$. 
Now, since $\gamma(b_i)$ are integers, we see that $J(f)$ is an
integer, so $J$ in fact belongs to $(T/T_0)(R)$. 
But it is easy to see that the map $T(R)\to (T/T_0)(R)$ is
surjective. Hence, there exists $J\in T(R)$ such that
$(\Phi_1,{\rm ev}_1)^J=(\Phi_2,{\rm ev}_2)$,
as desired. The theorem is proved.

(ii) The proof is the same, since, as we showed, 
the integrality of $J\in T/T_0$ 
follows already from the equation $(\Phi_1,{\rm ev}_1)^J=(\Phi_2,{\rm ev}_2)$.
\end{proof} 

\begin{theorem}\label{faithfun} 
Let $\mC,\mC'$ be fusion categories defined
over $R$, such that the fusion category 
$\mC/I$ is nondegenerate.
Let $F_1,F_2: \mC\to \mC'$ 
be tensor functors, which become isomorphic after tensoring with
$\Bbb F$. Then $F_1,F_2$ 
(and hence their reductions modulo $I$) are isomorphic. 
\end{theorem}

\begin{proof} The proof is analogous to the proof of Theorem \ref{faithlif}.
As in Theorem \ref{faithlif}, we may assume that $F_1=F_2=F$ 
as additive functors, with two different tensor structures
$J_1,J_2$. Let $G$ be the group scheme of invertible elements of
$C^1_F(\mC)$. Let $S$ be the scheme of tensor structures 
on $F$. The group scheme $G$ naturally acts 
on $S$ by ``gauge transformations'', $J\to J^g$. We know there exists 
$g\in G(\Bbb F)$ such that $J_1^g=J_2$, and need to show that
this $g$ is in fact in $G(R)$. 

Using the fact that 
$H^1_F(\mC)=H^2_F(\mC)=0$, one finds that the action map 
$\chi:G\to S$ given by $g\to J_1^g$ is finite. Now the proof 
that $g\in G(R)$ runs analogously to the proof that $J\in
(T/T_0)(R)$ in Theorem \ref{faithlif}. 
Theorem \ref{faithfun} is proved. 
\end{proof}

Now let $K$ be the algebraic closure of the field of quotients of
$W(k)$. 

\begin{corollary} \label{faithlif1}
(i) If two nondegenerate 
fusion (braided, symmetric) categories over $k$ have equivalent
liftings
to $K$, then they
are equivalent. 

(ii) If two tensor functors 
between nondegenerate fusion categories over $k$ have isomorphic
liftings to $K$, 
then they are isomorphic. 
\end{corollary}

\begin{proof} (i) follows easily from Theorem \ref{faithlif}, (ii) from
  Theorem \ref{faithfun}.
\end{proof}

To give the next corollary, recall that in \cite{EG1}
it was proved that a semisimple cosemisimple Hopf algebra 
over $k$ has a unique lifting to $K$. 
In fact, this is a special case of the lifting theory developed in this
paper, since a semisimple cosemisimple Hopf algebra
is nothing but a pair $(\mC,F)$, where $\mC$ is a nondegenerate
fusion category over $k$, and $F$ a fiber functor on $\mC$. 

\begin{corollary} \label{faithhopf}
If two semisimple cosemisimple Hopf algebras over $k$ 
have isomorphic liftings to $K$ then they are isomorphic. 
\end{corollary}

\begin{proof} We can regard these two Hopf algebras as pairs
  $(\mC_1,F_1)$ and $(\mC_2,F_2)$. 
By theorem \ref{faithlif1}(i), we can assume that $\mC_1=\mC_2$; 
then by theorem \ref{faithlif1}(ii), $F_1$ is isomorphic to
$F_2$. The Corollary is proved. 
\end{proof}

The same theorem applies to quasitriangular, triangular Hopf
algebras. The proof is parallel. 

\vskip .1in

{\bf Warning.} (Added in April 2017.) Unfortunately, the proof of Lemma \ref{fin} (and hence of the main results on faithfulness of the lifting given in 
this subsection) is incomplete; namely, the Nakayama lemma is not sufficient here, and one needs to use the fact that $T$ is a reductive group, see \cite{Et}. 
A complete proof of Lemma \ref{fin}, as well as alternative proofs of the results of this subsection are given in \cite{Et}, Section 4. This paper also gives an affirmative answer to the question in Remark \ref{nondq} below.  

\subsection{Some applications of lifting}

\begin{corollary}\label{fpnonzero}
Let $\mC$ be a nondegenerate fusion category 
over $k$. Then its Frobenius-Perron dimension $\Delta$
is not divisible by $p$. 
\end{corollary}

\begin{remark} We do not know if the converse statement is true. 
For Hopf algebras the converse statement is the well known
conjecture that if $H$ is a semisimple Hopf algebra whose
dimension is not divisible by $p$ then $H$ is cosemisimple. \end{remark}

\begin{proof} Assume that $\Delta$ is divisible by $p$. 
Let $\tilde\mC$ be the lifting of $\mC$, and $\hat
\mC=\tilde\mC\otimes_{W(k)}K$, where $K$ is the algebraic closure
of the field of fractions of $W(k)$. Then by Theorem 
\ref{divi}, the global dimension $D$ 
of $\hat\mC$ is divisible by $\Delta$, hence by $p$. 
So the global dimension of $\mC$ is zero. Contradiction. 
\end{proof}

\begin{corollary}\label{fpsubcat}
Let $\mD$ be a nondegenerate fusion category over $k$, 
and $\mC$ be a full subcategory of $\mD$ of integer
Frobenius-Perron dimension $N$. Then $\mC$ is nondegenerate.  
\end{corollary}

\begin{remark}\label{nondq} We do not know if this result is valid without
assuming that $N$ is an integer. \end{remark}

\begin{proof} Consider the lifting $\tilde \mD$ of $\mD$. 
This lifting necessarily contains a lifting 
$\tilde\mC$ of $\mC$. Consider the corresponding categories 
$\hat\mC$, $\hat\mD$ over the algebraically closed field $K$ of
zero characteristic. Let $D$ be the global dimension of
$\hat\mD$, $\Delta$ its Frobenius-Perron dimension.
Then $D$ is divisible by $\Delta$ by Proposition \ref{divi}.
Also, the Frobenius-Perron dimension of $\hat\mC$ is $N$ (an
integer), hence it is equal to the global dimension of $\hat\mC$ 
by Proposition \ref{intFP}. But by Theorem \ref{subdivi}, 
$N$ divides $\Delta$, and hence $D$. But $D$ is not divisible by
$p$, hence $N$ is not divisible by $p$. Thus, 
$\mC$ is nondegenerate, as desired. 
\end{proof} 

\section{Appendix: Galois properties of S-matrix}
The remarkable result due to J.~de~Boere, J.~Goeree, A.~Coste and T.~Gannon
states that the entries of the S-matrix of a semisimple modular
category lie in a cyclotomic field, see \cite{dBG, CG}.
This result is used in Section 8. For reader's
convenience in this Appendix we reproduce a proof of this result.

\begin{theorem}\label{app} Let $S=(s_{ij})_{i,j\in I}$ be the S-matrix of a
modular fusion category $\mC$. There exists a root of unity $\xi$ such that
$s_{ij}\in {\mathbb Q}(\xi)$.
\end{theorem}

\begin{proof}
Let $\{ X_i\}_{i\in I}$ be the representatives of isomorphism classes of
simple objects of $\mC$; let $0\in I$ be such that $X_0$ is the unit object
of $\mC$ and the involution $i\mapsto i^*$ of $I$ be defined by $X_i^*\cong
X_{i^*}$.

 By the definition of modularity, any homomorphism $f: K(\mC)\to {\mathbb C}$ is of
the form $f([X_i])=s_{ij}/s_{0j}$ for some well defined $j\in I$. Hence for
any automorphism $g$ of ${\mathbb C}$ one has $g(s_{ij}/s_{0j})=
s_{ig(j)}/s_{0g(j)}$ for a well defined action of $g$ on $I$.

Now we are going to use the following identities (see \cite{BaKi}):

(i) $\sum_js_{ij}s_{jk}=\delta_{ik^*}$;

(ii) $s_{ij}=s_{ji}$;

(iii) $s_{0i^*}=s_{0i}\ne 0$.

Thus, $\sum_js_{ij}s_{ji^*}=1$ and hence $(1/s_{0i})^2=\sum_j(s_{ji}/s_{0i})
(s_{ji^*}/s_{0i^*})$. Applying the automorphism $g$ to this equation we get
\begin{equation}
g(\frac{1}{s_{0i}^2})=g(\sum_j\frac{s_{ji}}{s_{0i}}\frac{s_{ji^*}}{s_{0i^*}})=
\sum_j\frac{s_{jg(i)}}{s_{0g(i)}}\frac{s_{jg(i^*)}}{s_{0g(i^*)}}=
\frac{\delta_{g(i)g(i^*)^*}}{s_{0g(i)}s_{0g(i^*)}}.
\end{equation}
It follows that $g(i^*)=g(i)^*$ and $g((s_{0i})^2)=(s_{0g(i)})^2$. Hence
$$g((s_{ij})^2)=g((s_{ij}/s_{0j})^2\cdot (s_{0j})^2)=(s_{ig(j)})^2.$$ Thus
$g(s_{ij})=\pm s_{ig(j)}$. Moreover the sign $\epsilon_g(i)=\pm 1$ such that
$g(s_{0i})=\epsilon_g(i)s_{0g(i)}$ is well defined since $s_{0i}\ne 0$, and
$g(s_{ij})=g((s_{ij}/s_{0j})s_{0j})=\epsilon_g(j)s_{ig(j)}=\epsilon_g(i)
s_{g(i)j}$. In particular, the extension $L$ of ${\mathbb Q}$ generated by all
entries $s_{ij}$ is finite and normal, that is Galois extension.

Now let $h$ be another automorphism of ${\mathbb C}$. We have
$$gh(s_{ij})=g(\epsilon_h(j)s_{ih(j)})=\epsilon_g(i)\epsilon_h(j)s_{g(i)h(j)}$$and
$$hg(s_{ij})=h(\epsilon_g(i)s_{g(i)j})=\epsilon_h(j)\epsilon_g(i)s_{g(i)h(j)}=
gh(s_{ij})$$
and the Galois group of $L$ over ${\mathbb Q}$ is abelian. Now the
Kronecker-Weber theorem (see e.g. \cite{Ca}) implies the result.
\end{proof}




\bibliographystyle{ams-alpha}
  
\end{document}